\documentclass[final]{siamart220329}

\RequirePackage{amsmath, amsfonts, amssymb, color, hyperref, enumerate, mathtools}

\newtheorem*{inftheorem}{``Theorem''}
\newsiamremark{defn}{Definition}
\newsiamremark{condition}{Condition}
\newsiamremark{assumption}{Assumption}
\newsiamremark{remark}{Remark}

\newcommand{\ud}{\,\mathrm{d}}
\newcommand{\R}{\mathbb{R}}
\newcommand{\Rd}{\mathbb{R}^d}
\newcommand{\G}{\mathcal{G}}

\newcommand{\ind}{\mathbf{1}}
\newcommand{\N}{\mathcal{N}}

\RequirePackage{pgfplots,neuralnetwork}
\pgfplotsset{compat=1.16}
\DeclarePairedDelimiter{\braces}{\{}{\}}

\begin{document}
	\title{Deep Gaussian Process Priors for Bayesian Inference in Nonlinear Inverse Problems\thanks{Submitted to SIAM JUQ. \funding{Supported by the EPSRC Programme Grant on the Mathematics of Deep Learning under the project EP/V026259/1.}}}
	
	\author{
		Kweku Abraham\thanks{University of Cambridge, Statistical Laboratory, Wilberforce Road, Cambridge CB3 0WB, UK. \email{lkwa2@cam.ac.uk} and \email{and30@cam.ac.uk}.} \and Neil Deo\footnotemark[2] }
	\maketitle

\headers{Deep GP priors for nonlinear inverse problems}{K.~Abraham and N.~Deo}	
	\begin{abstract}
		We study the use of a deep Gaussian process (DGP) prior in a general nonlinear inverse problem satisfying certain regularity conditions. We prove that when the data arises from a true parameter $\theta^*$ with a compositional structure, the posterior induced by the DGP prior concentrates around $\theta^*$ as the number of observations increases. The DGP prior accounts for the unknown compositional structure through the use of a hierarchical structure prior. As examples, we show that our results apply to Darcy’s problem of recovering the scalar diffusivity from a steady-state heat equation and the problem of determining the attenuation potential in a steady-state Schr\"{o}dinger equation. 
		We further provide a lower bound, proving in Darcy’s problem that typical Gaussian priors based on Whittle-Mat\'{e}rn processes (which ignore compositional structure) contract at a polynomially slower rate than the DGP prior for certain diffusivities arising from a generalised additive model. 
	\end{abstract}
	
			

	\tableofcontents
	
	
	\section{Introduction}

	Deep learning now provides state-of-the-art empirical performance in a wide range of complex tasks: image classification, speech recognition and medical imaging among others. Yet despite far-reaching empirical success, the theoretical performance of deep learning methods is not well understood. Recently, some progress has been made in obtaining statistical guarantees for deep neural networks in a nonparametric regression model in \cite{schmidt-hieberNonparametricRegressionUsing2020}, where it was shown that suitably calibrated networks achieved fast convergence rates when the signal has a compositional form.
	
	We instead consider a Bayesian deep learning method: that of \emph{deep Gaussian processes} (DGPs), used as a Bayesian prior. Gaussian process priors are some of the most widely used priors in Bayesian nonparametrics and in many instances offer optimal performance \cite[Chapter 11]{ghosalFundamentalsNonparametricBayesian2017}. Deep Gaussian processes, introduced in \cite{damianouDeepGaussianProcesses2013}, are formed by suitably iterating Gaussian processes, for example by composition. The resulting DGP can then have highly non-stationary behaviour even if the underlying Gaussian processes have smooth, stationary covariance kernels: see \cite{svendsenDeepGaussianProcesses2020,muirDeepGaussianProcess2023} for some applications to biogeophysical models and seismology. Moreover, the posterior distribution induced by the DGP prior (see \eqref{Eq: Bayes' formula} below) provides a method for uncertainty quantification, a typical benefit of Bayesian procedures. In \cite{finocchioPosteriorContractionDeep2023}, it was shown that a DGP prior achieves fast convergence rates in a nonparametric regression model when the signal has a compositional structure.
	In contrast, Gaussian priors model compositional functions poorly: in \cite{giordanoInabilityGaussianProcess2022}, it was shown in a `direct' regression problem with white noise that if $f$ arises from a \emph{generalised additive model} of the form
	\begin{equation}\label{Eq: intro generalised additive model}
		f(x) = F\left(g_1(x_1) + \ldots + g_d(x_d)\right), \quad x\in\mathcal{O},\, F,g_1,\ldots,g_d:\R\to\R
	\end{equation}
	where $F,g_1,\ldots,g_d$ are unknown, then \emph{any} mean-zero Gaussian process prior achieves a suboptimal rate. For Gaussian priors based on a random wavelet expansion, there is even a severe curse of dimensionality, in the sense that the contraction rate becomes arbitrarily slow as $d\to\infty$.
	
	Since the influential work of Andrew Stuart \cite{stuartInverseProblemsBayesian2010}, Bayesian methods have been particularly popular for solving inverse problems arising from partial differential equations (PDEs). A prototypical example is \emph{Darcy's problem}, where one seeks to recover the non-negative diffusivity $f$ from observing the solution $u$ to the PDE 
	\begin{equation}\label{Eq: intro Darcy's problem}
		\begin{aligned}
			\nabla\cdot(f\nabla u) &= g \quad \text{ on }\mathcal{O}, \\
			u&=0 \quad \text{ on }\partial\mathcal{O},
		\end{aligned}
	\end{equation}
	 where it is assumed that the source term $g$ is known. This problem has applications to subsurface hydrology, with $f$ describing the permeability of the medium through which groundwater is flowing: see \cite{yehReviewParameterIdentification1986,stuartInverseProblemsBayesian2010}. 	
	 
	
	When $f$ and $g$ are positive and sufficiently regular, equation \eqref{Eq: intro Darcy's problem} has a unique solution $u=u_f$. Let $G$ denote the solution map $f\mapsto u_f$; we consider observations $D_n:=(Y_i,X_i)_{i=1}^n$ of the form
	\begin{equation}\label{Eq: intro observation scheme}
		Y_i = G(f)(X_i) + \epsilon_i, \quad 1\leq i\leq n,
	\end{equation}
	where $X_i\overset{\mathrm{i.i.d.}}{\sim}\mathrm{Uniform}(\mathcal{O})$ and $\epsilon_i\overset{\mathrm{i.i.d.}}{\sim} N(0,1)$ independently of the $X_i$. Write $P_f$ for the law of $D_n$ under \eqref{Eq: intro observation scheme}, with associated expectation operator $E_f$.
	
	The map $f\mapsto G(f)$ is nonlinear, and so the negative log-likelihood function arising from \eqref{Eq: intro observation scheme} is possibly non-convex in $f$; as a result, optimisation-based methods such as maximum likelihood estimation or Tikhonov regularisation cannot be reliably implemented. Sampling from the Bayesian posterior, which can be done using Markov chain Monte Carlo methods (see \cite{cotterMCMCMethodsFunctions2013,hairerSpectralGapsMetropolis2014,beskosGeometricMCMCInfinitedimensional2017,nicklPolynomialtimeComputationHighdimensional2022,bohrLogconcaveApproximationsHighdimensional2023}), can avoid these shortcomings. Moreover, since $G$ arises from an elliptic PDE, it has regularity properties which can be leveraged to obtain frequentist guarantees stating that when the data arises from some fixed $f^*$, the posterior concentrates around $f^*$ in the large sample limit. This is usually expressed by a \emph{posterior contraction rate} for a suitable prior $\Pi$, which is a sequence $r_n\to0$ such that when the data $D_n$ arise from the parameter $f^*$ in \eqref{Eq: intro observation scheme},
	$$ E_{f^*} \Pi\left( f: \|f-f^*\|_{L^2} > r_n \mid D_n \right) \to 0$$
	as $n\to\infty$. One desires such a guarantee to hold uniformly over all parameters $f^*$ indexing the statistical model.
	
	Posterior contraction rates have been obtained for a variety of PDE-constrained nonlinear inverse problems in  \cite{monardConsistentInversionNoisy2021,abrahamStatisticalCalderonProblems2019,nicklBernsteinMisesTheorems2020a,giordanoConsistencyBayesianInference2020,kekkonenConsistencyBayesianInference2022}, mostly for priors based on Gaussian processes. We take the approach of \cite{monardStatisticalGuaranteesBayesian2021,nicklBayesianNonlinearStatistical2023} and study a general forward map $G$ satisfying regularity conditions; this framework encompasses both Darcy's problem \eqref{Eq: Darcy's problem} and the problem of identifying the potential from a steady-state Schr\"{o}dinger equation studied in \cite{nicklBernsteinMisesTheorems2020a} (see Section \ref{Subsection: PDE examples} below), among others. 
	
	
	The motivation for this article is to weave together these two strands of research: that is, to obtain theoretical guarantees for a deep Gaussian process prior in a nonlinear inverse problem. We show that in a general elliptic PDE inverse problem satisfying certain regularity conditions, of which Darcy's problem \eqref{Eq: intro Darcy's problem} is an instance, the DGP prior provides a method for consistent reconstruction with polynomial convergence rates. Moreover, we show that it outperforms certain Gaussian process priors (by a polynomial factor) when the true parameter arises from a generalised additive model \eqref{Eq: intro generalised additive model}. A key message of the paper is summarised in the following informal theorem.
	\begin{inftheorem}
		Consider Darcy's problem \eqref{Eq: intro Darcy's problem} with data arising from the observation model \eqref{Eq: intro observation scheme}. Let $\Pi$ be the DGP prior of \eqref{Eq: deep GP prior defn}. Let $\alpha$ be an integer such that $\alpha>d/2+2$, and let $\tilde{\Pi}$ be the prior from \eqref{Eq: specific rescaled GP prior} based on a rescaled Whittle-Mat\'{e}rn process (which provides a canonical choice of prior in this problem if it is only known that $f^*\in C^{\alpha}(\mathcal{O})$, see \cite{giordanoConsistencyBayesianInference2020}).		
		Then there exists $f^*$ of the form \eqref{Eq: intro generalised additive model} with $F\in C^{\alpha}(\R)$ and $g_1,\ldots,g_d\in C^{\infty}(\R)$ such that
		$$ E_{f^*}  \Pi\left( f: \|f-f^*\|_{L^2} > n^{-a} \mid D_n \right) \to 0, \quad  E_{f^*}  \tilde{\Pi}\left( f: \|f-f^*\|_{L^2} \leq n^{-b} \mid D_n \right) \to 0, $$
		for exponents $a,b>0$ depending on $\alpha,d$ which, for sufficiently large $d$, satisfy $n^{-a}\ll n^{-b}$.
	\end{inftheorem}
	This statement is implied by Corollary \ref{Cor: specific rates} and Theorem \ref{Thm: contraction rate lower bound} below. The upper bound for the DGP prior holds uniformly over all such choices of $f^*$ with $\|F\|_{C^{\alpha}}\leq K$, and is achieved without knowledge of the structure \eqref{Eq: intro generalised additive model} or the precise value of $\alpha$. While $\tilde{\Pi}$ depends on knowing $\alpha$, the lower bound also holds for hierarchical priors with randomised smoothness: see Remark \ref{Remark: lb for hierarchical Bayes}. We see that if the dimension $d$ is large enough, asymptotically the posterior arising from the DGP prior places almost all of its mass inside an $L^2$-ball of radius $n^{-a}$ centred at the true $f^*$, while the Gaussian process prior $\tilde{\Pi}$ induces a posterior which places almost all of its mass outside of a larger neighbourhood, with radius $n^{-b}$. So in this case, the DGP prior outperforms the rescaled Whittle-Mat\'{e}rn process prior $\tilde{\Pi}$.
	
	The usual choice of parameter space for $f$ is a ball in a suitably regular Sobolev or H\"{o}lder space, as in \cite{nicklBayesianNonlinearStatistical2023}; over such classes, the Gaussian-based prior $\tilde{\Pi}$ performs well in a minimax sense. As these parameter spaces are special cases of the compositional classes introduced in Section \ref{Subsection: compositional functions}, the DGP prior achieves fast contraction rates over these parameter spaces, though not as fast as $\tilde{\Pi}$: see Remark \ref{Remark: can't use Sobolev} and the discussion after Corollary \ref{Cor: specific rates}. Our results indicate that if one is willing to pay the additional computational cost to use the DGP prior (see Section \ref{Subsection: DGP prior discussion}) instead of a Gaussian-based prior such as $\tilde{\Pi}$, then the reward is fast convergence rates that reflect the compositional structure of the unknown parameter $f^*$, which typical Gaussian-based priors are unable to leverage.
	
	
	The paper is structured as follows: Section \ref{Section: setting} introduces the general inverse problem we study, as well as the compositional classes of functions which provide our parameter spaces. Section \ref{Section: DGP Prior} introduces the DGP prior, while Section \ref{Section: contraction rates} contains the contraction rate results for this prior. In Section \ref{Section: crlbs} we explore the sub-optimality of particular Gaussian process priors for modelling compositional functions, and compare their performance to that of the DGP prior. Section \ref{Section: discussion} contains some broader discussion on deep Gaussian processes. Proofs are deferred to Appendix \ref{Appendix: proofs}, while Appendix \ref{Appendix: PDE facts} reviews theory for the two specific PDE inverse problems we have discussed.
	
	\section{Setting} \label{Section: setting}
	
	\subsection{Notation}
	
	In this section, $\mathcal{X}$ stands for either a smooth domain $\mathcal{O}\subset\Rd$ (that is, a non-empty, open, bounded set with smooth boundary $\partial\mathcal{O}$) or the unit cube $[-1,1]^d$.
	
	We respectively define $C(\mathcal{X})$ and $L^{\infty}(\mathcal{X})$ to be the sets of all bounded continuous and essentially bounded measurable functions $\mathcal{X}\to\R$, each endowed with the supremum norm $\|\cdot\|_{\infty}$. Let $L^2(\mathcal{X}) = H^0(\mathcal{X})$ denote the usual space of square-integrable functions on $\mathcal{X}$, endowed with its norm $\|\cdot\|_{L^2}$. For $\beta>0$, let $C^{\beta}(\mathcal{X})$ and $H^{\beta}(\mathcal{X})$ respectively denote the usual H\"{o}lder and Sobolev spaces over $\mathcal{X}$, see Appendix \ref{Appendix: PDE facts} for details. We recall the Sobolev embedding $H^{\beta}(\mathcal{X}) \subset C^{\beta - d/2}(\mathcal{X})$ which holds for all $\beta>d/2$.
	
	Let $A(\mathcal{X})$ be any of the above function spaces. We write $B_{A(\mathcal{X})}(R)$ to denote the norm-ball of radius $R$ in $A(\mathcal{X})$. When $\mathcal{X}=\mathcal{O}$, for any compactly contained subset $\mathcal{K}\subset\mathcal{O}$ we write $A_{\mathcal{K}}(\mathcal{O})$ for the subspace of functions in $A(\mathcal{O})$ whose support is contained in $K$. We also write $A_c(\mathcal{O})$ to denote the functions in $A(\mathcal{O})$ whose support is compactly contained in $\mathcal{O}$. When $\mathcal{X}=[-1,1]^d$, we write $A_d = A([-1,1]^d)$ (for example, $L^2_d = L^2([-1,1]^d))$. 
	
	Throughout the paper we use $\lesssim$ and $\gtrsim$ to denote inequalities holding up to a constant, whose dependence on model parameters will be specified. For sequences $a_n, b_n$, we write $a_n\simeq b_n$ if $a_n\lesssim b_n$ and $b_n\lesssim a_n$. Finally, we denote by $\mathcal{L}(Z)$ the law of the random variable $Z$.

	\subsection[General Inverse Problem]{A General Statistical Non-Linear Inverse Problem}
	
	Fix a smooth domain $\mathcal{O}\subset\Rd$, and let $\Theta$ be a measurable subset of $L^2(\mathcal{O})$. Suppose we are given a \emph{forward map} $\G:\Theta\to L^2(\mathcal{O})$. Our goal is to recover $\theta\in\Theta$ given noisy observations of $\G(\theta)$: we observe independent and identically distributed (i.i.d.) pairs $(Y_i,X_i)_{i=1}^n$ from the model
	\begin{equation}\label{Eq: observation model}
		Y_i = \G(\theta)(X_i) + \epsilon_i, \quad \epsilon_i\overset{\mathrm{i.i.d.}}{\sim} N(0,1), \quad 1\leq i \leq n,
	\end{equation}
	where the covariates $X_i$ are i.i.d. draws from the uniform distribution $\mu$ on $\mathcal{O}$, independent of the $\epsilon_i$. We write $P_{\theta}$ for the law of $(Y_1,X_1)$; denoting by $\ud y$ the Lebesgue measure on $\R$, $P_{\theta}$ has Radon-Nikodym density with respect to $\ud y\times \ud\mu$ given by
	\begin{equation}\label{Eq: single obs likelihood}
		p_{\theta}(y,x) := \frac{\ud P_{\theta}}{\ud y\times\ud\mu}(y,x) = \frac{1}{\sqrt{2\pi}}\exp\left\{-\frac{1}{2}[y-\G(\theta)(x)]^2\right\}.
	\end{equation}
	Denote the full data vector $(Y_i,X_i)_{i=1}^n$ by $D_n$; by a slight abuse of notation, we also denote by $P_{\theta}$ the law of $D_n$, and by $E_{\theta}$ the corresponding expectation.
	
	Let $\Pi$ be a prior (i.e. a Borel probability measure) supported on the Banach space $C(\mathcal{O})$. Then the map $(\theta,(y,x))\mapsto p_{\theta}(y,x)$ is jointly measurable and so by Bayes' formula (a version of) the posterior is given by
	\begin{equation}\label{Eq: Bayes' formula}
		\Pi(B\mid D_n) = \frac{\int_B e^{\ell_n(\theta)}\,\ud\Pi(\theta)}{\int_{C(\mathcal{O})} e^{\ell_n(\theta)}\,\ud\Pi(\theta)}, \quad \text{ any measurable } B\subset C(\mathcal{O}),
	\end{equation}
	where the joint log-likelihood function is (up to an additive constant) given by
	\begin{equation}\label{Eq: log-likelihood}
		\ell_n(\theta) = -\frac{1}{2}\sum_{i=1}^n[Y_i - \G(\theta)(X_i)]^2,
	\end{equation}
	see p.7 of \cite{ghosalFundamentalsNonparametricBayesian2017}.
	
	We impose the following requirements on the forward map $\G$, adapted from the conditions on $\G$ from Chapter 2 of \cite{nicklBayesianNonlinearStatistical2023}. The first condition says that the forward map is uniformly bounded over $\Theta\times\mathcal{O}$.
	\begin{condition}[Uniform Boundedness of $\G$]\label{Cond: G unif boundedness}
		Assume that there exists a constant $U<\infty$ depending on $\G,\Theta,\mathcal{O}$ such that
		\begin{equation}\label{Eq: G unif boundedness}
			\sup_{\theta\in\Theta} \| \G(\theta) \|_{\infty} \leq U.
		\end{equation}
	\end{condition}
	The next condition imposes Lipschitz continuity of $\G$ over a suitable subset of regular functions.
	\begin{condition}[Lipschitz Continuity of $\G$]\label{Cond: forward Lipschitz}
		Assume that there exists $\beta\geq0$ such that for all $M>0$, there exists a constant $L>0$ (possibly depending on $\G,\Theta,\mathcal{O}$ and $M$) such that 
		\begin{equation}\label{Eq: forward Lipschitz condition}
			\|\G(\theta_1) - \G(\theta_2)\|_{L^2} \leq L\|\theta_1-\theta_2\|_{\infty} \quad \forall \theta_1,\theta_2\in\Theta\cap B_{C^{\beta}(\mathcal{O})}(M).
		\end{equation}
	\end{condition}
	
	The final condition is a \emph{stability estimate} for $\G$, which provides quantitative control of the injectivity of the forward map.
	\begin{condition}[Stability Estimate] \label{Cond: stability estimate}
		Assume there exists $\beta\geq0,L'>0,\xi>0$ and $\zeta>0$ such that for all $M>0$ and all $\delta>0$ sufficiently small,
		\begin{equation}\label{Eq: stability estimate}
			\sup\bigg\{ \|\theta-\theta^*\|_{L^2}:\theta\in\Theta\cap B_{C^{\beta}(\mathcal{O})}(M), \|\G(\theta)-\G(\theta^*)\|_{L^2} \leq \delta \bigg\} \leq L' M^{\xi}\delta^{\zeta}.
		\end{equation}
	\end{condition}
	
	The left-hand side of \eqref{Eq: forward Lipschitz condition} is known as the \emph{($L^2$-)prediction risk}. Under Condition \ref{Cond: G unif boundedness}, important information theoretic quantities such as the Kullback-Leibler divergence, the Kullback-Leibler variation and the Hellinger distance are all dominated by the prediction risk (c.f. Proposition 1.3.1 in \cite{nicklBayesianNonlinearStatistical2023}). Since the general theory of posterior contraction rates yields results for the Hellinger distance (see Chapter 8 of \cite{ghosalFundamentalsNonparametricBayesian2017}), we first obtain a contraction rate in prediction risk and then apply the stability estimate from Condition \ref{Cond: stability estimate} to convert this into an $L^2$-contraction rate for $\theta$. The forward Lipschitz estimate from Condition \ref{Cond: forward Lipschitz} is used to verify small ball and metric entropy conditions central to the general theory of posterior contraction rates.
	
	\begin{remark}[Compositional priors cannot leverage forward smoothing]\label{Remark: can't use Sobolev}
		Typically one can prove a better Lipschitz estimate for $\mathcal{G}$ than \eqref{Eq: forward Lipschitz condition}, with a weak Sobolev norm in place of the supremum norm on the right-hand side: for example Condition 2.1.1 in \cite{nicklBayesianNonlinearStatistical2023}, which is then verified for Darcy's problem with the $\left(H^1\right)^*$-norm (Proposition 2.1.3, ibid) and for the Schr\"{o}dinger problem with the $\left(H^2\right)^*$-norm (Exercise 2.4.1, ibid). Reproducing kernel Hilbert spaces describing the covariance structure of Gaussian priors have some compatibility with these dual Sobolev norms that enables the use of these refined Lipschitz estimates to prove fast contraction rates.
		
		However, when using a prior whose draws are compositional functions (such as the DGP prior), one must use `pointwise' norms since these are the only norms which behave well with respect to composition: there is no analogue of the key technical tool Lemma \ref{Lemma: continuity of composition} for the $(H^1)^*$- or $(H^2)^*-$norms. It therefore seems unlikely that one could use the DGP prior and still leverage the forward smoothing property of $\mathcal{G}$.
	\end{remark}
	
	\subsection[Examples]{Examples: Darcy's Problem and the Steady-State Schr\"{o}dinger Equation}\label{Subsection: PDE examples}
	
	We define the two specific PDE-constrained inverse problems we consider and give a summary of the above conditions for the associated forward maps. See Appendix \ref{Appendix: PDE facts} for a more detailed confirmation of Conditions \ref{Cond: G unif boundedness}-\ref{Cond: stability estimate}.
	
	\textbf{Darcy's Problem} Let $\mathcal{O}\subset\Rd$ be a given smooth domain. We wish to recover the scalar diffusivity function $f\in C^{\gamma}(\mathcal{O})$ ($\gamma>1$) from observations of the solution $u$ to the PDE
	\begin{equation}\label{Eq: Darcy's problem}
		\begin{aligned}
			\nabla\cdot(f\nabla u) &= g \quad \text{ on }\mathcal{O}, \\
			u&=0 \quad \text{ on }\partial\mathcal{O}.
		\end{aligned}
	\end{equation}
	The source term $g$ is known and assumed to be smooth and satisfy $g\geq g_{\min}$ on $\mathcal{O}$ for some $g_{\min}>0$. One may view \eqref{Eq: Darcy's problem} as a steady-state heat equation, where $f$ is the diffusivity and $u$ describes the temperature; alternatively, \eqref{Eq: Darcy's problem} describes a steady-state groundwater flow problem where $u$ is the distribution of water through $\mathcal{O}$ and $f$ is the permeability. Darcy's problem has been studied extensively in the inverse problems literature: see  \cite{bonitoDiffusionCoefficientsEstimation2017,giordanoConsistencyBayesianInference2020} and references therein.
	Assuming that $f\geq K_{\min}>0$ on $\mathcal{O}$, standard elliptic theory (e.g. \cite{gilbargEllipticPartialDifferential2001}) tells us that the solution $u_f$ to \eqref{Eq: Darcy's problem} is unique and lies in $C^{\gamma+1}(\mathcal{O})$. Define the solution map $G:f\mapsto u_f$.
	
	The condition $f\geq K_{\min}$ is not compatible with placing a Gaussian prior on $f$ directly. We therefore use a \emph{link function}: given $\theta\in\Theta\subset C^{\gamma}(\mathcal{O})$, we set
	\begin{equation}\label{Eq: link function}
		f_{\theta} = K_{\min} + e^{\theta},
	\end{equation}
	and define the forward map as
	\begin{equation}\label{Eq: forward map}
		\G:\Theta\to L^2(\mathcal{O}), \quad \G(\theta) = G(f_{\theta}).
	\end{equation}
	The following properties of $\mathcal{G}$ are established in Appendix \ref{Appendix: PDE facts}.
	\begin{lemma*}
		The forward map $\mathcal{G}$ defined in \eqref{Eq: forward map} satisfies Condition \ref{Cond: G unif boundedness}, Condition \ref{Cond: forward Lipschitz} for any $\beta\geq1$ and Condition \ref{Cond: stability estimate} for any integer $\beta>1$ with $\xi = \beta(\beta+1)$ and $\zeta = \frac{\beta-1}{\beta+1}$.
	\end{lemma*}
	We will state our contraction rate results for $\theta$. Due to the smoothness of the link function \eqref{Eq: link function}, these imply the same contraction rates for $f_{\theta}$ (c.f. Lemma \ref{Lemma: link function properties}). 
	
	\vspace{10pt}
	
	\textbf{Steady-State Schr\"{o}dinger Equation} Let $\mathcal{O}\subset\Rd$ be a given smooth domain. We wish to recover the `absorption potential' $f\geq0$ from observations of the solution $u$ to the equation
	\begin{equation}\label{Eq: Schrodinger problem}
		\begin{aligned}
			\frac{1}{2}\Delta u - fu &= 0 \quad \text{ on }\mathcal{O}, \\
			u &= h \quad \text{ on }\partial\mathcal{O}.
		\end{aligned}
	\end{equation}
	The boundary temperatures $h$ are assumed to be known and smooth, and to satisfy $h\geq h_{\min}>0$ on $\partial\mathcal{O}$. This is a steady-state version of the time-dependent Schr\"{o}dinger equation ubiquitous in quantum physics, where $f$ describes some attenuation effect. This problem has been studied from a Bayesian point of view in \cite{nicklBernsteinMisesTheorems2020a,nicklConvergenceRatesPenalized2020,kekkonenConsistencyBayesianInference2022}.
	
	So long as $f\in C^{\gamma}(\mathcal{O})$ for some $\gamma>0$ and $f\geq0$, again by the usual elliptic PDE theory there exists a unique solution $u_f\in C^{\gamma+2}(\mathcal{O})$. Similarly to the previous problem, the non-negativity constraint on $f$ means that we cannot place a prior whose support is a linear space directly on $f$. Instead, we use the link function
	\begin{equation}\label{Eq: link function 2}
		f_{\theta} = e^{\theta}, \quad \theta\in \Theta \subset C^{\gamma}(\mathcal{O}).
	\end{equation}
	Define the forward map $\mathcal{G}$ as in \eqref{Eq: forward map} with $G:f\mapsto u_f$ the solution map for \eqref{Eq: Schrodinger problem}. 
	\begin{lemma*}
		The above forward map $\mathcal{G}$ satisfies Condition \ref{Cond: G unif boundedness}, Condition \ref{Cond: forward Lipschitz} for any $\beta\geq0$, and Condition \ref{Cond: stability estimate} for any choice of $\beta>0$ with $\xi=\beta/2 + 1$ and $\zeta=\frac{\beta}{\beta+2}$.
	\end{lemma*}
	Again, by Lemma \ref{Lemma: link function properties} contraction rates for $\theta$ carry over to $f_{\theta}$.

	\subsection{Compositional Functions}\label{Subsection: compositional functions}
	
	In previous works studying inverse problems of the type described above (such as \cite{nicklBayesianNonlinearStatistical2023}), it is often assumed that the parameter $\theta$ lies in some Sobolev or H\"{o}lder space. Instead, we will model $\theta$ as a compositional function, in the manner of \cite{schmidt-hieberNonparametricRegressionUsing2020,finocchioPosteriorContractionDeep2023}. The next two subsections will give examples of compositional structures in the two example inverse problems. Assume that for some integer $q$, we can write $\theta$ in the form
	\begin{equation}\label{Eq: compositional function}
		\theta =\bar{\theta}_q\circ\dots\circ\bar{\theta}_0,
	\end{equation}
	i.e. a composition of $(q+1)$ functions $\bar{\theta}_i$. The $\bar{\theta}_i$ have the following domains and codomains:
	\begin{align*}
		&\bar{\theta}_0: \mathcal{O}\to[-1,1]^{d_1}, \\ 
		&\bar{\theta}_i:[-1,1]^{d_i}\to[-1,1]^{d_{i+1}}, \quad 1\leq i\leq q-1, \\
		&\bar{\theta}_q:[-1,1]^{d_q}\to\R.
	\end{align*}
	The choice of the cubes $[-1,1]^{d_i}$ is not restrictive, since the final function can take values in the whole of $\R$. Moreover, we assume without loss of generality that $d_i\leq d$ for all $i$, since we will eventually assume that each $\bar{\theta}_i$ is continuous and so its domain can always be embedded into a $d$-dimensional manifold (namely $\bar{\theta}_{i-1}\circ\cdots\circ\bar{\theta}_0(\mathcal{O})$). For each $i$, we write $\bar{\theta} = (\bar{\theta}_{ij})_{j=1}^{d_{i+1}}$, where $d_{q+1}=1$. Each of the $\bar{\theta}_{ij}$ takes values in the interval $[-1,1]$, with the exception of $\bar{\theta}_{q1}$ which takes values in $\R$.
	
	Of course, any function can be written in this form with $q=0$. The value of the compositional representation \eqref{Eq: compositional function} will come from reducing the dimensionality of the problem, or more precisely from allowing layers to trade off sparsity against smoothness. To that end, we will assume that each function $\bar{\theta}_{ij}$ only depends on a subset of its inputs $\mathcal{S}_{ij}\subset\{1,\ldots,d_i\}$ (here $d_0=d$; also, $d_{q+1}=1$). Write $t_i=\max_j \lvert \mathcal{S}_{ij} \rvert$  for the maximum size of such a subset. Note that $t_i\leq d_i$; we may assume that $t_i$ is the same for all $1\leq j\leq d_{i+1}$, although the sets $\mathcal{S}_{ij}$ can vary with $j$, since one can simply allow certain $\bar{\theta}_{ij}$ to `depend' on redundant variables. For any subset $S$ of indices, let $(\cdot)_S:x\mapsto x_S = (x_i)_{i\in S}$, and (understanding by abuse of notation that by $[-1,1]^{t_0}$ we mean the domain $\mathcal{O}$) define
	$$ \theta_{ij}:[-1,1]^{t_i} \to [-1,1], \quad x_{\mathcal{S}_{ij}}\mapsto\bar{\theta}_{ij}(x_{\mathcal{S}_{ij}},x_{\mathcal{S}_{ij}^c}), $$
	which is well-defined as $\bar{\theta}_{ij}$ does not depend on $x_{\mathcal{S}_{ij}^c}$ (for $i=q$ the codomain should strictly be $\R$, but often we will leave this to be understood by the reader for the sake of conciseness). Note that to specify $\bar{\theta}_{ij}$, it suffices to specify the function $\theta_{ij}$ which takes $t_i$ inputs, and the set $\mathcal{S}_{ij}$ identifying the $t_i$ relevant inputs.
	
	\begin{figure}
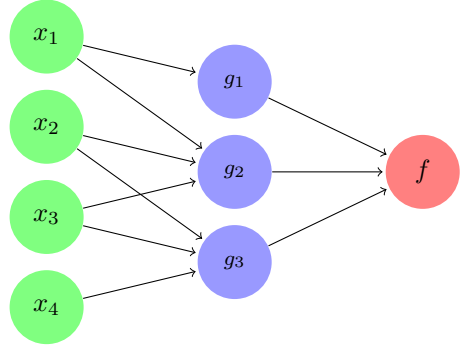

		\centering
		\begin{neuralnetwork}[height=4,nodesize=28pt,nodespacing = 12mm]
			\newcommand{\nodetextclear}[2]{}
			\newcommand{\nodetextg}[2]{\footnotesize \ifnum1=#2 {$g_1$} \fi \ifnum2=#2 {$g_2$} \fi \ifnum3=#2 {$g_3$} \fi}
			\newcommand{\nodetextxnb}[2]{\ifnum0=#2 \else $x_#2$ \fi}
			\newcommand{\logiclabel}[1]{\,{$\scriptstyle#1$}\,}
			\newcommand{\nodetextY}[2]{$y$}
			\setdefaultnodetext{\nodetextclear}
			
			\newcommand{\x}[2]{$x_#2$}
			\newcommand{\y}[2]{$f$}
			\newcommand{\hfirst}[2]{\small $h^{(1)}_#2$}
			\newcommand{\hsecond}[2]{\small $h^{(2)}_#2$}
			\inputlayer[count=4, bias=false, 
			text=\x		]
			\hiddenlayer[count=3, bias=false, 
			text=\nodetextg
			] 
			\link[from layer=0, to layer=1, from node=1, to node=1, style=black]
			\link[from layer=0, to layer=1, from node=1, to node=2, style=black]
			\link[from layer=0, to layer=1, from node=3, to node=2, style=black]
			\link[from layer=0, to layer=1, from node=2, to node=2, style=black]
			\link[from layer=0, to layer=1, from node=2, to node=3, style=black]
			\link[from layer=0, to layer=1, from node=3, to node=3, style=black]
			\link[from layer=0, to layer=1, from node=4, to node=3, style=black]
			\outputlayer[count=1, text=\y] \linklayers[style=black]
		\end{neuralnetwork}
		\caption{A schematic representing the function $\theta(x)= f(g_1(x_1),g_2(x_1,x_2,x_3),g_3(x_2,x_3,x_4))$. For this function we have $\theta_1=f$, $\theta_{0j}=g_j$, $d_0=4$, $t_0=3$,  $d_1=t_1=3$, $\mathcal{S}_{01}=\braces{1},$ $\mathcal{S}_{02}=\braces{1,2,3}$, $\mathcal{S}_{03}=\braces{2,3,4}$.}
		\label{Fig: graph of function}
	\end{figure}
	
	In summary, to specify such a function requires choosing the following parameters:
	\begin{itemize}
		\item a \emph{depth} $q\in\mathbb{N}$;
		\item a vector of \emph{dimensions} $\mathbf{d}\in\mathbb{N}^{q+1}$ such that $d_i\leq d$, where $d_0=d$ and $d_q=1$;
		\item a vector of \emph{intrinsic dimensions} $\mathbf{t}\in\mathbb{N}^{q+1}$ such that $t_i\leq d_i$;
		\item for each $i,j$, an \emph{active set} $\mathcal{S}_{ij}\subset\{1,\ldots,d_{i+1}\}$ of size $t_i$. Denote by $\mathcal{S}$ the set of all active sets;
		\item for each $i,j$, a function $\theta_{ij}:[-1,1]^{t_i}\to[-1,1]$.
	\end{itemize}
	See Figure \ref{Fig: graph of function} for an example of such a function; below we also discuss some concrete examples in inverse problems.
	
	We combine the first four structural parameters into a single parameter, called the \emph{graph} of the compositional function $\theta$, defined as
	\begin{equation}\label{Eq: graph of a function}
		\lambda:= (q,\mathbf{d},\mathbf{t},\mathcal{S}).
	\end{equation}
	
	The set of all possible graphs is denoted $\Lambda$. Once a graph is chosen, the compositional function $\theta$ can then be specified by choosing functions $\theta_{ij}$ for all relevant pairs $i,j$. Let $\boldsymbol{\alpha}=(\alpha_0,\ldots,\alpha_q)\in(0,\infty)^{q+1}$ be a vector of smoothnesses.
	Assuming that $\alpha_i > (1/2)t_i$ for all $i$, we define the parameter set
	\begin{equation}\label{Eq: compositional class defn}
		\Theta(\lambda,\boldsymbol{\alpha}) = \left\{ \theta \text{ of the form \eqref{Eq: compositional function}} : \theta \text{ has graph }\lambda,\,\, \theta_{ij}\in H^{\alpha_i}_{t_i}\,\, \forall{i,j}\right\}.
	\end{equation}
	(Recall the notational convention $H^{\alpha}_t = H^{\alpha}\left([-1,1]^t\right)$.) The condition on $\boldsymbol{\alpha}$ ensures (by Sobolev embedding) that the functions $\theta_{ij}$ are defined pointwise. Given a constant $K>0$, we also define
	\begin{equation}\label{Eq: bounded compositional class defn}
		\Theta(\lambda,\boldsymbol{\alpha},K) = \left\{ \theta\in\Theta(\lambda,\boldsymbol{\alpha}) : \theta_{ij}\in B_{H^{\alpha_i}_{t_i}}(K)\,\, \forall{i,j}\right\}.
	\end{equation}
	For $\theta\in\Theta(\lambda,\boldsymbol{\alpha})$, we combine the graph and smoothness parameters to form a new parameter
	\begin{equation}\label{Eq: structure parameter definition}
		\eta:=(\lambda,\boldsymbol{\alpha})
	\end{equation}
	which we call the \emph{structure} of $\theta$. The structure of a compositional function was shown to determine the minimax estimation rate (in a regression problem) in \cite{schmidt-hieberNonparametricRegressionUsing2020}. We denote by $\Omega$ the set of all structures.
	
	We assume that our true parameter $\theta^*$ lies in $\Theta(\lambda^*,\boldsymbol{\alpha}^*,K)$ for some structure $(\lambda^*,\boldsymbol{\alpha}^*)\in\Omega$ (which may be unknown) and some known $K>0$. Note that the representation of $\theta$ described by \eqref{Eq: compositional class defn} is not unique, and there may be several valid structures $\eta$ for a single function $\theta$. Our results should be interpreted as holding for whichever structure $\eta$ provides the best convergence rate.

	\subsubsection{Darcy's problem with layer structure}
	As briefly noted, Darcy's problem can be used to model groundwater flow, where the goal is to recover the scalar permeability function $f$. Permeability within a fixed type of rock varies relatively little, while different rocks and soils have permeability spanning multiple orders of magnitude, e.g.\ Chapter 6 of \cite{blunt_2017}. As such, a plausible approximate model for the permeability $f = f_{\theta}$ is that it is piecewise constant on (potentially unknown) regions. Such functions are modelled by the compositional structure of Section \ref{Subsection: compositional functions} (up to relaxing the Sobolev constraint, or taking a smooth approximation to the indicator functions in the below).

	Consider a collection of (hyper-)planes $\{ x \in \R^d :  \langle a_i,x \rangle = c_i\}_{i\leq k}$, where $a_i$ is a unit vector and $c_i\in \R$.
	Define $\bar{\theta}_{0,i}(x) = \langle a_i, x \rangle$, define $\bar{\theta}_{1,i}(x) = \ind\{x_i < c_i\}$ and define $\bar{\theta}_2(x) = \sum_{i=1}^k \alpha_i \prod_{j=1}^i x_j$ for some $\alpha_i\in \R$.
	Then $\bar{\theta}_2\circ\bar{\theta}_1\circ\bar{ \theta}_0$ is a function which is piecewise constant on the regions bounded by the planes, and can be made to take arbitrary values on each such region by choosing $\alpha_i$ appropriately.

	(Note the assumption in Section \ref{Subsection: compositional functions} that $d_i\leq d$ for all $i$ means that we can have at most $k=d$ such bounding planes; by adjusting the prior, this assumption can be relaxed to simply having a known bound on all $d_i$ in the optimal compositional structure and hence accommodate also functions most parsimoniously expressed as piecewise constant on regions bounded by $k> d$ planes. Also note that while described here for regions separated by planes, any boundary surface described by an equation $g(x)=c$ for a suitably smooth function $g$ is accommodated similarly.)
	
	Another way to model layer structure is to have $f$ constant in some directions. For example, if the soil consists of a single material, of density varying with depth, we may expect the permeability to depend only on the depth. This is also captured by the compositional model through taking $d=3$ and $\bar{\theta}_0(x)=x_3$, with $\bar{\theta}_1$ arbitrary. Let us emphasise once more that our model does not require prior knowledge of which form of layer structure, or any other compositional structure, is appropriate, rather picking up on this structure automatically.
	
	\subsubsection{Schr\"odinger equation with spherical symmetry}
		One of the first uses of the Schr\"odinger equation covered in introductory textbooks on quantum mechanics is for modelling a particle in a spherically symmetric potential, e.g.\ Chapter 2 of \cite{weinberg_2015}. The spherically symmetric Schr\"odinger equation so obtained is for example solved to find the energy levels of a hydrogen atom. Solving in this way requires prior knowledge of the symmetry; this structure can then be imposed directly in a Bayesian prior to achieve a fast, `one-dimensional' convergence rate.
	
		In contrast, the DGP method can discover unknown spherical symmetry. Indeed, any symmetric potential $f(x) = F(\lVert x\rVert^2)$ falls within the compositional class \eqref{Eq: compositional function} so long as $F$ is sufficiently smooth. Specifically, we may take $\bar{\theta}_0 = \lVert \cdot\rVert^2 \in C^\infty(\R^d),$ and $\bar{\theta}_1 = F \in C(\R)$. The contraction rate for $f$ in this setting (given in Corollary \ref{Cor: specific rates}) is also one-dimensional, achieved without prior knowledge of this symmetry.
		
		Note that we can also express $f$ as a generalised additive model \eqref{Eq: generalised additive model}, with $F=g$, $g_i(u)=u^2$ for all $i$; a version of Theorem \ref{Thm: contraction rate lower bound} concerning the non-optimality of a typical (non-deep) GP prior will hold, showing that this GP prior cannot take advantage of spherical symmetry.
		
	\section{Deep Gaussian Process Prior}\label{Section: DGP Prior}
	
	We construct a DGP prior which models compositional functions. We will first select a structure from a suitable hyperprior, and then draw each component function from an `elementary' process prior based on a Gaussian process. In both stages, a crucial role will be played by the convergence rates we are aiming to achieve. Given a dimension $t\in\mathbb{N}$ and a smoothness $\alpha>0$, define the rate
	\begin{equation}\label{Eq: rates defn}
		\varepsilon^{\alpha,t}_n := n^{-\frac{\alpha}{2\alpha+t}}.
	\end{equation}
	For vectors of smoothnesses $\boldsymbol{\alpha}$ and intrinsic dimensions $\mathbf{t}$, define the rate $\varepsilon^{\boldsymbol{\alpha},\mathbf{t}}_n := \max_{0\leq i\leq q} \varepsilon^{\alpha_i,t_i}_n$. Given a structure $\eta = (q,\mathbf{d},\mathbf{t},\mathcal{S},\boldsymbol{\alpha})$, we write
	\begin{equation}\label{Eq: compositional rate defn}
		\varepsilon^{\eta}_n = \varepsilon^{\boldsymbol{\alpha},\mathbf{t}}_n. 
	\end{equation}
	We also fix some smoothness $\beta>0$ such that Conditions \ref{Cond: forward Lipschitz} and \ref{Cond: stability estimate} hold for the forward map $\G$.
	
	\subsection{Elementary Process}
	
	We first introduce the `elementary' process, which is the prior distribution of each component $\theta_{ij}$ conditional on the structure parameter $\eta$. Our elementary process is based on the rescaled Gaussian priors first used for inverse problems in \cite{monardConsistentInversionNoisy2021}.
	
	For the moment, we consider the intermediate layers of the composition in \eqref{Eq: compositional function} (i.e. $1\leq i <q$). Assume a given intrinsic dimension $t\in\mathbb{N}$ and smoothness $\alpha>t/2$. Let $\Pi'_{\alpha,t}$ denote the law of a centred Gaussian process whose RKHS $\mathcal{H}$ embeds into $H^{\alpha}([-1,1]^t) = H^{\alpha}_t$ with equivalent norms. For example, we can use a suitable truncated series prior, or a Whittle-Mat\'{e}rn process: see Chapter 11 of \cite{ghosalFundamentalsNonparametricBayesian2017}. Define the rescaled prior
	\begin{equation*}
		\bar{\Pi}_{\alpha,t} = \mathcal{L}\left((\sqrt{n}\varepsilon^{\alpha,t}_n)^{-1}Z'\right), \quad Z'\sim\Pi'_{\alpha,t}.
	\end{equation*}
	
	Finally, we condition this process so that samples take values in $[-1,1]$ and are sufficiently regular to behave well under composition. For a constant $M_0\geq1$ to be chosen below, we obtain a prior
	\begin{equation}\label{Eq: elementary rescaled prior}
		\Pi_{\alpha,t} = \mathcal{L}\left(Z\mid \|Z\|_{\infty}\leq 1, \|Z\|_{C^{\beta}}\leq M_0 \right), \quad Z\sim\bar{\Pi}_{\alpha,t}.
	\end{equation}
	Note that due to the conditioning step, $\Pi_{\alpha,t}$ is not a Gaussian process. However, it is based on the Gaussian process prior $\bar{\Pi}_{\alpha,t}$, and the conditioning does not hugely alter the process since $\bar{\Pi}_{\alpha,t}$ concentrates on the conditioning set with high probability: as in the proof of Lemma 16 in \cite{giordanoConsistencyBayesianInference2020}, an application of the Borell-Sudakov-Tsirelson inequality (\cite{gineMathematicalFoundationsInfinitedimensional2016}, Theorem 2.5.8) gives that
	\begin{equation}\label{Eq: Borell-TIS on rescaled prior}
		\bar{\Pi}_{\alpha,t}\left(\|Z\|_{\infty}\leq 1, \|Z\|_{C^{\beta}}\leq M_0 \right) \geq 1 - \exp\left\{ -C_{\alpha,t}M_0^2n(\varepsilon_n^{\alpha,t})^2 \right\},
	\end{equation}
	for all $\alpha>t/2 + \beta$, where the constant $C_{\alpha,t}$ is decreasing in $\alpha$ and $t$. Some conditioning is necessary to achieve adaptive results using the techniques herein, as otherwise one does not achieve sufficiently good control of the $C^{\beta}$-norm on smooth models (the rate in the exponential inequality \eqref{Eq: Borell-TIS on rescaled prior} is not fast enough). Such control is required in Lemma \ref{Lemma: continuity of composition} below to control the effect of composing several such processes, as well to apply the Lipschitz and stability estimates \eqref{Eq: forward Lipschitz condition}, \eqref{Eq: stability estimate} .
	
	For the final layer, we wish to model a function $\theta_q:[-1,1]^t\to\R$. The construction is almost identical to the above, except that we do not condition on the event $\|Z\|_{\infty}\leq 1$. That is, the elementary process prior in the case of the final layer of the composition is
	$$ \Pi_{\alpha,t} = \mathcal{L}\left(Z\mid \|Z\|_{C^{\beta}}\leq M_0 \right), \quad Z\sim\bar{\Pi}_{\alpha,t}. $$
	
	For the first layer, to avoid technicalities associated with modelling the function near the boundary $\partial\mathcal{O}$, we assume that $\theta^*$ is supported on a known compact subset $\mathcal{K}\subset\mathcal{O}$; as $\mathcal{O}$ is open, $\mathcal{K}$ has some fixed positive distance from the boundary $\partial\mathcal{O}$. We can then model the components of the first layer $\theta_{0j}$ using (for example) a Whittle-Mat\'{e}rn process on $\mathcal{O}$ multiplied by a smooth cutoff function which equals 1 on $\mathcal{K}$; see Example 25 in \cite{giordanoConsistencyBayesianInference2020} for details. We then condition the process as in \eqref{Eq: elementary rescaled prior}. By an abuse of notation, we will simply write $\theta_{0j}:[-1,1]^d\to[-1,1]$ and leave these details to be understood by the reader.
	
	\subsection{Structure Hyperprior}
	
	We now describe the construction of the hyperprior on the structure of the function.
	
	Any probability density $\gamma$ on the set of structures $\Omega$ is fully determined by the conditional probability formula
	$$ \gamma(\eta) = \gamma(q)\gamma(\mathbf{d}\mid q)\gamma(\mathbf{t}\mid \mathbf{d},q)\gamma(\mathcal{S}\mid \mathbf{t},\mathbf{d},q)\gamma(\boldsymbol{\alpha}\mid\lambda). $$
	
	
	Fix a maximal smoothness $\alpha^+>\beta + d/2$, and define the interval $I(t_i):=[\beta + t_i/2,\alpha^+]$ (this interval is non-empty as $t_i\leq d_i\leq d$). Write $\Omega'\subset\Omega$ for the subset of structures satisfying $d_i\leq d$ and $\alpha_i\in I(t_i)$ for all $i$. We make the following assumption on $\gamma$, based on Assumption 1 in \cite{finocchioPosteriorContractionDeep2023}. 
	\begin{assumption}\label{Assumption 1}
		Assume that for any $\lambda\in\Lambda$, the distribution of smoothnesses $\gamma(\cdot\mid\lambda)$ equals the law of each $\alpha_i$ drawn independently and uniformly at random from the interval $I(t_i)$. Moreover, we assume that $\gamma$ is independent of $n$, $\gamma$ is supported on $\Omega'$, $\gamma(\eta)>0$ for all $\eta\in\Omega'$, and $\int_{\Omega'}\sqrt{\gamma(\eta)}\,\ud\eta <\infty$.
	\end{assumption}
	(We insist above that $\gamma(\cdot\mid\lambda)$ is uniform in order to simplify our proofs; however, if one chooses any density $\gamma$ on $\Omega'$ that is bounded, bounded away from zero, and such that $\gamma$ satisfies the square-root integrability condition, then the proofs of Theorems \ref{Thm: prediction risk contraction rate} and \ref{Thm: L2 contraction rate} still work.)

	We will not use $\gamma$ directly as our structure hyperprior, but rather a penalised version which ensures that with high probability we draw structures that are not too complex. 
	We then consider the hyperprior
	\begin{equation}\label{Eq: structure hyperprior}
		\pi(\eta) \propto e^{-\Psi_n(\eta)}\gamma(\eta), \quad \Psi_n(\eta) := n(\varepsilon^{\eta}_n)^2 + e^{e^{|\mathbf{d}|_1}},
	\end{equation}
	where $|\mathbf{d}|_1 = \sum_{i}|d_i|$ is the $\ell^1$-norm of $\mathbf{d}$. Note that $\pi$ is well-defined since $0<e^{-\Psi_n(\eta)}\leq 1$ and $\int\gamma(\eta)\,\ud\eta = 1$. The normalising constant of proportionality is therefore bounded above by 1.
	
	\subsection{Construction of the DGP Prior}
	
	Given a structure $\eta\in\Omega$, we construct the deep process as follows: for $0\leq i\leq q, 1\leq j\leq d_{i+1}$, we take $Z_{ij}$ to be independent draws from the elementary process prior $\Pi_{\alpha_i,t_i}$ as defined in \eqref{Eq: elementary rescaled prior} (with the necessary modifications for $i=0,q$). We set $Z_i = (Z_{ij})_{1\leq j\leq d_{i+1}}$ and finally $Z=Z_q\circ\dots\circ Z_0$. The law of this resulting $Z$ is denoted $\Pi(\cdot\mid\eta)$.
	
	The overall DGP prior is the measure $\Pi$, where
	\begin{equation}\label{Eq: deep GP prior defn}
		\Pi \mid \eta = \Pi(\cdot\mid \eta), \quad \eta \sim \pi.
	\end{equation}
	The deep GP prior depends on $n$ (both through the penalisation term in $\pi$ and the rescalings in $\Pi(\cdot\mid\eta)$) but we leave this implicit. Note that since $\pi$ is supported on structures in $\Omega'$, by Sobolev embedding and the fact that the composition of continuous functions is continuous, $\Pi$ is supported on $C(\mathcal{O})$. Thus Bayes' formula \eqref{Eq: Bayes' formula} holds for the DGP prior $\Pi$.
	
	\section{Contraction Rates} \label{Section: contraction rates}
	
	Fix some $\beta\geq1$ such that Conditions \ref{Cond: forward Lipschitz} and \ref{Cond: stability estimate} hold for the forward map $\G$. Let $\Pi$ be the DGP prior constructed in the previous section. Recall the definition of the structure parameter $\eta$ from \eqref{Eq: structure parameter definition}, and let $D_n\sim P^n_{\theta^*}$ be data generated according to \eqref{Eq: observation model}, where $\theta^*\in\Theta(\eta^*)$ for some $\eta^*\in\Omega'$. We let $\Pi(\cdot\mid D_n)$ be the posterior distribution based on $D_n$, as defined through \eqref{Eq: Bayes' formula}.
	
	Let $\mathcal{K}\subset\mathcal{O}$ be a known compact set, and write $\Theta_{\mathcal{K}}(\eta^*,K) = \{\theta\in\Theta(\eta^*,K): \theta\mid_{\mathcal{K}^c}\equiv 0\}$. Our first result establishes a contraction rate in prediction risk, which holds uniformly for $\theta^*\in\Theta_{\mathcal{K}}(\eta^*,K)$. Moreover, it shows that with high probability, posterior draws have controlled $C^{\beta}$-norm.
	\begin{theorem}\label{Thm: prediction risk contraction rate}
		Let $\Pi$ be the DGP prior as constructed above. Assume that $\eta^*\in\Omega'$, and let $K>0$. If $M_0$ in \eqref{Eq: elementary rescaled prior} is chosen sufficiently large depending only on $K$, then for any $\delta>\log{M_0}$ we have that
		$$ \sup_{\theta^*\in \Theta_{\mathcal{K}}(\eta^*,K)}E_{\theta^*} \Pi\left( \theta: \|\G(\theta) - \G(\theta^*)\|_{L^2(\mathcal{O})}\geq (\log{n})^{\delta}\varepsilon_n^{\eta^*} \mid D_n\right) \to 0$$
		as $n\to\infty$, where $\varepsilon_n^{\eta^*}$ is defined in \eqref{Eq: compositional rate defn}. Moreover, as $n\to\infty$ we have that
		$$ \sup_{\theta^*\in \Theta_{\mathcal{K}}(\eta^*,K)}E_{\theta^*} \Pi\left( \theta: \|\theta\|_{C^{\beta}(\mathcal{O})}  \geq (\log{n})^{\delta} \mid D_n\right) \to 0 .$$
	\end{theorem}
	The proof of Theorem \ref{Thm: prediction risk contraction rate} is given in Appendix \ref{Appendix: proofs}, and uses ideas from Theorems 1 and 2 from \cite{finocchioPosteriorContractionDeep2023} together with techniques from the Bayesian approach to nonlinear inverse problems described in \cite{nicklBayesianNonlinearStatistical2023}.
	As posterior draws have bounded $C^{\beta}$-norm with high probability, the stability estimate \eqref{Eq: stability estimate} immediately yields a contraction rate for $\theta$ in the $L^2$-distance.
	\begin{theorem}\label{Thm: L2 contraction rate}
		Assume $\beta\geq1$ is an integer. Under the conditions of Theorem \ref{Thm: prediction risk contraction rate}, we have for the constants $L'>0,\xi>0,\zeta>0$ from \eqref{Eq: stability estimate} that
		$$ \sup_{\theta^*\in \Theta_{\mathcal{K}}(\eta^*,K)} E_{\theta^*} \Pi\left( \theta: \|\theta - \theta^*\|_{L^2(\mathcal{O})}\geq L'(\log{n})^{\delta(\xi+\zeta)}(\varepsilon_n^{\eta^*})^{\zeta} \mid D_n\right) \to 0$$
		as $n\to\infty$.
	\end{theorem}
	
	The logarithmic terms in Theorems \ref{Thm: prediction risk contraction rate} and \ref{Thm: L2 contraction rate} are needed to control the unboundedness of the depth $q$ in the structure hyperprior; if the true depth $q^*$ is known or $q$ is assumed to be bounded above then these terms can be replaced by (large) constants.
	
	The general results above yield the following contraction rates in our two specific inverse problems: Darcy's problem and the Schr\"{o}dinger problem.
	\begin{corollary} \label{Cor: specific rates}
		The DGP prior attains the following convergence rates in specific inverse problems:
		\begin{enumerate}
			\item Consider Darcy's problem, defined in (\ref{Eq: Darcy's problem}-\ref{Eq: forward map}). Fix an integer $\beta>1$, a compact $\mathcal{K}\subset\mathcal{O}$, some $K>0$ and let $\Pi$ be the DGP prior given by \eqref{Eq: deep GP prior defn} with the constant $M_0$ in \eqref{Eq: elementary rescaled prior} chosen as in Theorem \ref{Thm: prediction risk contraction rate}. For any $\eta^*\in\Omega'$ and any $\delta>\log{M_0}$, there exists a constant $C>0$ such that for any $\theta^*\in\Theta_{\mathcal{K}}(\eta^*,K)$,
			$$ E_{\theta^*} \Pi\left( \theta: \|\theta - \theta^*\|_{L^2(\mathcal{O})}\geq C(\log{n})^{\bar{\delta}}(\varepsilon_n^{\eta^*})^{\frac{\beta-1}{\beta+1}} \mid D_n\right) \to 0 $$
			as $n\to\infty$, where $\bar{\delta} = \delta\left(\beta(\beta+1) + \frac{\beta-1}{\beta+1}\right)$. In the special case that $\theta^*\in B_{H^{\alpha}_{\mathcal{K}}(\mathcal{O})}(K)$, this becomes
			$$ E_{\theta^*} \Pi\left( \theta: \|\theta - \theta^*\|_{L^2(\mathcal{O})}\geq C(\log{n})^{\bar{\delta}}n^{-\frac{\alpha(\beta-1)}{(2\alpha+d)(\beta+1)}} \mid D_n\right) \to 0 .$$
			\item Consider the Schr\"{o}dinger equation problem defined in (\ref{Eq: Schrodinger problem}-\ref{Eq: link function 2}). Fix an integer $\beta>0$, a compact $\mathcal{K}\subset\mathcal{O}$, some $K>0$ and let $\Pi$ be the DGP prior as before. For any $\eta^*\in\Omega'$ and any $\delta>\log{M_0}$, there exists a constant $C>0$ such that for any $\theta^*\in\Theta_{\mathcal{K}}(\eta^*,K)$,
			$$ E_{\theta^*} \Pi\left( \theta: \|\theta - \theta^*\|_{L^2(\mathcal{O})}\geq C(\log{n})^{\bar{\delta}}(\varepsilon_n^{\eta^*})^{\frac{\beta}{\beta+2}} \mid D_n\right) \to 0 $$
			as $n\to\infty$, where $\bar{\delta} = \delta\left(\frac{\beta+2}{2}+\frac{\beta}{\beta+2}\right)$. In the special case where we have $\theta^*\in B_{H^{\alpha}_{\mathcal{K}}(\mathcal{O})}(K)$, this becomes
			$$ E_{\theta^*} \Pi\left( \theta: \|\theta - \theta^*\|_{L^2(\mathcal{O})}\geq C(\log{n})^{\bar{\delta}}n^{-\frac{\alpha\beta}{(2\alpha+d)(\beta+2)}} \mid D_n\right) \to 0 .$$
		\end{enumerate}
	\end{corollary}
	
	We may compare our contraction rates, both in prediction risk and $L^2$ risk, with those obtained elsewhere in the literature when $\theta^*$ lies in a Sobolev or H\"{o}lder ball. 
	
	Let us first consider Darcy's problem: over a ball in $H^{\alpha}(\mathcal{O})$, the optimal prediction risk rate derived from \cite[Theorem 10]{nicklConvergenceRatesPenalized2020} is $n^{-\frac{\alpha+1}{2\alpha+2+d}}$. In \cite{giordanoConsistencyBayesianInference2020}, the rescaled prior $\bar{\Pi}_{\alpha^*,d}$ is shown to attain this rate in prediction risk (Theorem 4, ibid), and this rate to the power $\frac{\beta-1}{\beta+1}$ in $L^2$ risk by virtue of the stability estimate \eqref{Eq: stability estimate} (Theorem 5, ibid). In fact, for a specific choice of Gaussian process prior, one may improve the exponent $\zeta$ to $\frac{\alpha-1}{\alpha+1}$: see Exercise 2.4.3 in \cite{nicklBayesianNonlinearStatistical2023}. The optimal $L^2$ recovery rate for $\theta$ in Darcy's problem is not currently known.	More precise results exist for the steady-state Schr\"{o}dinger equation problem described in (\ref{Eq: Schrodinger problem}-\ref{Eq: link function 2}). A Bayesian approach to solving this problem was studied in \cite{nicklBernsteinMisesTheorems2020a}; there it was shown that over a ball in the H\"{o}lder space $C^{\alpha}(\mathcal{O})$, $\alpha>2+d/2$, the minimax $L^2$-risk for recovering the parameter $\theta^*$ is $n^{-\frac{\alpha}{2\alpha+2+d}}$ (Proposition 2, ibid). Moreover, a prior based on a random wavelet expansion was constructed which contracts about any $\theta^*\in C_c^{\alpha}(\mathcal{O})$ at this rate up to a logarithmic factor (Theorem 1, ibid); see Exercises 2.4.1 and 2.4.3 in \cite{nicklBayesianNonlinearStatistical2023} for a Gaussian prior which contracts at the minimax rate. In both problems, these contraction rates are faster than those achieved by the DGP prior in Corollary \ref{Cor: specific rates} for the reason discussed in Remark \ref{Remark: can't use Sobolev}.
	
	The gap between these rates and the DGP contraction rates from Corollary \ref{Cor: specific rates} suggests that there is a cost to adapting to arbitrary compositional structures $\eta$. We note that the rates achieved by the DGP prior are still `fast' and if $\theta^*$ is very smooth, one may choose $\beta$ to be large so that the contraction rates are both close to $n^{-1/2}$. Moreover, the DGP prior is able to adapt to an unknown structure $\eta$ by virtue of the carefully chosen structure hyperprior $\pi$ in \eqref{Eq: structure hyperprior}. As we shall see in the next section, when $\theta^*$ has the prototypical compositional structure of a generalised additive model \eqref{Eq: generalised additive model}, using a Gaussian process prior which ignores this structure can lead to a substantially slower contraction rate.  
	
	\section{Sub-Optimality of Gaussian Priors in Compositional Models} \label{Section: crlbs}
	
	In this section, we work in Darcy's problem defined in (\ref{Eq: Darcy's problem}-\ref{Eq: forward map}) to fix ideas; analogous results hold in other settings.
	
	We have seen that the deep Gaussian process prior can successfully leverage compositional structure of the underlying true parameter $\theta^*$ to attain fast convergence rates. Even over Sobolev balls with no additional compositional structure, the DGP prior achieves a polynomial contraction rate, almost as fast as a specifically tailored (i.e. non-adaptive) Gaussian process prior.
	
	A natural question is how well standard Gaussian process priors perform when the true parameter has a compositional structure. This question can be addressed by proving a lower bound on the contraction rate, that is a sequence $\zeta_n\to0$ such that for a given prior $\Pi$,
	$$ \Pi(\|\theta - \theta^*\|_{L^2}\leq \zeta_n\mid D_n) \overset{P_{\theta^*}}{\to} 0 $$
	as $n\to\infty$ when $D_n\sim P_{\theta^*}^n$. In proving such a result, we may assume that the structure $\eta^*$ is known. In fact, we will assume that $\theta^*$ comes from a \emph{generalised additive model} of the form
	\begin{equation}\label{Eq: generalised additive model}
		\theta^*(x) = F^*\left(g_1(x_1) + \ldots + g_d(x_d)\right), \quad x\in\mathcal{O}, F^*,g_1,\ldots,g_d\in C(\R).
	\end{equation}
	Generalised additive models are a popular and flexible class of models, used frequently in function estimation (see \cite{hastieGeneralizedAdditiveModels1990}). This setting was studied for a regression problem in \cite{giordanoInabilityGaussianProcess2022}, where it was shown that any mean-zero Gaussian process prior based on a wavelet series expansion suffers a severe curse of dimensionality (Theorem 3 of that reference). Moreover (Theorem 1, ibid), any Gaussian process prior has worse performance than that of the DGP prior in \cite{finocchioPosteriorContractionDeep2023}. However, proofs in \cite{giordanoInabilityGaussianProcess2022} hinge on the conjugacy of the Gaussian regression model, and in particular rely on a closed form expression for the posterior mean and variance. In our inverse problem setting, the nonlinear forward map $\G$ results in a non-Gaussian posterior, rendering such an approach impossible. We therefore restrict our attention to proving contraction rate lower bounds for specific Gaussian process priors, using ideas from \cite{castilloLowerBoundsPosterior2008}.
	
	For lower bounds, it suffices to consider the special case of \eqref{Eq: generalised additive model} where $g_1 = \ldots = g_d = \mathrm{id}:\R\to\R$. To reduce technicalities, suppose that $\mathcal{O}\supset [-1,1]^d$ and that $\theta^*$ has known smoothness $\alpha>\beta+d/2$ and is supported in the cube $[-1,1]^d$. Thus we restrict our attention to parameters $\theta^*$ of the form
	\begin{equation}\label{Eq: simple GAM}
		\theta^*(x) = F^*\left(x_1 + \ldots + x_d\right), \quad x\in\mathcal{O}, F^*\in H^{\alpha}(\R),\, \mathrm{supp}(F^*)\subset [-d,d].
	\end{equation}
	
	
	Let $\Pi^{\alpha}$ be the law of an $\alpha$-smooth Whittle-Mat\'{e}rn process on $\mathcal{O}$ multiplied by a smooth cutoff function which is supported inside $\mathcal{O}$ and equals 1 on the cube $[-1,1]^d$. We then define the prior
	\begin{equation}\label{Eq: Darcy rescaled GP prior}
		\tilde{\Pi} = \mathcal{L}\left((\sqrt{n}\delta_n)^{-1}Z\right), \quad Z\sim \Pi^{\alpha}
	\end{equation}
	where the rate $\delta_n$ is defined as
	\begin{equation}\label{Eq: Darcy small ball rate}
		\delta_n = n^{-\frac{\alpha+1}{2\alpha + 2 + d}}.
	\end{equation}
	The upper bounds in \cite{giordanoConsistencyBayesianInference2020} suggest that the posterior induced by the prior $\tilde{\Pi}$ contracts around any $\theta^*$ of the form \eqref{Eq: simple GAM} at a rate $ \delta_n^{(\beta-1)/(\beta+1)}$.
%
%
	Observe that by Theorem \ref{Thm: L2 contraction rate}, the DGP prior $\Pi$ defined in Section \ref{Section: DGP Prior} uniformly attains the contraction rate
	$$( \log{n})^{\bar{\delta}}n^{-\frac{\alpha(\beta-1)}{(2\alpha + 1)(\beta+1)}} \ll \delta_n^{\frac{\beta-1}{\beta+1}} = n^{-\frac{(\alpha+1)(\beta-1)}{(2\alpha + 2 + d)(\beta+1)}} $$
	if $d$ is large. In other words, the DGP prior can leverage the additive structure of $\theta^*$ to achieve a `one-dimensional' rate, whereas the prior $\tilde{\Pi}$ only attains a $d$-dimensional rate. We will prove that this effect is genuine by establishing a contraction rate lower bound for $\tilde{\Pi}$, for a given $\theta^*$ of the form \eqref{Eq: simple GAM}.
	
	We consider a family of rescaled Gaussian process priors, of which the above $\tilde{\Pi}$ is a special instance. For any $\tau>\beta+d/2$, let $\Pi^{\tau}$ denote the law of a $\tau$-smooth Whittle-Mat\'{e}rn process on $\mathcal{O}$ multiplied by a smooth cutoff function equalling 1 on $[-1,1]^d$ as before. Define the rescaled prior
	\begin{equation}\label{Eq: specific rescaled GP prior}
		\tilde{\Pi}^{\tau} = \mathcal{L}\left( n^{-\frac{d}{4\tau+4+2d}} Z \right), \quad Z \sim\Pi^{\tau}. 
	\end{equation}
	The choice of rescaling in the prior is in some sense canonical for modelling a $\tau$-smooth function in this inverse problem (see \cite{giordanoConsistencyBayesianInference2020}). The prior $\tilde{\Pi}$ from \eqref{Eq: Darcy rescaled GP prior} is the case $\tau=\alpha$.
	\begin{theorem}\label{Thm: contraction rate lower bound}
		Let $\G$ be the forward map in Darcy's problem given by \eqref{Eq: forward map}, with $\mathcal{O}\supset[-1,1]^d$. Fix an integer $\beta>1$ and an integer smoothness $\alpha>\beta + d/2$. Let $\tau>\beta+d/2$ be an integer, and consider the prior $\tilde{\Pi}^{\tau}$ defined in \eqref{Eq: specific rescaled GP prior}. Then for any $K>0$, for all $n$ sufficiently large there exists $\theta^*$ of the form \eqref{Eq: simple GAM} with $F^*\in B_{H^{\alpha}(\R)}(K)$, $\mathrm{supp}(F^*)\subset[-d,d]$ such that for some sufficiently small constant $a>0$ (depending on $\tau,\alpha,\beta,d,K$),
		\begin{equation}\label{Eq: crlb}
			E_{\theta^*}\tilde{\Pi}^{\tau}\left(\theta: \|\theta - \theta^*\|_{L^2} \leq a\delta_n^{\frac{\alpha}{\alpha+1}} \mid D_n \right) \to 0.
		\end{equation}
	\end{theorem}
	\begin{remark}
		From the proof of the theorem, it can be seen that this lower bound is only sharp for $\tau=\alpha$: if $\tau\neq\alpha$, one can deduce an even slower rate. Proposition \ref{Prop: specific crlbs} considers a greater variety of rescaling rates in the prior, all of which are subject to the lower bound \eqref{Eq: crlb}.
	\end{remark}
	\begin{remark}\label{Remark: lb for hierarchical Bayes}
		The lower bound also holds for any hierarchical prior $\tilde{\Pi}$ defined by first drawing $\tau$ from a hyperprior with compact support in $(\beta+d/2,\infty)$, and then setting $\tilde{\Pi}\mid\tau = \tilde{\Pi}^{\tau}$: the result of Theorem \ref{Thm: contraction rate lower bound} is not due to the non-adaptivity of the prior.
	\end{remark}
	
	The proof of Theorem \ref{Thm: contraction rate lower bound} can be found in Appendix \ref{Appendix: proofs}. The theorem says that the posterior induced by $\tilde{\Pi}^{\tau}$ asymptotically places almost all of its mass \emph{outside} the $L^2$-ball around $\theta^*$ of radius proportional to $\delta_n^{\frac{\alpha}{\alpha+1}}$. Meanwhile, as previously discussed the DGP prior $\Pi$ induces a posterior which asymptotically puts all of its mass \emph{inside} an $L^2$ neighbourhood of $\theta^*$ of radius
	$$( \log{n})^{\bar{\delta}}n^{-\frac{\alpha(\beta-1)}{(2\alpha + 1)(\beta+1)}} \ll \delta_n^{\frac{\alpha}{\alpha+1}} = n^{-\frac{\alpha}{2\alpha + 2 + d}} $$
	if $d$ is large. Thus in this regime, the DGP prior outperforms the Gaussian process prior $\tilde{\Pi}^{\tau,\rho}$ by a polynomial factor. Note that the Gaussian process prior does not suffer from a curse of dimensionality in the strict sense (i.e. the contraction rate does not become arbitrarily slow as $d\to\infty$) since we must impose minimum smoothness requirements in order to solve the inverse problem. However, the gap between these two rates can be considerable when $d$ is large. 
	
	Intuitively, the contraction rate lower bound arises due to one of two issues: either the prior is simply too rough to concentrate quickly around the truth, or the RKHS of the prior does not suitably approximate the truth. This should be understood as a bias-variance tradeoff: a smoother prior will concentrate faster, but has a smaller RKHS which may approximate $\theta^*$ poorly. When the limiting factor is the quality of RKHS approximation, it is interesting to consider the particular choice of $\theta^*$ for which the lower bound \eqref{Eq: crlb} holds. We choose $F^*$ to be a `spike' with the correct smoothness, and the additive structure in \eqref{Eq: simple GAM} then propagates this spike in all directions, which results in $\theta^*$ having a large number of non-negligible coefficients when expressed in a wavelet basis: see \eqref{Eq: hard to approximate GAM} below. Since the RKHS norm of the Gaussian prior is equivalent to a Sobolev norm which can be characterised in terms of wavelet coefficients, this leads to a fundamental limit to the quality of approximation. However, the DGP prior seems to be able to `learn', or at least exploit, the structural symmetry of such a $\theta^*$, resulting in a `one-dimensional' rate as discussed above.
	
	\section{Further Discussion}\label{Section: discussion}
	
	\subsection{The DGP Prior}\label{Subsection: DGP prior discussion}
	
	Links between Gaussian processes and other deep learning methods, such as deep neural networks and Bayesian neural networks, are drawn throughout much of the machine learning literature. Deep Gaussian processes may be added to this conversation when considering `bottlenecked' deep neural networks. Rather than give a survey here, we refer to Section 7 of \cite{finocchioPosteriorContractionDeep2023}, and we instead discuss our DGP prior in the context of other DGP prior constructions for which there exist theoretical guarantees, namely \cite{bachocPosteriorContractionRates2021} and \cite{finocchioPosteriorContractionDeep2023}.
	
	Like the authors of \cite{finocchioPosteriorContractionDeep2023}, we view our DGP prior more as a proof of concept than an implementable algorithm. In particular, the randomisation over structures $\eta$ incurs a massive computational cost due to the combinatorial explosion of the number of parameters as the depth increases. As discussed in Section 7.1 of \cite{finocchioPosteriorContractionDeep2023}, this effect can be reduced as many structures lead to equivalent contraction rates, and it suffices to consider equivalence classes of structures rather than all possible structures.
	
	One area where our DGP prior improves on that of \cite{finocchioPosteriorContractionDeep2023} is the conditioning step \eqref{Eq: elementary rescaled prior}. In \cite{finocchioPosteriorContractionDeep2023}, the set conditioned on is an $L^{\infty}$-widening of a Sobolev ball; as noted in that paper, it is a challenging computational problem to actually confirm that a draw from a Gaussian process belongs to such a set. Instead, our conditioning set is the intersection of a $L^{\infty}$ ball and a $C^{\beta}$ ball: it is very easy to check that a draw from a Gaussian process satisfies these conditions, and so a simple accept-reject step can be added to perform the conditioning. Moreover, as shown in \eqref{Eq: Borell-TIS on rescaled prior}, our specific choice of Gaussian process $\bar{\Pi}_{\alpha,t}$ means that draws lie in the conditioning set with high probability and so (for sufficiently large $n$) the probability of rejection in this accept-reject algorithm is low. Our conditioning is similar to that used in \cite{bachocPosteriorContractionRates2021}, which considered density estimation and classification problems for compositional parameter spaces with known structure parameter. However, the problem of adapting to the structure $\eta$ is not considered in \cite{bachocPosteriorContractionRates2021}.
	
	\subsection{Posterior Computation}
	
	The posterior arising from the DGP prior is potentially very complex and multimodal, due to the complexity of the prior and the non-concavity of the log-likelihood \eqref{Eq: log-likelihood}. Moreover, computing the posterior itself is computationally intractable, due to the normalising factor $\int e^{\ell_n(\theta)}\,\ud\Pi(\theta)$. A variational Bayes approach is discussed in \cite{finocchioPosteriorContractionDeep2023}, wherein for a fixed structure $\eta$ the posterior $\Pi(\cdot\mid D_n,\eta)$ is approximated by a composition of super-smooth Gaussian processes; one can then sample from the full posterior by first sampling a structure $\eta\sim\pi(\eta)$ where $\pi$ is the structure hyperprior defined in \eqref{Eq: structure hyperprior}, and then using the variational approximation for $\Pi(\cdot\mid D_n,\eta)$.
	
	Alternatively, Markov-chain Monte Carlo (MCMC) methods are commonly used in nonlinear inverse problems to approximate Bayesian posteriors, and if carefully calibrated come with attractive computational guarantees: see Chapter 5 of \cite{nicklBayesianNonlinearStatistical2023} and references therein. However, these results are all for Gaussian priors, for which Gaussian proposal kernels in Metropolis-Hastings algorithms are a natural choice. Moreover, infinite-dimensional Gaussian process priors possess a natural finite-dimensional approximation through truncating their Karhunen-Lo\`{e}ve expansion (see e.g. \cite[Theorem 2.6.10]{gineMathematicalFoundationsInfinitedimensional2016}); however, simply composing these truncations may not lead to a good approximation of a deep Gaussian process. It is therefore not clear what a suitable proposal kernel for the DGP prior could be in such algorithms, even for a fixed structure. These questions are left to future research.

	\subsection{Compositional Structures, Depth and Non-Stationarity}
	
	In the deep learning literature, depth is typically a proxy for `expressivity': the ability of a procedure to reconstruct complicated or irregular functions. For example, in the case of deep neural networks, adding additional layers of neurons enriches the class of functions expressible by the network. However, as shown in \cite[Theorem 4]{dunlopHowDeepAre2018}, repeated composition of Gaussian processes eventually leads to trivial behaviour. Thus in the case of deep Gaussian processes, the role of depth should be thought of somewhat differently.
	
	One way to do this is to consider compositional classes of functions such as $\Theta(\eta^*)$ for $\eta^*\in\Omega$, which were introduced in \cite{schmidt-hieberNonparametricRegressionUsing2020}. Here, the depth $q$ plays a role much the same as any other structure parameter measuring smoothness or dimension. However, the non-identifiability of the compositional representation \eqref{Eq: compositional function} somewhat complicates the proofs, since there is not a `correct' structure around which the (marginal) posterior concentrates as occurs elsewhere in the hierarchical Bayes literature, for example \cite{knapikBayesProceduresAdaptive2016,szaboEmpiricalBayesScaling2013}. The penalisation term in \eqref{Eq: structure hyperprior} suggests that posterior draws should have not too large a depth, and thereby a somewhat simple or efficient structure is typically selected.
	
	An alternative use of deep Gaussian processes has been to generate non-stationary behaviour from covariance kernels, by using a small depth greater than 1. While realisations of Gaussian processes from many widely used covariance kernels (square exponential, Whittle-Mat\'{e}rn) have paths with a global prescribed smoothness, in many applications it is desirable to generate draws which are very regular in some places and more irregular in others. See \cite{sauerNonstationaryGaussianProcess2023} for a survey of such methods, including both non-stationary covariance kernels and deep Gaussian processes. Our analysis applies to this setting insofar as compositional classes model functions with differing degrees of local smoothness. It would be interesting to see if fast contraction rates can be proved for suitably defined classes of functions with variable local smoothness.
		
	\section{Acknowledgements}
	The authors would like to thanks Conor Osborne for many productive conversations, as well as Richard Nickl, Randolf Altmeyer and Kolyan Ray for several suggestions which improved this article.	
	\bibliographystyle{siamplain}
	\bibliography{inverseproblems}
	
\newpage
	\appendix
	
	\section{Proofs}\label{Appendix: proofs} 
	
	\subsection{Proof of Theorem \ref{Thm: prediction risk contraction rate}}\label{Subsection: UB proof}
	
	\subsubsection{Information-Theoretic Distances, Scheme of the Proof}
	
	Here we gather some facts about the relation between information-theoretic distances in this model and the prediction risk distance, and give an overview of the proof of the prediction risk contraction rate result, Theorem \ref{Thm: prediction risk contraction rate}.
	
	Recall the Kullback-Leibler divergence and variation from $P_{\theta_1}$ to $P_{\theta_2}$, defined respectively as
	$$ K(P_{\theta_1},P_{\theta_2}) = E_{\theta_1}\log{\frac{p_{\theta_1}(Y_1,X_1)}{p_{\theta_2}(Y_1,X_1)}}, \quad V(P_{\theta_1},P_{\theta_2}) = E_{\theta_1}\left(\log{\frac{p_{\theta_1}(Y_1,X_1)}{p_{\theta_2}(Y_1,X_1)}}\right)^2.$$
	We also recall the Hellinger distance: given two probability measures $P_{\theta_1},P_{\theta_2}$ on $\R\times\mathcal{O}$ with respective Lebesgue densities $p_{\theta_1},p_{\theta_2}$, this is defined as
	$$ h^2(P_{\theta_1},P_{\theta_2}) = h^2(p_{\theta_1},p_{\theta_2}) = \int_{\R\times\mathcal{O}}\left(\sqrt{p_{\theta_1}(y,x)}-\sqrt{p_{\theta_2}(y,x)}\right)^2\,\ud y\ud x. $$
	By Proposition 1.3.1 in \cite{nicklBayesianNonlinearStatistical2023}, Condition \ref{Cond: G unif boundedness} implies the following inequalities:
	\begin{align}
		K(P_{\theta_1},P_{\theta_2}) &= \frac{1}{2}\|\G(\theta_1) - \G(\theta_2)\|^2_{L^2(\mathcal{O})}, \label{Eq: KL-pred risk equivalence} \\
		V(P_{\theta_1},P_{\theta_2}) &\leq C_1(U)\|\G(\theta_1) - \G(\theta_2)\|^2_{L^2(\mathcal{O})} \label{Eq: KLVar- pred risk bound}, \\
		C_2(U)\|\G(\theta_1) - \G(\theta_2)\|^2_{L^2(\mathcal{O})} &\leq h^2(P_{\theta_1},P_{\theta_2}) \leq \frac{1}{4}\|\G(\theta_1) - \G(\theta_2)\|^2_{L^2(\mathcal{O})}, \label{Eq: Hellinger-pred risk equivalence}
	\end{align}
	where $U$ is the constant from \eqref{Eq: G unif boundedness} and $C_i(U)>0$ are constants depending on $U$ only. We further define the Kullback-Leibler neighbourhood
	\begin{equation}\label{Eq: KL neighbourhood def}
		\mathcal{B}_2(P_{\theta^*},\varepsilon) = \left\{\theta: K(P_{\theta^*},P_{\theta})\leq \varepsilon^2, V(P_{\theta^*},P_{\theta})\leq \varepsilon^2 \right\}.
	\end{equation}
	
	The proof of Theorem \ref{Thm: prediction risk contraction rate} follows standard Bayesian contraction rate ideas and methods, although it is complicated by the use of compositional functions. In particular, we will use a partition entropy argument (see Theorem 8.14 in \cite{ghosalFundamentalsNonparametricBayesian2017}): this is necessitated by the fact that unlike in many settings, the marginal posterior on structures $\eta$ will not concentrate on or close to the true structure parameter $\eta^*$. Instead, various structures $\eta$ induce convergence rates (almost) as fast as $\eta^*$ itself, due to several factors including various forms of redundance in the compositional model and the fact that arbitrarily deep structures can approximate all functions well. However, the penalisation of our structure hyperprior (see \eqref{Eq: structure hyperprior}) forces the posterior to concentrate on `simple' models capable of obtaining fast rates; see Lemma \ref{Lemma: model selection} below. The partition entropy argument then ensures that the posterior concentrates about the true $\theta^*$ in prediction risk on these simple models.
	
	\subsubsection{Small Ball Probability}
	
	As is typical in contraction rate proofs, we first verify a small ball condition for the deep GP prior $\Pi$ with convergence rate $\varepsilon^{\eta^*}_n$ as defined in \eqref{Eq: rates defn}, where $\eta^*$ is the structure parameter (see \eqref{Eq: structure parameter definition}) of the true $\theta^*$. This is a bound of the form
	\begin{equation}\label{Eq: small ball condition}
		\Pi\left(\mathcal{B}_2\left(P_{\theta^*},(\varepsilon_n^{\eta^*})^2\right)\right) \geq a\exp\left(-An(\varepsilon_n^{\eta^*})^2\right),
	\end{equation} 
	for some constants $a,A>0$.
	Note that by \eqref{Eq: KL-pred risk equivalence} and \eqref{Eq: KLVar- pred risk bound}, we have that for any $\varepsilon>0$
	\begin{equation}\label{Eq: KL nbhd contains pred risk ball}
		\mathcal{B}_2(P_{\theta^*},\varepsilon) \supseteq \left\{\theta: \|\G(\theta^*) - \G(\theta)\|_{L^2}\leq c_U\varepsilon \right\} 
	\end{equation}
	for a constant $c_U$ depending only on $U$. So it suffices to check that there exist constants $a,A>0$ such that for all sufficiently large $n$,
	\begin{equation}\label{Eq: STP small ball condition}
		\Pi\left(\theta: \|\G(\theta)-\G(\theta^*)\|_{L^2}\leq c_U\varepsilon_n^{\eta^*}\right) \geq a\exp\left\{-An(\varepsilon_n^{\eta^*})^2\right\}.
	\end{equation}
	We first localise around the true structure $\eta^*$. Given a smoothness $\alpha^*>0$, define the interval $I_n(\alpha^*) = [\alpha^*-1/\log{n},\alpha^*]$, and let $I_n^* = I_n(\boldsymbol{\alpha}^*) = \prod_i I_n(\alpha^*_i)$ be the hypercube of smoothnesses close to $\boldsymbol{\alpha}^*$. For sufficiently large $n$ this interval is contained within the marginal support of $\alpha$ under the hyperprior $\pi$. Some simple algebra shows that for all $\alpha'\in I_n(\alpha^*)$, we have that
	\begin{equation}\label{Eq: equivalence of rates for nearby smoothnesses}
		\varepsilon^{\alpha^*,t}_n \leq \varepsilon^{\alpha',t}_n \leq 3\varepsilon^{\alpha^*,t}_n.
	\end{equation}
	Then
	\begin{align}
		&\Pi\left(\theta: \|\G(\theta)-\G(\theta^*)\|_{L^2}\leq c_U\varepsilon_n^{\eta^*}\right) \nonumber\\
		\geq& \int_{\{\lambda^*\}\times I_n^*}\Pi\left(\theta: \|\G(\theta)-\G(\theta^*)\|_{L^2}\leq c_U\varepsilon_n^{\eta^*}\mid\eta\right)\,\ud\pi(\eta) \nonumber\\
		\gtrsim& e^{-\Psi_n(\eta^*)}\gamma(\lambda^*)\int_{I_n^*}\Pi\left(\theta: \|\G(\theta)-\G(\theta^*)\|_{L^2}\leq c_U\varepsilon_n^{\eta^*}\mid\lambda^*,\boldsymbol{\alpha}\right)\,\ud\gamma(\eta\mid\lambda^*) \label{Eq: lemma smallball stage 1}
	\end{align}
	where the constant is universal, using \eqref{Eq: equivalence of rates for nearby smoothnesses}. In order to develop the integrand, we appeal to the Lipschitz estimate \eqref{Eq: forward Lipschitz condition} together with the following lemma, which establishes that function composition behaves continuously with respect to the supremum norm.
	\begin{lemma}[Lemma 13, \cite{finocchioPosteriorContractionDeep2023}] \label{Lemma: continuity of composition}
		Let $h_{ij}:[-1,1]^{t_i}\to[-1,1]$, and $h_i = (h_{ij})_j$. Assume that for some $M>0$, $h_{ij}\in B_{C^1_{t_i}}(M)$. Then for any $\tilde{h}_{ij}:[-1,1]^{t_i}\to[-1,1]$, $\tilde{h}_i = (\tilde{h}_{ij})_j$, we have that
		$$ \left\| h_q\circ\cdots\circ h_0 - \tilde{h}_q\circ\cdots\circ \tilde{h}_0 \right\|_{\infty} \leq M^q \sum_{i=0}^q \left\| \max_{1\leq j\leq t_{i+1}} |h_{ij} - \tilde{h}_{ij}| \right\|_{\infty}. $$
	\end{lemma}
	The result is not affected by letting $h_q:[-1,1]^{t_q}\to\R$ and $h_0:\mathcal{O}\to[-1,1]^{d_1}$.
	
	We may now return to bounding \eqref{Eq: lemma smallball stage 1}. Fix $\boldsymbol{\alpha}$ for the moment. Note that conditioned on $\lambda^*$, a draw $\theta$ from the prior may be expressed as $\theta = \theta_{q^*}\circ\dots\circ\theta_0$. Also, due to the conditioning step in \eqref{Eq: elementary rescaled prior} and the fact that $\beta\geq1$, for every $i$, $\|\theta_i\|_{C^1}\leq M_0$. Assume that  $M_0\geq 2\max_{i,j}\|\theta^*_{ij}\|_{C^{\beta}}$ (one may choose $M_0<\infty$ since $\eta^*\in\Omega'$, so for all $i$, $\alpha^*_i>\beta + t^*_i/2$, and hence by Sobolev embedding, $\|\theta^*_{ij}\|_{C^{\beta}}<\infty$). Also, conditionally on $\lambda^*$ the $C^{\beta}$-norm of prior draws is bounded by $M_0^{(q^*+1)}$. Since $\beta\geq1$, the $C^{\beta}$-norm dominates the $C^1$-norm. Then by the forward Lipschitz estimate \eqref{Eq: forward Lipschitz condition} and Lemma \ref{Lemma: continuity of composition}, we have that
	\begin{equation}\label{Eq: continuity of composition}
		\|\mathcal{G}(\theta)-\mathcal{G}(\theta^*)\|_{L^2} \lesssim \|\theta - \theta^*\|_{L^{\infty}} \leq M_0^{q^*} \sum_{i=0}^{q^*}\left\| \max_{1\leq j\leq d^*_{i+1}}| \theta_{ij} - \theta^*_{ij}|\right\|_{\infty},
	\end{equation}
	where the first constant depends on $M_0$ and $q^*$. Since the $\theta_{ij}$ are drawn independently under $\Pi$, there is a constant $c_1<\infty$ which depends on $c_U, M_0,q^*$ such that
	\begin{align}\label{Eq: lemma product LB}
		&\Pi\left(\theta: \|\mathcal{G}(\theta)-\mathcal{G}(\theta^*)\|_{L^2} \nonumber \leq c_U\varepsilon_n^{\eta^*}\mid\lambda^*,\boldsymbol{\alpha}\right)\\
		\geq& \prod_{i=0}^{q^*} \prod_{j=1}^{d^*_{i+1}}\Pi\left(\|\theta_{ij} - \theta^*_{ij}\|_{\infty} \leq c_1 \varepsilon_n^{\eta^*} \mid\lambda^*,\boldsymbol{\alpha}\right).
	\end{align}
	We continue to lower bound the term in the product as
	\begin{align*}
		& \Pi_{\alpha_i,t_i^*}\left(\|\theta_{ij} - \theta^*_{ij}\|_{\infty} \leq c_1 \varepsilon_n^{\eta^*}\right) \\
		=& \frac{\bar{\Pi}_{\alpha_i,t_i^*}\left(\|\theta_{ij} - \theta^*_{ij}\|_{\infty} \leq c_1 \varepsilon_n^{\eta^*}, \, \|\theta_{ij}\|_{\infty}\leq 1, \, \|\theta_{ij}\|_{ C^\beta } \leq M_0  \right)}{\bar{\Pi}_{\alpha_i,t_i^*}\left( \|\theta_{ij}\|_{\infty}\leq 1, \, \|\theta_{ij}\|_{ C^\beta} \leq M_0  \right)} \\
		\geq&\bar{\Pi}_{\alpha_i,t_i^*}\left(\|\theta_{ij} - \theta^*_{ij}\|_{\infty} \leq \min\left\{ c_1 \varepsilon_n^{\eta^*}, 1-\|\theta^*_{ij}\|_{\infty} \right\}  , \, \|\theta_{ij}-\theta^*_{ij}\|_{C^1} \leq \frac{M_0}{2} \right)
	\end{align*}
	where we assume that $\|\theta^*_{ij}\|_{\infty}\leq 1-\delta$ for some fixed and known $\delta>0$ (this is not problematic as we can always scale by a constant and just transfer it into the final layer, whose codomain is $\R$; if necessary, we therefore make $K$ a little larger). The argument for the final layer is analogous except there is no restriction that $\|\theta_{ij}\|_{\infty}\leq1$. We have also used that $M_0\geq 2\max_{i,j}\|\theta^*_{ij}\|_{C^{\beta}}$, so that $\|\theta_{ij}\|_{C^1}\leq M_0$ is implied by $\|\theta_{ij}-\theta^*_{ij}\|_{C^1}\leq M_0/2$ (via the triangle inequality). For all sufficiently large $n$, the second term in the minimum exceeds the first and we may therefore lower bound this quantity by
	\begin{align}
			& \bar{\Pi}_{\alpha_i,t_i^*}\left(\|\theta_{ij} - \theta^*_{ij}\|_{\infty} \leq c_1 \varepsilon_n^{\eta^*},  \, \|\theta_{ij}-\theta^*_{ij}\|_{C^1} \leq \frac{M_0}{2} \right) \nonumber \\
			\geq& e^{-\frac{1}{2}n(\varepsilon_n^{\alpha_i,t_i^*})^2\|\theta^*_{ij}\|_{\mathcal{H}(\alpha_i,t_i^*)}^2}\bar{\Pi}_{\alpha_i,t_i^*}\left(\|\theta_{ij}\|_{\infty} \leq c_1 \varepsilon_n^{\eta^*},  \, \|\theta_{ij}\|_{C^1} \leq \frac{M_0}{2} \right) \nonumber \\
			\geq& e^{-\frac{1}{2}n(\varepsilon_n^{\alpha_i,t_i^*})^2\|\theta^*_{ij}\|_{\mathcal{H}(\alpha_i,t_i^*)}^2}\bar{\Pi}_{\alpha_i,t_i^*}\left(\|\theta_{ij}\|_{\infty} \leq c_1 \varepsilon_n^{\eta^*} \right) \bar{\Pi}_{\alpha_i,t_i^*}\left(\|\theta_{ij}\|_{C^1} \leq \frac{M_0}{2} \right) \label{Eq: SBC single component bound}
	\end{align}
	where $\mathcal{H}(\alpha_i,t_i^*)$ is the RKHS of $\Pi'_{\alpha_i,t_i^*}$, whose norm is equivalent to the $H^{\alpha_i}([-1,1]^{t_i^*})$ norm, with universal embedding constants; here, we have used the Cameron-Martin theorem (e.g. Corollary 2.6.18 in \cite{gineMathematicalFoundationsInfinitedimensional2016}) and then the Gaussian correlation inequality (e.g. Theorem 6.2.2 in \cite{nicklBayesianNonlinearStatistical2023}) to establish this lower bound. By \eqref{Eq: Borell-TIS on rescaled prior}, for all $n$ sufficiently large the final probability is at least $1/2$, and so it remains to bound the first probability. 
	
	Theorem 1.2 from \cite{liApproximationMetricEntropy1999} establishes that (in the manner of equation (A15) from \cite{giordanoConsistencyBayesianInference2020})
	\begin{align}
		-\log \bar{\Pi}_{\alpha_i,t_i^*}\left(\|Z\|_{\infty}\leq c_1\varepsilon_n^{\alpha_i,t_i^*}\right) &\simeq \left(\sqrt{n}(\varepsilon_n^{\alpha_i,t_i^*})^2\right)^{-\frac{2t_i^*}{2\alpha_i - t_i^*}} \nonumber \\
		&= n(\varepsilon_n^{\alpha_i,t_i^*})^2, \label{Eq: Li Linde small ball}
	\end{align}
	where the constants (can be chosen to) depend continuously on $\alpha_i$. Plugging this back into \eqref{Eq: SBC single component bound}, we obtain (for all sufficiently large $n$) the lower bound
	$$ \frac{1}{2}\exp\left\{ -c_{ij}n(\varepsilon_n^{\alpha_i,t_i^*})^2 \right\},$$
	for a constant $c_{ij}>0$ depending only on $\alpha_i$ (continuously) and $\|\theta^*_{ij}\|_{H^{\alpha_i}_{t_i^*}}$. In fact, by taking the supremum over $\boldsymbol{\alpha}\in I_n^*$ and using $\|\theta^*_{ij}\|_{H^{\alpha_i}_{t_i^*}}\leq K$, we may choose the constants $c_{ij}$ independent of $\alpha_i$ and $\|\theta^*_{ij}\|_{H^{\alpha_i}_{t_i^*}}$, instead depending only on $K, \boldsymbol{\alpha}^*$ and $c_1$.
	
	Substituting the previous bound into \eqref{Eq: lemma product LB}, one obtains for $\boldsymbol{\alpha}\in I_n^*$ that
	\begin{align}
		\Pi\left(\theta: \|\G(\theta)-\G(\theta^*)\|_{L^2} \nonumber \leq c_U\varepsilon_n^{\eta^*}\mid\lambda^*,\boldsymbol{\alpha}\right)
		\geq& \prod_{i=0}^{q^*} \prod_{j=1}^{d^*_{i+1}}\frac{1}{2}\exp\left\{ -c_{ij}n(\varepsilon_n^{\alpha_i,t_i^*})^2 \right\} \nonumber \\
		\geq&\frac{1}{2}\exp\left\{ -c_2|\mathbf{d}^*|_1n(\varepsilon_n^{\eta^*})^2 \right\}  \label{Eq: local STP small ball}
	\end{align}
	for a constant $c_2$ which depends on $K,\boldsymbol{\alpha}^*,c_1$ only. Here we have used that $\boldsymbol{\alpha}\in I_n^*$ and the inequality \eqref{Eq: equivalence of rates for nearby smoothnesses}. As this lower bound is uniform over $\boldsymbol{\alpha}\in I_n^*$, substituting it into \eqref{Eq: lemma smallball stage 1} yields
	\begin{align}
		\Pi\left(\mathcal{B}_2(P_{\theta^*},\varepsilon_n^{\eta^*})\right) 
		&\gtrsim e^{-\Psi_n(\eta^*)}\gamma(\lambda^*)\gamma(I_n^*\mid\lambda^*)\exp\left\{ -c_2|\mathbf{d}^*|_1n(\varepsilon_n^{\eta^*})^2 \right\} \nonumber \\
		&\gtrsim e^{-\Psi_n(\eta^*)}\gamma(\lambda^*)(\log{n})^{-(q^*+1)}\exp\left\{ -c_2|\mathbf{d}^*|_1n(\varepsilon_n^{\eta^*})^2 \right\} \nonumber \\
		&\gtrsim e^{-\Psi_n(\eta^*)}\gamma(\lambda^*)\exp\left\{ -\underbrace{(c_2+1)|\mathbf{d}^*|_1}_{=:A}n(\varepsilon_n^{\eta^*})^2 \right\}  \label{Eq: small ball final bound}
	\end{align}
	where we have used that $\gamma(\cdot\mid\lambda)$ is the uniform distribution over $[\alpha^-,\alpha^+]^{q+1}$, and that $n(\varepsilon_n^{\eta^*})^2\to\infty$ polynomially fast. The multiplicative constant depends only on the choice of $\gamma$ and $\alpha^+$. This establishes the required small ball condition \eqref{Eq: small ball condition} with constants $A = (c_2+1)|\mathbf{d}^*|_1$ and $a = C(\gamma,\alpha^+)e^{-\Psi_n(\eta^*)}\gamma(\lambda^*)$. Note that $A$ depends only on the parameters which define the class $\Theta(\lambda,\alpha,K)$ (defined in \eqref{Eq: bounded compositional class defn}) and $a$ depends only on these parameters as well as the choice of $\gamma$ and $\alpha^+$, which are part of the definition of the prior.
	
	\subsubsection{Model Selection}
	
	The next stage of the proof is to establish that the posterior concentrates on models which are, in a sense, not too complex. This is done exactly as in \cite{finocchioPosteriorContractionDeep2023}, using the penalisation in \eqref{Eq: structure hyperprior}.
	
	First, we define our set of `simple' models: for any $R>0$, define the set of structures
	\begin{equation}\label{Eq: set of good models}
		\mathcal{M}_n(R) := \left\{\eta: \varepsilon_n^{\eta} \leq R\varepsilon_n^{\eta^*}\right\}\cap \left\{\eta:|\mathbf{d}|_1\leq \log\log n\right\}.
	\end{equation}
	These are models which permit a small ball rate at least as fast as (a constant multiple of) $\varepsilon_n^{\eta^*}$, and whose graphs are not too complicated, in the sense that $|\mathbf{d}|_1$ (which is a measure of how many component processes are required) cannot grow too quickly. We will show that so long as $R$ is chosen sufficiently large, the posterior concentrates on $\mathcal{M}_n(R)$.
	
	The key technical tool is the following, reproduced here for the convenience of the reader.
	\begin{lemma}[Lemma 14, \cite{finocchioPosteriorContractionDeep2023}] \label{Lemma: prior condition for posterior contraction}
		Let $(A_n)$ be a sequence of events and $(a_n)$ be some positive sequence such that $na_n^2\to\infty$. Let $\Pi$ be a generic prior and denote the associated posterior by $\Pi(\cdot\mid D_n)$. Suppose that as $n\to\infty$,
		\begin{equation}\label{Eq: prior condition for posterior contraction}
			e^{2na_n^2}\frac{\Pi(A_n)}{\Pi(\mathcal{B}_2(P_{\theta^*},a_n))} \to 0. 
		\end{equation}
		Then as $n\to\infty$,
		$$ E_{\theta^*} \left[ \Pi(A_n\mid D_n) \right] \to 0.$$
	\end{lemma}
	Write $Z_n = \int_{C(\mathcal{O})} e^{\ell_n(\theta)}\,\ud\Pi(\theta)$. By a slight abuse of notation, for any subset of structures $\mathcal{M}\subset\Omega$ we define
	$$ \Pi(\eta \in \mathcal{M}\mid D_n) := Z_n^{-1} \int_{\mathcal{M}} \int_{C(\mathcal{O})} e^{\ell_n(\theta)}\ud\Pi(\theta\mid\eta)\,\pi(\eta)\,\ud\eta, $$
	which is the contribution to the posterior mass from structures in $\mathcal{M}$.
	\begin{lemma}[Model Selection] \label{Lemma: model selection}
		For $R>0$ chosen sufficiently large depending on $K$ and $\eta^*$ only, we have that
		$$ E_{\theta^*}\left[ \Pi\left(\eta\not\in\mathcal{M}_n(R)\mid D_n\right) \right] = O_{P_{\theta^*}}\left(e^{-cn(\varepsilon_n^{\eta^*})^2}\right) $$
		as $n\to\infty$, where the constant $c>0$ depends only on $K, \eta^*$ and $R$.
		\begin{proof}
			We apply Lemma \ref{Lemma: prior condition for posterior contraction} with $a_n = \varepsilon_n^{\eta^*}$ and $A_n = \mathcal{M}_n^c(R)$, where $R$ is to be chosen below. By \eqref{Eq: small ball final bound}, we see that 
			\begin{equation}\label{Eq: prior condition on bad models}
				e^{2n(\varepsilon_n^{\eta^*})^2}\Pi(\mathcal{B}_2(P_{\theta^*},\varepsilon_n^{\eta^*}))^{-1} \lesssim e^{Cn(\varepsilon_n^{\eta^*})^2}
			\end{equation}
			as $n\to\infty$, for some constant $C>0$. To confirm \eqref{Eq: prior condition for posterior contraction}, it therefore suffices to show that
			$$ \Pi\left(\mathcal{M}_n^c(R)\right) \lesssim e^{-Dn(\varepsilon_n^{\eta^*})^2},$$
			for a constant $D>0$ which we may make as large as desired (through our choice of $R$). 
			
			Recall the penalisation term $\Psi_n(\eta) = n(\varepsilon_n^{\eta})^2 + e^{e^{|\mathbf{d}|_1}}$ defined in \eqref{Eq: structure hyperprior}. If $\eta\in\mathcal{M}_n^c(R)$, then one of the following must be true:
			\begin{itemize}
				\item $\varepsilon_n^{\eta}>R\varepsilon_n^{\eta^*}$, in which case $\Psi_n(\eta)\geq R^2n(\varepsilon_n^{\eta^*})^2$;
				\item $|\mathbf{d}|_1 > \log\log{n}$, in which case
				$$ \Psi_n(\eta) > n > R^2n(\varepsilon_n^{\eta^*})^2 $$
				for all large $n$, for any fixed choice of $R$.
			\end{itemize}
			So for $\eta\in\mathcal{M}_n^c(R)$, we have that $\Psi_n(\eta)\geq R^2n(\varepsilon_n^{\eta^*})^2$. Then we may bound the prior mass of $\mathcal{M}_n^c(R)$ as
			\begin{align*}
				\Pi(\mathcal{M}_n^c(R)) &= \int_{\mathcal{M}_n^c(R)}\pi(\eta)\,\ud\eta \\
				&\leq \int_{\mathcal{M}_n^c(R)}e^{-\Psi_n(\eta)}\gamma(\eta)\,\ud\eta \\
				&\leq e^{-R^2n(\varepsilon_n^{\eta^*})^2}\int_{\mathcal{M}_n^c(R)}\gamma(\eta)\,\ud\eta \\
				&\leq e^{-R^2n(\varepsilon_n^{\eta^*})^2}.
			\end{align*}
			Choosing $R>0$ sufficiently large confirms \eqref{Eq: prior condition for posterior contraction} for a suitable constant $D>0$, and concludes the proof.
		\end{proof}
	\end{lemma}

	\subsubsection{Partition Entropy Argument}
	
	We want to show the first claim of Theorem \ref{Thm: prediction risk contraction rate}, namely that the posterior concentrates on the prediction risk ball
	$$ \mathcal{A}_n := \left\{ \theta: \|\G(\theta) - \G(\theta^*)\|_{L^2}  \leq (\log{n})^{\delta}\varepsilon_n^{\eta^*}\right\}, $$
	for some suitable choice of $\delta>0$ to be specified below. Note that by \eqref{Eq: Hellinger-pred risk equivalence}, $\mathcal{A}_n\supset\tilde{\mathcal{A}}_n$, where $\tilde{\mathcal{A}}_n$ is the Hellinger ball
	$$ \tilde{\mathcal{A}}_n := \left\{ \theta: h(p_{\theta},p_{\theta^*}) \leq C_2(U)^{-1}(\log{n})^{\delta}\varepsilon_n^{\eta^*}\right\}. $$
	By Lemma \ref{Lemma: model selection}, it suffices to prove that
	$$ E_{\theta^*} \left[ \Pi(\tilde{\mathcal{A}}_n^c\cap \mathcal{M}_n(R)\mid D_n )\right] \to 0 $$
	as $n\to\infty$, where $R>0$ is chosen large enough that the lemma holds. 
	
	Write $\mathcal{B}_n = \mathcal{B}_2(P_{\theta^*}, \varepsilon_n^{\eta^*})$ for the Kullback-Leibler neighbourhood for which the small ball condition \eqref{Eq: small ball final bound} holds. Further, define the events
	$$ B_n^* := \left\{ \int \frac{p_{\theta}}{p_{\theta^*}}(D_n)\,\ud\Pi(\theta) \geq \Pi(\mathcal{B}_n)e^{-2n(\varepsilon_n^{\eta^*})^2} \right\}; $$
	by Lemma 8.10 in \cite{ghosalFundamentalsNonparametricBayesian2017}, $P_{\theta^*}(B_n^*) \to 1$ as $n\to\infty$, uniformly in $\theta^*$. We have the bound
	\begin{equation}\label{Eq: contraction rate proof stage 1}
		E_{\theta^*} \left[ \Pi(\mathcal{A}_n^c\cap \mathcal{M}_n(R)\mid D_n )\right] \leq 	P_{\theta^*}((B_n^*)^c)  + E_{\theta^*}\left[ \ind_{B_n^*} \Pi(\mathcal{A}_n^c\cap \mathcal{M}_n(R)\mid D_n )\right];
	\end{equation}
	the first term vanishes as $n\to\infty$ and so it remains to control the second.
	
	We first partition the set of models $\mathcal{M}_n(R)$. Let $\Lambda_n(R) = \{\lambda(\eta) : \eta\in\mathcal{M}_n(R) \}$ be the set of graphs represented in $\mathcal{M}_n(R)$. Given $\lambda\in\Lambda_n(R)$, there exists a vector of smoothnesses $\boldsymbol{\alpha}_{\min}$ such that $(\lambda,\boldsymbol{\alpha}_{\min})\in\mathcal{M}_n(R)$ and, letting $\boldsymbol{\alpha}^+=(\alpha^+,\ldots,\alpha^+)$, such that
	$$ (\lambda,\boldsymbol{\alpha})\in\mathcal{M}_n(R) \quad \Rightarrow \quad \boldsymbol{\alpha}_{\min}\leq \boldsymbol{\alpha} \leq \boldsymbol{\alpha}^+ . $$
	Partition the hyperrectangle $[\boldsymbol{\alpha}_{\min},\boldsymbol{\alpha}^+]$ into hypercubes of side length $1/\log{n}$; call these $A_1(\lambda),\ldots,A_{N(\lambda)}(\lambda)$. Then we can partition $\mathcal{M}_n(R)$ as
	\begin{equation}\label{Eq: partition of good models}
		\mathcal{M}_n(R) =  \bigcup_{\lambda\in\Lambda_n(R)}\bigcup_{k=1}^{N(\lambda)}A_k(\lambda).
	\end{equation}
	Consequently, we can write $ \ind_{B_n^*} \Pi(\mathcal{A}_n^c\cap \mathcal{M}_n(R)\mid D_n )$ as
	$$ \ind_{B_n^*} \frac{\sum_{\lambda\in\Lambda_n(R)}\sum_{k=1}^{N(\lambda)}\int_{\mathcal{A}_n^c\cap A_k(\lambda)}\frac{p_{\theta}}{p_{\theta^*}}(D_n)\,\ud\Pi(\theta) }{\int \frac{p_{\theta}}{p_{\theta^*}}(D_n)\,\ud\Pi(\theta)}. $$
	For each pair $(\lambda, k)$, we introduce a test $\phi_{n,\lambda,k}$, i.e. a measurable function of the data $D_n$ taking values in $[0,1]$. We will specify the tests later. Using that $\phi_{n,\lambda,k} + (1-\phi_{n,\lambda,k}) = 1$, we can upper bound the previous quantity by
	\begin{align}
		&\sum_{\lambda\in\Lambda_n(R)}\sum_{k=1}^{N(\lambda)} \phi_{n,\lambda,k}  \nonumber \\
		&+  \ind_{B_n^*} \frac{\sum_{\lambda\in\Lambda_n(R)}\sum_{k=1}^{N(\lambda)}\int_{A_k(\lambda)}\pi(\lambda,\boldsymbol{\alpha})\int_{\mathcal{A}_n^c}\frac{p_{\theta}}{p_{\theta^*}}(D_n)(1-\phi_{n,\lambda,k})\,\ud\Pi(\theta\mid \lambda,\boldsymbol{\alpha})\ud(\lambda,\boldsymbol{\alpha}) }{\int \frac{p_{\theta}}{p_{\theta^*}}(D_n)\,\ud\Pi(\theta)} \nonumber \\
		=: & \,\, T_1 + T_2. \label{Eq: T_1 + T_2 decomp}
	\end{align}
	Let us set aside $T_1$ for the moment and develop the term $T_2$. Using the definition of the event $B_n^*$ and the small ball condition \eqref{Eq: small ball final bound}, we have that
	\begin{equation}\label{Eq: T_2 bound 1}
		E_{\theta^*}[T_2] \lesssim C(\eta^*)e^{ (A+2)n(\varepsilon_n^{\eta^*})^2}\sum_{\lambda\in\Lambda_n(R)}\sum_{k=1}^{N(\lambda)}\int_{A_k(\lambda)}\pi(\lambda,\boldsymbol{\alpha})\int_{\mathcal{A}_n^c}E_{\theta}[1-\phi_{n,\lambda,k}]\,\ud\Pi(\theta\mid \lambda,\boldsymbol{\alpha})\ud(\lambda,\boldsymbol{\alpha}). 
	\end{equation}
	Next, we introduce sieves $\Theta_{n,\lambda,k}$. Fix $\lambda,A_k(\lambda)$, and let $\boldsymbol{\alpha}$ be the minimal smoothness in $A_k(\lambda)$. Observe that since $A_k(\lambda)$ has side length $(\log{n})^{-1}$, by virtue of \eqref{Eq: equivalence of rates for nearby smoothnesses} any $\boldsymbol{\alpha'}\in A_k(\lambda)$ induces the same rates $\varepsilon_n^{\alpha_i,t_i}$ as $\boldsymbol{\alpha}$ up to universal constants. For any $L_1,L_2>0$, define
	\begin{equation}
		\Theta_{n,\lambda,k}(L_1,L_2) := \left\{ \theta\in\Theta(\lambda,\boldsymbol{\alpha}): \forall i,j, \, \theta_{ij}\in B_{H^{\alpha_i}_{t_i}}\left(L_1\frac{\varepsilon_n^{\eta^*}}{\varepsilon_n^{\alpha_i,t_i}}\right) + B_{L^{\infty}_{t_i}}\left(L_2\varepsilon_n^{\eta^*}\left[\frac{\varepsilon_n^{\alpha_i,t_i}}{\varepsilon_n^{\eta^*}}\right]^{2\alpha_i/t_i}\right), \|\theta_{ij}\|_{C^{\beta}}\leq M_0 \right\}.
	\end{equation}
	
	\begin{lemma}\label{Lemma: sieve properties}
		Fix $\lambda,k$ such that $(\lambda,A_k(\lambda))\subset \mathcal{M}_n(R)$. Given any $C>0$, we may choose $L_1,L_2$ such that
		$$ \Pi\left(\Theta_{n,\lambda,k}(L_1,L_2) \mid \lambda, A_k(\lambda)\right) \geq 1 - \exp\left\{ -Cn(\varepsilon_n^{\eta^*})^2 \right\}. $$
		Moreover, if $L_1,L_2$ are chosen as above, for any $\delta>\log{M_0}$ we have for all sufficiently large $n$ depending on $\alpha^+,\alpha^-,L_1,L_2$ and $\delta$ that
		$$ \log \N(\Theta_{n,\lambda,k}(L_1,L_2),\|\cdot\|_{\infty},(\log{n})^{\delta}\varepsilon_n^{\eta^*}) \leq R^2n(\varepsilon_n^{\eta^*})^2. $$
		\begin{proof}
			For the covering number bound, note that by \eqref{Eq: continuity of composition}, to form an $\varepsilon$-covering of $\Theta_{n,\lambda,k}(L_1,L_2)$, it suffices to cover each component part at a radius of $(qM_0^q)^{-1}\varepsilon$. Thus
			\begin{equation}\label{Eq: product covering number bound}
				\begin{aligned}
					&\N\left(\Theta_{n,\lambda,k}(L_1,L_2),\|\cdot\|_{\infty},(\log{n})^{\delta}\varepsilon_n^{\eta^*}\right)\\
					\leq &\prod_{i,j} \N\left(B_{H^{\alpha_i}_{t_i}}\left(L_1\frac{\varepsilon_n^{\eta^*}}{\varepsilon_n^{\alpha_i,t_i}}\right) + B_{L^{\infty}_{t_i}}\left(L_2\varepsilon_n^{\eta^*}\left[\frac{\varepsilon_n^{\alpha_i,t_i}}{\varepsilon_n^{\eta^*}}\right]^{2\alpha_i/t_i}\right), \|\cdot\|_{\infty},\frac{(\log{n})^{\delta}}{qM_0^q}\varepsilon_n^{\eta^*}\right). 
				\end{aligned}
			\end{equation}
			Note that by definition of $\mathcal{M}_n(R)$,
			$$ \left[\frac{\varepsilon_n^{\alpha_i,t_i}}{\varepsilon_n^{\eta^*}}\right]^{2\alpha_i/t_i} \leq R^{2\alpha^+}; $$
			meanwhile, since $q \leq |\mathbf{d}|_1 \leq \log\log n$ we have that for $\delta>\log{M_0}$,
			$$ \frac{(\log{n})^{\delta}}{qM_0^q} \to \infty $$
			as $n\to\infty$. Hence eventually, the radii of the $L^{\infty}_{t_i}$-balls in \eqref{Eq: product covering number bound} are all less than $\frac{(\log{n})^{\delta}}{2qM_0^q}\varepsilon_n^{\eta^*}$. Thus to control these covering numbers, it suffices to control
			$$ \N\left(B_{H^{\alpha_i}_{t_i}}\left(L_1\frac{\varepsilon_n^{\eta^*}}{\varepsilon_n^{\alpha_i,t_i}}\right) , \|\cdot\|_{\infty},\frac{(\log{n})^{\delta}}{2qM_0^q}\varepsilon_n^{\eta^*}\right). $$
			For this, we appeal to standard covering number bounds for Sobolev spaces, for example Proposition C.7 in \cite{ghosalFundamentalsNonparametricBayesian2017} or Theorem 4.3.36 in \cite{gineMathematicalFoundationsInfinitedimensional2016} (the former is more directly applicable; for the first layer, we must use the analogous bound for the domain $\mathcal{O}$, which also holds: see Chapter 3 in \cite{edmundsFunctionSpacesEntropy1996}): for any $\alpha>t/2$, there exists a constant $C_{\alpha}>0$ depending only on $\alpha$ in a continuous fashion such that
			\begin{equation}\label{Eq: Sobolev-sup norm covering numbers}
				\log\N\left(B_{H^{\alpha}_t}(r),\|\cdot\|_{\infty},\varepsilon \right) \leq C_{\alpha}\left( \frac{r}{\varepsilon} \right)^{\frac{t}{\alpha}}.
			\end{equation}
			
			Thus, again using that $M_0^q \leq (\log{n})^{\log{M_0}}$ and the definition of $\varepsilon_n^{\alpha_i,t_i}$, we have the upper bound
			\begin{align*}
				& \log \N\left(B_{H^{\alpha_i}_{t_i}}\left(L_1\frac{\varepsilon_n^{\eta^*}}{\varepsilon_n^{\alpha_i,t_i}}\right) , \|\cdot\|_{\infty},\frac{(\log{n})^{\delta}}{2qM_0^q}\varepsilon_n^{\eta^*}\right) \\
				\leq & C_{\alpha_i} \left(\frac{L_1qM_0^q}{(\log{n})^{\delta}\varepsilon_n^{\alpha_i,t_i}}\right)^{\frac{t_i}{\alpha_i}} \\
				\leq & C_{\alpha_i}	L_1^{t_i/\alpha_i}(\log\log{n})^{t_i/\alpha_i}(\log{n})^{(\log{M_0}-\delta)t_i/\alpha_i} \times n(\varepsilon_n^{\alpha_i,t_i})^2.
			\end{align*}
			Since $\delta>\log{M_0}$, what precedes the $\times$ tends to 0 as $n\to\infty$, and so for all sufficiently large $n$ (depending on $\alpha^+,\alpha^-$ and the constants in the above inequality) we have that
			$$ \ldots \leq n(\varepsilon_n^{\alpha_i,t_i})^2 \leq R^2n(\varepsilon_n^{\eta^*})^2. $$
			Plugging this bound into \eqref{Eq: product covering number bound}, on the logarithmic level the product becomes a sum; each summand is bounded as above and there are $O(|\mathbf{d}|_1)$ terms, which is at most a multiple of $\log\log{n}$. Hence for $n$ sufficiently large this term is absorbed as before and so 
			$$ \log \N\left(\Theta_{n,\lambda,k}(L_1,L_2),\|\cdot\|_{\infty},(\log{n})^{\delta}\varepsilon_n^{\eta^*}\right) \leq R^2n(\varepsilon_n^{\eta^*})^2, $$
			as required.
			
			Next, we prove the prior probability result. We begin by considering a single component function $\theta_{ij}$. Fix a smoothness $\alpha$ and a dimension $t$. As we work conditionally on $A_k(\lambda)$, WLOG $\alpha$ is the smoothness used to define the sieve; if not, we appeal to \eqref{Eq: equivalence of rates for nearby smoothnesses} and alter the constants $L_1,L_2$ by a universal multiplicative factor if necessary. Note that by the conditioning step in the definition of $\Pi_{\alpha,t}$, we can ignore the $C^{\beta}$-norm condition in the definition of the sieve. Thus we wish to derive a lower bound for
			$$ \Pi_{\alpha,t}\left( B_{H^{\alpha}_{t}}\left(L_1\frac{\varepsilon_n^{\eta^*}}{\varepsilon_n^{\alpha,t}}\right) + B_{L^{\infty}_{t}}\left(L_2\varepsilon_n^{\eta^*}\left[\frac{\varepsilon_n^{\alpha,t}}{\varepsilon_n^{\eta^*}}\right]^{2\alpha/t}\right)  \right). $$
			Note that this set is convex and symmetric; by the Gaussian correlation inequality (e.g. \cite[Theorem 6.2.2]{nicklBayesianNonlinearStatistical2023}), the above probability is therefore bounded below by
			$$ \bar{\Pi}_{\alpha,t}\left( B_{H^{\alpha}_{t}}\left(L_1\frac{\varepsilon_n^{\eta^*}}{\varepsilon_n^{\alpha,t}}\right) + B_{L^{\infty}_{t}}\left(L_2\varepsilon_n^{\eta^*}\left[\frac{\varepsilon_n^{\alpha,t}}{\varepsilon_n^{\eta^*}}\right]^{2\alpha/t}\right)  \right), $$
			as the conditioning set in the definition of $\Pi_{\alpha,t}$ is also convex and symmetric. Recall that $\bar{\Pi}_{\alpha,t}$ is the rescaled version of the process $\Pi'_{\alpha,t}$ and so this probability is equal to
			$$ \Pi'_{\alpha,t}\left( B_{H^{\alpha}_{t}}\left(L_1\sqrt{n}\varepsilon_n^{\eta^*}\right) + B_{L^{\infty}_{t}}\left(L_2\left(\sqrt{n}\varepsilon_n^{\eta^*}\right)^{-\frac{2\alpha-t}{t}}\right)  \right). $$
			Applying Borell's inequality (Proposition 11.19 in \cite{ghosalFundamentalsNonparametricBayesian2017}), we obtain that (multiplying $L_1$ by a universal embedding constant)
			\begin{align}
				&  \Pi'_{\alpha,t}\left( B_{H^{\alpha}_{t}}\left(L_1\sqrt{n}\varepsilon_n^{\eta^*}\right) + B_{L^{\infty}_{t}}\left(L_2\left(\sqrt{n}\varepsilon_n^{\eta^*}\right)^{-\frac{2\alpha-t}{t}}\right)  \right) \nonumber \\
				\geq& \Phi\left( \Phi^{-1}\left(e^{-\varphi_0^{\alpha,t}\left(L_2(\sqrt{n}\varepsilon_n^{\eta^*})^{-(2\alpha-t)/t}\right)}\right) + L_1\sqrt{n}\varepsilon_n^{\eta^*} \right), \label{Eq: Borell inequality for exotic sieve}
			\end{align}
			where $\Phi$ is the standard normal cdf and $\varphi_0^{\alpha,t}$ is the small ball exponent of $\Pi'_{\alpha,t}$, defined by $\Pi'_{\alpha,t}\left(\|Z\|_{\infty}\leq\varepsilon\right) = e^{-\varphi_0^{\alpha,t}(\varepsilon)}$. As in \eqref{Eq: Li Linde small ball}, the covering number bound \eqref{Eq: Sobolev-sup norm covering numbers} together with Theorem 1.2 of \cite{liApproximationMetricEntropy1999} establish that
			$$ \varphi_0^{\alpha,t}(\varepsilon) \lesssim \varepsilon^{-\frac{2t}{2\alpha - t}} $$
			for a constant which can be chosen to depend continuously on $\alpha,t$. Thus the first term in the argument of $\Phi$ can be lower bounded as
			\begin{align*}
				\Phi^{-1}\left(e^{-\varphi_0^{\alpha,t}\left(L_2(\sqrt{n}\varepsilon_n^{\eta^*})^{-(2\alpha-t)/t}\right)}\right) &\geq \Phi^{-1}\left(\exp\left\{- c L_2^{-\frac{2t}{2\alpha - t}}n(\varepsilon_n^{\eta^*})^2\right\} \right) \\
				&\gtrsim  c' L_{2}^{-\frac{t}{2\alpha - t}}\sqrt{n}\varepsilon_n^{\eta^*},
			\end{align*}
			where $c,c'>0$ depend continuously on $\alpha,t$, and the constant in the second inequality is universal, by Lemma K.6 in \cite{ghosalFundamentalsNonparametricBayesian2017}. Thus for any $L_1>0$, by choosing $L_2$ sufficiently large (depending on $L_1,t,\alpha$), the argument of $\Phi$ in \eqref{Eq: Borell inequality for exotic sieve} is at least $(L_1/2)\sqrt{n}\varepsilon_n^{\eta^*}$. To ensure that 
			$$ \Phi\left(\frac{L_1}{2}\sqrt{n}\varepsilon_n^{\eta^*}\right) \geq 1 - \exp\left\{ -Cn(\varepsilon_n^{\eta^*})^2\right\}, $$
			we apply $\Phi^{-1}$ to both sides, and then, using Lemma K.6 of \cite{ghosalFundamentalsNonparametricBayesian2017} it suffices to check that
			$$ \frac{L_1}{2}\sqrt{n}\varepsilon_n^{\eta^*} \geq \frac{1}{2}\sqrt{C}\sqrt{n}\varepsilon_n^{\eta^*}.$$
			Clearly given any $C>0$, it suffices to choose $L_1>\sqrt{C}$ and then choose $L_2$ accordingly. We have established that for any $C>0$, one may choose $L_1,L_2>0$ uniformly over $\alpha$ in an interval of width $(1/\log{n})$ (recall \eqref{Eq: equivalence of rates for nearby smoothnesses}) such that
			\begin{equation}\label{Eq: elementary process sieve inequality}
				\Pi_{\alpha,t}\left( B_{H^{\alpha}_{t}}\left(L_1\frac{\varepsilon_n^{\eta^*}}{\varepsilon_n^{\alpha,t}}\right) + B_{L^{\infty}_{t}}\left(L_2\varepsilon_n^{\eta^*}\left[\frac{\varepsilon_n^{\alpha,t}}{\varepsilon_n^{\eta^*}}\right]^{2\alpha/t}\right)  \right) \geq 1 - \exp\left\{ -Cn(\varepsilon_n^{\eta^*})^2 \right\}.
			\end{equation}
			Having established the first result for a single component process, we return to the compositional process. For any $L_1,L_2$ such that \eqref{Eq: elementary process sieve inequality} holds (recalling that $A_k(\lambda)$ is a hypercube with side-length $(1/\log{n})$), we have that
			\begin{align*}
				\Pi\left(\Theta_{n,\lambda,k}^c(L_1,L_2) \mid \lambda, A_k(\lambda)\right) &\leq \sum_{i=0}^q\sum_{j=1}^{d_{i+1}} e^{-Cn(\varepsilon_n^{\eta^*})^2} \\
				&\leq |\mathbf{d}|_1 e^{-Cn(\varepsilon_n^{\eta^*})^2} \\
				&\leq (\log\log{n}) e^{-Cn(\varepsilon_n^{\eta^*})^2} \\
				&\leq e^{-(C/2)n(\varepsilon_n^{\eta^*})^2}
			\end{align*} 
			for $n$ sufficiently large; here we used the definition of $\mathcal{M}_n(R)$. This concludes the proof of the first bound.
		\end{proof}
	\end{lemma}
	
	Using Lemma \ref{Lemma: sieve properties} in conjunction with \eqref{Eq: T_2 bound 1}, we see that $E_{\theta^*}[T_2]$ is bounded above by
	\begin{align}
		&E_{\theta^*}[T_2] \lesssim C(\eta^*)e^{ (A+2)n(\varepsilon_n^{\eta^*})^2}\sum_{\lambda\in\Lambda_n(R)}\sum_{k=1}^{N(\lambda)}\int_{A_k(\lambda)}\pi(\lambda,\boldsymbol{\alpha})\int_{\mathcal{A}_n^c}E_{\theta}[1-\phi_{n,\lambda,k}]\,\ud\Pi(\theta\mid \lambda,\boldsymbol{\alpha})\ud(\lambda,\boldsymbol{\alpha}) \nonumber \\
		\leq& C(\eta^*)e^{ (c_2+2)|\mathbf{d}^*|_1n(\varepsilon_n^{\eta^*})^2}\sum_{\lambda\in\Lambda_n(R)}\sum_{k=1}^{N(\lambda)}\int_{A_k(\lambda)}\pi(\lambda,\boldsymbol{\alpha})\int_{\mathcal{A}_n^c\cap\Theta_{n,\lambda,k}(L_1,L_2)}E_{\theta}[1-\phi_{n,\lambda,k}]\,\ud\Pi(\theta\mid \lambda,\boldsymbol{\alpha})\ud(\lambda,\boldsymbol{\alpha}) \nonumber \\
		&+ C(\eta^*)e^{ (c_2+2)|\mathbf{d}^*|_1n(\varepsilon_n^{\eta^*})^2}\sum_{\lambda\in\Lambda_n(R)}\sum_{k=1}^{N(\lambda)}\int_{A_k(\lambda)}\pi(\lambda,\boldsymbol{\alpha})e^{-Cn(\varepsilon_n^{\eta^*})^2}\ud(\lambda,\boldsymbol{\alpha}),  \label{Eq: T_2 bound 2}
	\end{align}
	using that $E_{\theta}(1-\phi_{n,\lambda,k}) \leq 1$. Note that in Lemma \ref{Lemma: sieve properties}, $L_1$ is chosen depending on $C$ only and $L_2$ is chosen depending on $L_1$ and $\alpha$ only; thus we may choose $L_1,L_2$ such that for fixed but arbitrarily large $C>0$, the lemma holds uniformly over $\eta\in\mathcal{M}_n(R)$ (by using $\alpha^+$ in our choice of $L_2$). Therefore the second term is bounded by
	$$ C(\eta^*)e^{ (A+2)n(\varepsilon_n^{\eta^*})^2}e^{-Cn(\varepsilon_n^{\eta^*})^2} \to 0 $$
	as $n\to\infty$, assuming $C>0$ is chosen sufficiently large depending on $A$.
	
	To bound $E_{\theta^*}[T_1]$ and the remaining term in \eqref{Eq: T_2 bound 2}, we choose specific tests $\phi_{n,\lambda,k}$. For any $\lambda, k$, Theorem D.5 of \cite{ghosalFundamentalsNonparametricBayesian2017} gives a test $\phi_{n,\lambda,k}$ such that for some universal constant $\tilde{D}>0$, we have that
	\begin{align}
		E_{\theta^*}[\phi_{n,\lambda,k}] &\leq c_{\lambda,k}\N\left(\Theta_{n,\lambda,k}(L_1,L_2), \|\cdot\|_{\infty},   (\log{n})^{\delta}\varepsilon_n^{\eta^*}\right) \frac{e^{-4\tilde{D}(\log{n})^{2\delta}n(\varepsilon_n^{\eta^*})^2}}{1 - e^{-4\tilde{D}(\log{n})^{2\delta}n(\varepsilon_n^{\eta^*})^2}},  \label{Eq: tests type I error}\\
		E_{\theta}[1-\phi_{n,\lambda,k}] &\leq c_{\lambda,k}^{-1} e^{-4\tilde{D}
			(\log{n})^{2\delta}n(\varepsilon_n^{\eta^*})^2} \quad \forall \theta\in\Theta_{n,\lambda,k}(L_1,L_2)\cap\mathcal{A}_n^c, \label{Eq: tests type II error}
	\end{align}
	where we choose the constants $c_{\lambda,k}$ such that
	$$c_{\lambda,k}^2 := \frac{\pi(\lambda, A_k(\lambda))}{\N\left(\Theta_{n,\lambda,k}(L_1,L_2), \|\cdot\|_{\infty}, (\log{n})^{\delta}\varepsilon_n^{\eta^*}\right)}. $$
	Define the `local complexities'
	$$ \psi_{\lambda,k} := \sqrt{\pi(\lambda, A_k(\lambda))}\sqrt{\N\left(\Theta_{n,\lambda,k}(L_1,L_2), \|\cdot\|_{\infty}, (\log{n})^{\delta}\varepsilon_n^{\eta^*}\right)}. $$
	Then it is clear from \eqref{Eq: tests type I error} and \eqref{Eq: tests type II error} that
	$$ E_{\theta^*}[T_1] \leq \frac{e^{-4\tilde{D}(\log{n})^{2\delta}n(\varepsilon_n^{\eta^*})^2}}{1 - e^{-4\tilde{D}(\log{n})^{2\delta}n(\varepsilon_n^{\eta^*})^2}}\sum_{\lambda\in\Lambda_n(R)}\sum_{k=1}^{N(\lambda)}\psi_{\lambda,k} $$
	and, using \eqref{Eq: T_2 bound 2}, that
	$$ E_{\theta^*}[T_2] \leq C(\eta^*)e^{ (A+2)n(\varepsilon_n^{\eta^*})^2} e^{-4\tilde{D}(\log{n})^{2\delta}n(\varepsilon_n^{\eta^*})^2}\sum_{\lambda\in\Lambda_n(R)}\sum_{k=1}^{N(\lambda)}\psi_{\lambda,k}  +  o(1). $$
	To complete the proof, it suffices to establish that
	\begin{equation}\label{Eq: local complexity summability condition}
		\sum_{\lambda\in\Lambda_n(R)}\sum_{k=1}^{N(\lambda)}\psi_{\lambda,k} \lesssim 	e^{c_3n(\varepsilon_n^{\eta^*})^2}
	\end{equation}
	for some constant $c_3>0$, since then in view of the previous two displays, we have that $E_{\theta^*}[T_1+T_2] \to 0$ as $n\to\infty$. Now by Lemma \ref{Lemma: sieve properties}, for $L_1,L_2$ as chosen previously, uniformly over $\lambda,k$ we have that (for sufficiently large $n$)
	$$ \sqrt{\N\left(\Theta_{n,\lambda,k}(L_1,L_2), \|\cdot\|_{\infty}, (\log{n})^{\delta}\varepsilon_n^{\eta^*}\right)} \leq \exp\left\{ \frac{1}{2}R^2n(\varepsilon_n^{\eta^*})^2 \right\}, $$
	so to check \eqref{Eq: local complexity summability condition} it suffices to check (for a different but still arbitrary constant $c'_3$) that
	$$ \sum_{\lambda\in\Lambda_n(R)}\sum_{k=1}^{N(\lambda)}\sqrt{\pi(\lambda,A_k(\lambda))}\lesssim e^{c'_3n(\varepsilon_n^{\eta^*})^2}. $$
	
	Now, for $(\lambda,\boldsymbol{\alpha})\in\mathcal{M}_n(R)$, we have that
	$$ \pi(\lambda,\boldsymbol{\alpha}) = z_n^{-1}e^{-\Psi_n(\lambda,\boldsymbol{\alpha})}\gamma(\lambda,\boldsymbol{\alpha}) \leq z_n^{-1}\gamma(\lambda,\boldsymbol{\alpha}), $$
	where $z_n = \int_{\Omega}\pi(\eta)\,\ud\eta$. Recall $I_n^*$, the hypercube of smoothnesses close to $\boldsymbol{\alpha}^*$, which has side-length $(1/\log{n})$. Then using \eqref{Eq: equivalence of rates for nearby smoothnesses} and the fact that $\gamma(\cdot\mid \lambda)$ is uniform,
	\begin{align*}
		z_n &\geq e^{-\Psi_n(\eta^*)}\int_{\{\lambda^*\}\times I_n^*} \gamma(\lambda,\boldsymbol{\alpha})\, \ud(\lambda,\boldsymbol{\alpha}) \\
		&\gtrsim e^{-n(\varepsilon_n^{\eta^*})^2} \gamma(\lambda^*)(\log{n})^{-(q^*+1)} \\
		&\gtrsim e^{-c_4n(\varepsilon_n^{\eta^*})^2} 
	\end{align*}
	for some $c_4>0$ and a multiplicative constant only depending on $\eta^*$. Thus again using that $\gamma(\cdot\mid\lambda)$ is uniform, and letting $\boldsymbol{\alpha}_{k,\lambda}$ be any smoothness in $A_k(\lambda)$, we have that
	\begin{align*}
		\sum_{\lambda\in\Lambda_n(R)}\sum_{k=1}^{N(\lambda)}\sqrt{\pi(\lambda,A_k(\lambda))} &\lesssim e^{\frac{1}{2}c_4n(\varepsilon_n^{\eta^*})^2}\sum_{\lambda\in\Lambda_n(R)}\sum_{k=1}^{N(\lambda)}\sqrt{\gamma(\lambda,A_k(\lambda))} \\
		&= e^{\frac{1}{2}c_4n(\varepsilon_n^{\eta^*})^2}\sum_{\lambda\in\Lambda_n(R)}\sum_{k=1}^{N(\lambda)}\sqrt{|A_k(\lambda)|}\sqrt{\gamma(\lambda,\boldsymbol{\alpha}_{k,\lambda})} \\
		&= e^{\frac{1}{2}c_4n(\varepsilon_n^{\eta^*})^2}\sum_{\lambda\in\Lambda_n(R)}\sum_{k=1}^{N(\lambda)}\frac{1}{\sqrt{|A_k(\lambda)|}}|A_k(\lambda)|\sqrt{\gamma(\lambda,\boldsymbol{\alpha}_{k,\lambda})} 
	\end{align*}
	where $\boldsymbol{\alpha}_{k,\lambda}$ is any smoothness in $A_k(\lambda)$; this holds since $\gamma(\cdot\mid\lambda)$ is uniform. Since $A_k(\lambda)$ is a hypercube of side length $(\log{n})^{-1}$ and a subset of $\mathcal{M}_n(R)$, we have that
	$$ |A_k(\lambda)|^{-1} = (\log{n})^q \leq e^{(\log\log{n})^2} \ll e^{c_4n(\varepsilon_n^{\eta^*})^2}, $$
	and combining this with the previous bound gives
	$$ \sum_{\lambda\in\Lambda_n(R)}\sum_{k=1}^{N(\lambda)}\sqrt{\pi(\lambda,A_k(\lambda))}  \leq e^{c_4n(\varepsilon_n^{\eta^*})^2}\sum_{\lambda\in\Lambda_n(R)}\sum_{k=1}^{N(\lambda)}|A_k(\lambda)|\sqrt{\gamma(\lambda,\boldsymbol{\alpha}_{k,\lambda})} . $$
	Finally, observe that again using that $\gamma(\cdot\mid\lambda)$ is the uniform distribution, we see that
	$$ \sum_{\lambda\in\Lambda_n(R)}\sum_{k=1}^{N(\lambda)}|A_k(\lambda)|\sqrt{\gamma(\lambda,\boldsymbol{\alpha}_{k,\lambda})} \leq \sum_{\lambda\in\Lambda_n(R)}\int_{[\alpha^-,\alpha^+]^q} \sqrt{\gamma(\lambda,\boldsymbol{\alpha})}\,\ud(\lambda,\boldsymbol{\alpha}) \leq \int_{\Omega}\sqrt{\gamma(\eta)}\,\ud\eta <\infty $$
	by Assumption \ref{Assumption 1}. Thus we have established \eqref{Eq: local complexity summability condition}. Constants in all of the above inequalities only depended on $\theta^*$ through the parameters of the class $\Theta_{\mathcal{K}}(\eta^*,K)$ (defined at the start of Section \ref{Section: contraction rates}), so the result holds uniformly over this set. In summary, we have proved that
	$$ \sup_{\theta^*\in \Theta_{\mathcal{K}}(\eta^*,K)}E_{\theta^*} \Pi\left( \theta: \|\G(\theta) - \G(\theta^*)\|_{L^2(\mathcal{O})}\geq (\log{n})^{\delta}\varepsilon_n^{\eta^*}, \eta\not\in\mathcal{M}_n(R) \mid D_n\right) \to 0,$$
	which implies the first part of Theorem \ref{Thm: prediction risk contraction rate}.
	
	For the second part of the theorem, note that for $\eta\in\mathcal{M}_n(R)$, we have that the depth $q$ is also bounded by $\log\log{n}$. Draws $\theta$ from the posterior $\Pi(\cdot\mid D_n)$ therefore satisfy $\|\theta_{ij}\|_{C^{\beta}}\leq M_0$ for all $i,j$ (since prior draws have this property almost surely) and, with high probability, $q(\theta)\leq\log\log{n}$. By the chain rule, for integer $\beta$ this implies that $\|\theta\|_{C^{\beta}}\leq M_0^{\log\log{n}}$. Thus $\|\theta\|_{C^{\beta}}\leq (\log{n})^{\delta}$ with high probability, since $\delta>\log{M_0}$. \qed

	\subsection{Proof of Theorem \ref{Thm: contraction rate lower bound}}\label{Subsection: LB proofs}
	
	We consider a more general family of rescaled Gaussian priors than introduced in Section \ref{Section: crlbs}. Let $\rho = (\rho_n)$ be a sequence such that $\rho_n\to\infty$ and let $\tau>\beta+d/2$. Let $\Pi^{\tau}$ denote the law of a $\tau$-smooth Whittle-Mat\'{e}rn process on $\mathcal{O}$ multiplied by a smooth cutoff function equalling 1 on $[-1,1]^d$ as before. Given such $\rho = (\rho_n),\tau$, define the prior
	\begin{equation}\label{Eq: general rescaled GP prior}
		\tilde{\Pi}^{\tau,\rho} = \mathcal{L}\left(\rho_n^{-1}Z\right), \quad Z\sim\Pi^{\tau}.
	\end{equation}
	We insist on the condition $\rho_n\to\infty$ as this ensures that the prior $\tilde{\Pi}^{\tau,\rho}$ concentrates on the regularisation set $B_{C^{\beta}_d}(M)$ (see \eqref{Eq: Borell-TIS on rescaled prior} above), which is required to apply the stability estimate \eqref{Eq: stability estimate} and thereby achieve consistent reconstruction of $\theta^*$. The prior $\tilde{\Pi}^{\tau}$ from \eqref{Eq: specific rescaled GP prior} is equal to $\tilde{\Pi}^{\tau,\rho}$ for the choice $\rho_n = n^{\frac{d}{4\tau+4+2d}}$; thus Theorem \ref{Thm: contraction rate lower bound} is an immediate consequence of the following result.
	\begin{proposition}\label{Prop: specific crlbs}
		Let $\mathcal{O},\mathcal{G},\beta,\alpha$ be as in Theorem \ref{Thm: contraction rate lower bound}, and fix $K>0$. Let $\tau>\beta+d/2$ be an integer and $\rho_n\to\infty$; if $\tau\leq\alpha$, assume that $\rho_n\lesssim n^{\frac{d}{4\tau+4+2d}}$. Then for all $n$ sufficiently large, there exists $\theta^*$ of the form \eqref{Eq: generalised additive model} with $F^*\in H^{\alpha}(\R)$ such that $F^*$ is supported in $[-d,d]$ and $\|F^*\|_{H^{\alpha}(\R)}\leq K$ for which the contraction rate lower bound
		$$ E_{\theta^*}\tilde{\Pi}^{\tau,\rho}\left(\theta: \|\theta - \theta^*\|_{L^2} \leq a\zeta_n \mid D_n \right) \to 0 $$
		holds for the following rates $\zeta_n$ with a suitable choice of $a>0$:
		\begin{enumerate}[(i)]
			\item If $\tau\leq\alpha$, then $\zeta_n = n^{-\frac{\tau}{2\tau+2+d}}$;
			\item If $\alpha<\tau<\alpha+\frac{d}{2}$, then $\zeta_n = n^{-\frac{\alpha}{2\alpha+2+d}}$;
			\item if $\tau>\alpha+\frac{d}{2}$, then $\zeta_n = (\rho_n)^{\frac{\alpha}{\tau+1}} n^{-\frac{\alpha+1}{2\tau+2}}$.
		\end{enumerate}
		These rates result from choosing $\rho_n\to\infty$ to make $\zeta_n$ as fast as possible; sub-optimal choices of $\rho_n$ result in an even slower rate $\zeta_n$.
	\end{proposition}
	
	The key technical result is Theorem 1 of \cite{castilloLowerBoundsPosterior2008}, which we state here adapted to our setting for the convenience of the reader. Recall that the ($L^2$-)\emph{concentration function} at $\theta^*$ of the Gaussian process $Z$ with RKHS $\mathbb{H}$ is defined as
	\begin{equation}\label{Eq: concentration function definition}
		\varphi_{\theta^*}(\varepsilon) = \varphi^*(\varepsilon) = \inf_{h\in\mathbb{H},\|h - \theta^*\|_{L^2}\leq\varepsilon} \frac{1}{2}\|h\|_{\mathbb{H}}^2 - \log\Pr(\|Z\|_{L^2}\leq \varepsilon).
	\end{equation}
	\begin{proposition}\label{Prop: general contraction rate LB}
		Assume that for some $\theta^*\in C(\mathcal{O})$, $D_n\sim P_{\theta^*}^n$. Let $\Pi$ be the law of a Gaussian process supported on $C(\mathcal{O})$. Let $r_n\to0$ be a sequence such that $nr_n^2\to\infty$ and
		\begin{equation}\label{Eq: crlb small ball condition}
			\Pi\left(\mathcal{B}_2(P_{\theta^*},r_n)\right) \geq \exp\left(-cnr_n^2\right)
		\end{equation}
		for some constant $c>0$, where $\mathcal{B}_2(P_{\theta^*},r_n)$ is the Kullback-Leibler neighbourhood defined in \eqref{Eq: KL neighbourhood def}.
		
		Let $\varphi^*$ be the concentration function of $\Pi$ at $\theta^*$. Suppose that $\zeta_n\to0$ is such that
		$$ \varphi^*(\zeta_n) \geq (c+2)nr_n^2. $$
		Then
		$$ E_{\theta^*} \Pi(\|\theta - \theta^*\|_{L^2} \leq \zeta_n\mid D_n) \to 0. $$
	\end{proposition}
	
	To get the best possible lower bound, the rate $r_n$ in \eqref{Eq: crlb small ball condition} should be as fast as possible. To find a suitable sequence $\zeta_n$ it then suffices to lower bound either of the two terms in \eqref{Eq: concentration function definition} (with $\varepsilon=\zeta_n$), since both are non-negative.
	
	Consider the prior $\tilde{\Pi}^{\tau,\rho}$, which is based on the prior $\Pi^{\tau}$ whose RKHS is
	\begin{equation}\label{Eq: rescaled prior RKHS}
		\mathbb{H} = \left\{ \chi f: f\in H^{\tau}(\mathcal{O})\right\},
	\end{equation}
	where $\chi$ is the cutoff function used in the definition of $\Pi^{\tau}$ (see before \eqref{Eq: general rescaled GP prior}). For any $h\in\mathbb{H}$, there exists $f\in H^{\tau}(\mathcal{O})$ such that $h = \chi f$ and the RKHS norm satisfies
	$$ \|h\|_{\mathbb{H}} = \|f\|_{H^{\tau}(\mathcal{O})}. $$
	As a consequence, for $h\in\mathbb{H}$, we have that$ \|h\|_{H^{\tau}(\mathcal{O})}\lesssim\|h\|_{\mathbb{H}}$ (see Example 25 in \cite{giordanoConsistencyBayesianInference2020}). Also, since $\chi=1$ on $[-1,1]^d$,  (identifying $h$ with its restriction to the cube)  we have that $\|h\|_{H^\tau_d}\leq \|h\|_{\mathbb{H}}$. The RKHS of the rescaled version $\tilde{\Pi}^{\tau,\rho}$ is equal to $\mathbb{H}$ as a set but now the norm is rescaled by a factor of $\rho_n$. We will use these facts throughout the proof.
	
	The next lemma establishes the best possible rates $r_n$ for the prior $\tilde{\Pi}^{\tau,\rho}$, which should in turn give rise to the slowest possible contraction rate lower bounds $\zeta_n$.
	\begin{lemma}\label{Lemma: rescaled GP prior small ball rates}
		Let $\tilde{\Pi}^{\tau,\rho}$ be as in \eqref{Eq: general rescaled GP prior}, for $\tau>\beta+d/2$ and any sequence $\rho_n\to\infty$.
		Let $\theta^*\in H^{\alpha}(\mathcal{O})$ be supported in $[-1,1]^d$, where $\alpha>\beta+d/2$. Then for the choices of $\rho_n$ given below, the small ball condition \eqref{Eq: crlb small ball condition} is satisfied for the following rates $r_n$:
		\begin{enumerate}[(i)]
			\item $\tau\leq \alpha$: $ \rho_n \simeq n^{\frac{d}{4\tau + 4 + 2d}}, r_n \simeq n^{-\frac{\tau+1}{2\tau+2+d}}$;
			\item $\alpha<\tau<\alpha + \frac{d}{2}$: $\rho_n\simeq n^{\frac{d-2(\tau-\alpha)}{4\alpha+4+2d}}, r_n \simeq n^{-\frac{\alpha+1}{2\alpha+2+d}}$;
			\item $\tau\geq\alpha + \frac{d}{2}$: for any sequence $\rho_n\to\infty$, $r_n \simeq \rho_n^{\frac{\alpha+1}{\tau+1}}n^{-\frac{\alpha+1}{2\tau + 2}}$.
		\end{enumerate}
		\begin{remark}
			Observe that in case (iii), since $\tau\geq\alpha+\frac{d}{2}$ the rate is strictly slower than the rate in case (ii). One should think of $\rho_n\to\infty$ very slowly in this instance.
		\end{remark}
		\begin{proof}
			Using the Lipschitz estimate (23) in \cite{giordanoConsistencyBayesianInference2020} for $\mathcal{G}$ (which uses the weak Sobolev norm $\|\cdot\|_{(H^1)^*}$, the topological dual norm of $H^1_c(\mathcal{O})$, in place of the supremum norm in our general Lipschitz estimate \eqref{Eq: forward Lipschitz condition}, see Remark \ref{Remark: can't use Sobolev}), and \eqref{Eq: KL nbhd contains pred risk ball}, it suffices to show that
			$$ \tilde{\Pi}^{\tau,\rho}\left(\theta: \|\theta - \theta^*\|_{(H^1)^*}\leq cr_n, \|\theta\|_{C^{\beta}}\leq M , \|\theta^*\|_{C^{\beta}}\leq M\right) \gtrsim e^{-Anr_n^2} $$
			for a suitable rate $r_n$ and constants $c,M,A>0$. If we choose $M$ such that $\|\theta^*\|_{C^{\beta}}\leq M$, then note that the condition $\|\theta-\theta^*\|_{C^{\beta}}$ implies that $\|\theta\|_{C^{\beta}}\leq 2M$ by the triangle inequality, and we can apply the Lipschitz estimate. So it then suffices to show that
			$$ \tilde{\Pi}^{\tau,\rho}\left(\theta: \|\theta - \theta^*\|_{(H^1)^*}\leq cr_n, \|\theta - \theta^*\|_{C^{\beta}}\leq M \right) \gtrsim e^{-Anr_n^2}. $$
			Letting $\mathbb{H}^{\tau,\rho}$ be the RKHS of $\tilde{\Pi}^{\tau,\rho}$, for any $h\in \mathbb{H}^{\tau,\rho}$ with $\|h-\theta^*\|_{(H^1)^*}\leq \frac{cr_n}{2},\|h-\theta^*\|_{C^{\beta}}\leq \frac{M}{2}$ we have by the triangle inequality that
			\begin{align}
				&\tilde{\Pi}^{\tau,\rho}\left(\theta: \|\theta - \theta^*\|_{(H^1)^*}\leq cr_n, \|\theta-\theta^*\|_{C^{\beta}}\leq M \right) \nonumber \\
				\geq & \tilde{\Pi}^{\tau,\rho}\left(\theta: \|\theta - h\|_{(H^1)^*}\leq \frac{c}{2}r_n, \|\theta-h\|_{C^{\beta}}\leq \frac{M}{2} \right)  \nonumber \\
				\geq & e^{-\frac{1}{2}\|h\|_{\mathbb{H}^{\tau,\rho}}^2} \tilde{\Pi}^{\tau,\rho}\left(\|\theta \|_{(H^1)^*}\leq \frac{c}{2}r_n, \|\theta\|_{C^{\beta}}\leq \frac{M}{2} \right) \nonumber \\
				\geq & e^{-\frac{1}{2}\|h\|^2_{\mathbb{H}^{\tau,\rho}}} \tilde{\Pi}^{\tau,\rho}\left(\|\theta \|_{(H^1)^*}\leq \frac{c}{2}r_n \right)\tilde{\Pi}^{\tau,\rho}\left( \|\theta\|_{C^{\beta}} \leq \frac{M}{2} \right), \label{Eq: small ball rates decomposition 1}
			\end{align}
			using first the Cameron-Martin theorem (e.g. \cite[Corollary 2.6.18]{gineMathematicalFoundationsInfinitedimensional2016}) and then the Gaussian correlation inequality (e.g. \cite[Theorem 6.2.2]{nicklBayesianNonlinearStatistical2023}). As in \eqref{Eq: Borell-TIS on rescaled prior}, since $\rho_n\to\infty$ we have that for $M>0$ chosen sufficiently large the final term tends to 1 as $n\to\infty$. Thus to check \eqref{Eq: crlb small ball condition} it suffices to find a rate $r_n$ such that
			\begin{equation}\label{Eq: small ball concentration function}
				\inf_{\substack{h\in \mathbb{H}^{\tau,\rho}, \|h-\theta^*\|_{(H^1)^*}\leq \frac{cr_n}{2}, \\ \|h-\theta^*\|_{C^{\beta}}\leq \frac{M}{2}}} \left( \frac{1}{2}\|h\|_{\mathbb{H}^{\tau,\rho}}^2 \right) -\log \tilde{\Pi}^{\tau,\rho}\left(\|Z\|_{(H^1)^*}\leq \frac{c}{2}r_n\right) \lesssim nr_n^2.
			\end{equation}
			By Theorem 1.2 in \cite{liApproximationMetricEntropy1999} and equation (A14) in \cite{giordanoConsistencyBayesianInference2020}, we have that
			\begin{equation}\label{Eq: small ball exponent UB}
				-\log \tilde{\Pi}^{\tau,\rho}\left(\|Z\|_{(H^1)^*}\leq \frac{c}{2}r_n\right) \lesssim (\rho_nr_n)^{-\frac{2d}{2\tau+ 2 -d}}.
			\end{equation}
			Observe that we may upper bound the first term in \eqref{Eq: small ball concentration function} by choosing any valid $h$ instead of taking the infimum. We divide into subcases where $\tau\leq\alpha$ and $\tau>\alpha$. When $\tau\leq\alpha$, $\theta^*\in H^{\tau}_d\subset\mathbb{H}^{\tau,\rho}$ and so by choosing $h=\theta^*$, the first term on the left-hand side of \eqref{Eq: small ball concentration function} is bounded by $\frac{1}{2}\rho_n^2\|\theta^*\|_{H^{\alpha}}$. To achieve \eqref{Eq: small ball concentration function}, it suffices to choose $\rho_n,r_n$ such that
			\begin{equation}\label{Eq: small ball rate ineq, undersmooth prior}
				\rho_n^2 \lesssim nr_n^2, \quad (\rho_nr_n)^{-\frac{2d}{2\tau+ 2 -d}} \lesssim nr_n^2. 
			\end{equation}
			Recall that we want to choose $r_n$ as small as possible; thus from the second inequality above, we should choose $\rho_n$ as large as possible. By the first inequality, the correct choice is $\rho_n \simeq \sqrt{n}r_n$. Solving the second inequality then leads to
			$$ r_n = n^{-\frac{\tau+1}{2\tau+2+d}}, \quad \rho_n = n^{\frac{d}{4\tau + 4 + 2d}}, $$
			as in the statement of the lemma.
			
			When $\tau>\alpha$, $\theta^*$ no longer lies in the RKHS and the approximation term becomes significant as we can no longer just choose $h=\theta^*$ to control it. To approximate $\theta^*$, we introduce the \emph{spectral Sobolev spaces}; see Section 6.1.3 of \cite{nicklBayesianNonlinearStatistical2023}. Let $(\lambda_j,e_j)_{j\geq1}$ be the eigenpairs of the Laplacian $-\Delta$ on $\mathcal{O}$; we have the Weyl asymptotics $\lambda_j\simeq j^{2/d}$, and that $(e_j)_{j\geq1}$ is an orthonormal basis for $L^2(\mathcal{O})$. Then we define
			$$ \tilde{H}^s(\mathcal{O}) = \left\{ f: \|f\|_{\tilde{H}^s(\mathcal{O})}^2 = \sum_{j\geq1}\lambda_j^s\langle f,e_j \rangle_{L^2}^2 < \infty \right\}, \quad s\in\R; $$
			for $f\in\tilde{H}^s(\mathcal{O})$, we have the representation
			$$ f = \sum_{j\geq 1} \langle f,e_j \rangle_{L^2} e_j $$
			which converges in $\tilde{H}^s$, as well as the duality relation $\tilde{H}^s(\mathcal{O}) = \left(\tilde{H}^{-s}(\mathcal{O})\right)^*$. Moreover, we have the embeddings
			\begin{equation}\label{Eq: spectral-normal Sobolev inclusions}
				\tilde{H}^s(\mathcal{O}) \subset H^s_0(\mathcal{O}), \, s\in\mathbb{N}, \quad \tilde{H}^{-1}(\mathcal{O})\subset \left(H^1(\mathcal{O})\right)^*,
			\end{equation}
			with equivalent norms on $\tilde{H}^s$; see (6.21) in \cite{nicklBayesianNonlinearStatistical2023}.
			
			For any $l\geq 1$, define the projection
			$$ K_l(\theta) := \sum_{j=1}^l \langle \theta,e_j \rangle_{L^2} e_j. $$
			By a standard argument, we have that whenever $\theta\in\tilde{H}^{\alpha}(\mathcal{O})$,
			\begin{equation}\label{Eq: projection approximation}
				\left\| K_l(\theta) - \theta \right\|_{\tilde{H}^{-1}(\mathcal{O})} \lesssim l^{-\frac{\alpha+1}{d}}\|\theta\|_{\tilde{H}^{\alpha}(\mathcal{O})}.
			\end{equation}
			Also, $K_l(\theta)\in\tilde{H}^{\tau}(\mathcal{O})\subset H^{\tau}(\mathcal{O})$ and so $\chi K_l(\theta)\in\mathbb{H}^{\tau,\rho}$. Then
			$$ \|\chi K_l(\theta)\|_{\mathbb{H}^{\tau,\rho}} = \rho_n \| K_l(\theta)\|_{H^{\tau}(\mathcal{O})}\simeq \rho_n \| K_l(\theta)\|_{\tilde{H}^{\tau}(\mathcal{O})} \leq \rho_n\lambda_l^{\frac{(\tau-\alpha)}{2}}\|\theta\|_{\tilde{H}^{\alpha}(\mathcal{O})} .$$
			For $\theta$ compactly supported inside $\mathcal{O}$, by an extension argument (see p140, \cite{nicklBayesianNonlinearStatistical2023}) we have that $\|\theta\|_{\tilde{H}^{\alpha}(\mathcal{O})}\lesssim \|\theta\|_{H^{\alpha}(\mathcal{O})}$, and so this ultimately yields (using also the Weyl asymptotics) that
			\begin{equation}\label{Eq: projection norm bound}
				\|\chi K_l(\theta)\|_{\mathbb{H}^{\tau,\rho}}  \lesssim \rho_n l^{\frac{\tau-\alpha}{d}}\|\theta\|_{H^{\alpha}(\mathcal{O})}.
			\end{equation}
			Finally, for any function $f$ supported on $[-1,1]^d$ we have that $\|\chi f\|_{(H^1)^*} \lesssim \|f\|_{(H^1)^*}$, by considering the duality formula for the norm and using the multiplicative inequality \eqref{Eq: multiplicative inequality} below. Since $\theta^*$ is supported on $[-1,1]^d$, we thus have that 
			$$\|\chi K_l(\theta^*) - \theta^*\|_{\tilde{H}^{-1}(\mathcal{O})} = \|\chi(K_l(\theta^*) - \theta^*)\|_{\tilde{H}^{-1}(\mathcal{O})} \lesssim \|K_l(\theta^*) - \theta^*\|_{\tilde{H}^{-1}(\mathcal{O})};$$
			combining this with \eqref{Eq: spectral-normal Sobolev inclusions} and \eqref{Eq: projection approximation}, we see that for any $l\gtrsim r_n^{-\frac{d}{\alpha+1}}$, we have that
			$$\|\chi K_l(\theta^*) - \theta^*\|_{(H^1)^*} \lesssim r_n . $$
			Choose the minimal such $l$, so that $l\simeq r_n^{-\frac{d}{\alpha+1}}$. Then by \eqref{Eq: projection norm bound},
			$$  \|\chi K_l(\theta^*)\|_{\mathbb{H}^{\tau,\rho}}  \lesssim \rho_n r_n^{-\frac{\tau-\alpha}{\alpha+1}}, $$
			and, by choosing $h=\chi K_l(\theta^*)$, the approximation term in \eqref{Eq: small ball concentration function} is therefore bounded above by a constant multiple of $\rho_n^2 r_n^{-\frac{2(\tau-\alpha)}{\alpha+1}}$. Thus to achieve \eqref{Eq: small ball concentration function} we must choose $\rho_n,r_n$ such that
			\begin{equation}\label{Eq: small ball rate ineq, oversmooth prior}
				\rho_n^2(r_n)^{-{\frac{2(\tau-\alpha)}{\alpha+1}}} \lesssim nr_n^2,\quad (\rho_nr_n)^{-\frac{2d}{2\tau+ 2 -d}} \lesssim nr_n^2. 
			\end{equation}
			The best choice of $\rho_n$ should balance the two left-hand sides, and is
			$$ \rho_n \simeq r_n^{\frac{2(\tau-\alpha) - d}{2(\alpha+1)}}. $$
			However, we stipulated that $\rho_n\to\infty$; if $\tau<\alpha+d/2$ then this choice of $\rho_n$ is valid. Otherwise, when $\tau\geq\alpha+d/2$ we pick any slowly increasing $\rho_n\to\infty$.
			In the former case, one can solve the previous display to see that the best choice of $r_n$ is (a multiple of) $n^{-\frac{\alpha+1}{2\alpha+2+d}}$, while in the latter the approximation term dominates and so the best choice of $r_n$ is $\rho_n^{\frac{\alpha+1}{\tau+1}}n^{-\frac{\alpha+1}{2\tau + 2}}$. This concludes the proof.
		\end{proof}
	\end{lemma}
	
	We now continue with the proof of Proposition \ref{Prop: specific crlbs}. With the `small ball rates' $r_n$ in hand, it remains to find a sequence $\zeta_n\to0$ such that 
	\begin{equation}\label{Eq: crlb condition}
		\varphi^*(\zeta_n)\gtrsim nr_n^2.
	\end{equation} 
	Observe that the two terms which comprise $\varphi^*$ in \eqref{Eq: concentration function definition} are both nonnegative, and so to lower bound $\varphi^*$ it suffices to lower bound either of these two terms. For the moment, let us just consider the choice(s) of $\rho_n$ given in Lemma \ref{Lemma: rescaled GP prior small ball rates}; we will later establish that these rescaling rates are optimal.
	
	\textbf{Case 1: $\tau\leq\alpha$.} In this case, the concentration term dominates in $\varphi^*$; for the choice of $r_n,\rho_n$ in Lemma \ref{Lemma: rescaled GP prior small ball rates} (i), it suffices to choose $\zeta_n\to0$ such that
	\begin{equation} \label{Eq: undersmooth prior condition}
		-\log \tilde{\Pi}^{\tau,\rho}\left(\|\theta\|_{L^2}\leq \zeta_n\right) \gtrsim n^{\frac{d}{2\tau+2+d}}. 
	\end{equation}
	We will lower bound the left-hand side using a metric entropy argument.
	For any smooth domain $\mathcal{Y}\subset[-1,1]^d\subset\mathcal{O}$, using the fact that for $h\in H^{\tau}_0(\mathcal{Y})$ we may extend $h$ by zero to all of $\mathcal{O}$ and then $h=\chi h$, we have the chain of inclusions (by identifying functions on $\mathcal{O}$ with their restriction to $\mathcal{Y}$)
	$$ \mathbb{H}^{\tau,\rho} = \mathbb{H}^{\tau} \supset H^{\tau}_0(\mathcal{Y}) \supset \tilde{H}^{\tau}(\mathcal{Y}); $$
	only the first embedding constant depends on $\rho$ (see \eqref{Eq: spectral-normal Sobolev inclusions}). Thus for some constant $k>0$,
	$$ B_{\mathbb{H}^{\tau,\rho}}(1) = B_{\mathbb{H}^{\tau}}(\rho_n^{-1}) \supset B_{\tilde{H}^{\tau}(\mathcal{Y})}(k\rho_n^{-1}). $$
	Also, it is clear that $\|\cdot\|_{L^2(\mathcal{Y})} \leq \|\cdot\|_{L^2(\mathcal{O})}$ on $L^2(\mathcal{O})$. So to lower bound $\log N\left(B_{\mathbb{H}^{\tau,\rho}}(1), \|\cdot\|_{L^2(\mathcal{O})}, \zeta_n \right)$, it suffices to lower bound
	$$ \log N\left(B_{\tilde{H}^{\tau}(\mathcal{Y})}(k\rho_n^{-1}), \|\cdot\|_{L^2(\mathcal{Y})},\zeta_n\right) =  \log N\left(B_{\tilde{H}^{\tau}(\mathcal{Y})}(k), \|\cdot\|_{L^2(\mathcal{Y})},\rho_n\zeta_n\right) .$$

	By Remark 6.1.2 in \cite{nicklBayesianNonlinearStatistical2023} (see also Chapter 3 in \cite{edmundsFunctionSpacesEntropy1996}), we have the lower bound
	$$  \log N\left(B_{\tilde{H}^{\tau}(\mathcal{Y})}(k), \|\cdot\|_{L^2},\rho_n\zeta_n\right) \geq k' (\rho_n \zeta_n)^{-\frac{d}{\tau}}, $$
	where $k'$ depends on $k,d,\tau$. By \cite[Theorem 1.1]{liApproximationMetricEntropy1999}, this implies that
	$$ -\log \tilde{\Pi}^{\tau,\rho}\left(\|\theta\|_{L^2}\leq \zeta_n\right) \gtrsim \left(\rho_n \zeta_n\right)^{-\frac{2d}{2\tau-d}}; $$
	thus to obtain \eqref{Eq: undersmooth prior condition} it suffices for $\zeta_n$ to satisfy
	$$ \left(\rho_n \zeta_n \right)^{-\frac{2d}{2\tau-d}} \gtrsim n^{\frac{d}{2\tau+2+d}}, $$
	where $\rho_n=n^{\frac{d}{4\tau+4+2d}}$. The slowest such rate $\zeta_n$ is
	$$ \zeta_n \simeq n^{-\frac{\tau}{2\tau+2+d}}, $$
	as required.
	
	\textbf{Case 2: $\tau>\alpha$.} In this case, the approximation term dominates in $\varphi^*$. Thus for any $\varepsilon>0$, we use the bound
	$$ \varphi^*(\varepsilon) \geq \frac{1}{2}\inf_{h\in \mathbb{H}^{\tau,\rho},\|h - \theta^*\|_{L^2}\leq\varepsilon} \|h\|_{\mathbb{H}^{\tau,\rho}}^2. $$
	Observe that for any function $h \in \mathbb{H}^{\tau}$, identifying functions with their restriction to $[-1,1]^d$ we have that
	$$ \|h\|_{\mathbb{H}^{\tau}} = \|\chi g\|_{H^{\tau}(\mathcal{O})} \geq \|\chi g\|_{H^{\tau}_d} = \|h\|_{H^{\tau}_d}, $$
	where $g\in H^{\tau}(\mathcal{O})$ is such that $\|h\|_{\mathbb{H}^{\tau}} = \|\chi g\|_{H^{\tau}(\mathcal{O})}$. Thus
	$$ \|h\|_{\mathbb{H}^{\tau,\rho}} = \rho_n\|h\|_{\mathbb{H}^{\tau}} \geq \rho_n\|h\|_{H^{\tau}_d}. $$
	The value of this final quantity depends only on the values of $h$ over $[-1,1]^d$. We therefore do not increase the value of the previous infimum by replacing $\|\cdot\|_{\mathbb{H}^{\tau,\rho}}$ with $\rho_n\|\cdot\|_{H^{\tau}_d}$ and considering all $h\in H^{\tau}_d$ such that $\|h-\theta^*\|_{L^2}\leq\varepsilon$. Thus to apply Proposition \ref{Prop: general contraction rate LB},  it suffices to find $\zeta_n\to0$ such that
	$$ \frac{1}{2}\rho_n^2 \inf_{h\in H^{\tau}_d, \|h - \theta^*\|_{L^2}\leq\zeta_n} \|h\|^2_{H^{\tau}_d} \gtrsim nr_n^2. $$
	Fix $S>\tau$ and let $(\phi,\psi_{lk})_{l\geq0,0\leq k<2^{ld}}$ be a $S$-regular boundary-corrected wavelet basis of $L^2([-1,1]^d)$ (see Section 4.3.5 of \cite{gineMathematicalFoundationsInfinitedimensional2016} for details). By the wavelet characterisation of Sobolev spaces, for any $h\in H^{\tau}_d$, 
	$$\|h\|_{H^{\tau}_d}^2 \simeq |\langle h,\phi\rangle|^2 + \sum_{l\geq0}2^{2l\tau}\sum_{k=0}^{2^{ld}-1} \langle h,\psi_{lk} \rangle^2.$$
	Thus using the inequality $(x-y)^2\geq\frac{1}{2}x^2-y^2$ which holds for all $x,y\in\R$, we have for any $h\in H^{\tau}([-1,1]^d)$ and any $\theta^*\in H^{\alpha}$ that
	\begin{align}
		\|h\|_{H^{\tau}_d}^2 &\geq \sum_{l\geq0}2^{2l\tau}\sum_{k=0}^{2^{ld}-1} \langle h,\psi_{lk}\rangle^2 \nonumber \\
		&\geq \sum_{l=0}^j2^{2l\tau}\sum_{k=0}^{2^{ld}-1}\left[ \langle \theta^*,\psi_{lk}\rangle - \langle \theta^*-h,\psi_{lk}\rangle \right]^2\nonumber \\
		&\geq \frac{1}{2}\sum_{l=0}^j2^{2l\tau}\sum_{k=0}^{2^{ld}-1}\langle \theta^*,\psi_{lk}\rangle^2 - \sum_{l=0}^j2^{2l\tau}\sum_{k=0}^{2^{ld}-1}\langle \theta^*-h,\psi_{lk}\rangle^2\nonumber \\
		&\geq  \frac{1}{2}\sum_{l=0}^j2^{2l\tau}\sum_{k=0}^{2^{ld}-1}\langle \theta^*,\psi_{lk}\rangle^2 - 2^{2j\tau}\|\theta^*-h\|^2_{L^2} \label{Eq: Sobolev approx lower bound}
	\end{align}
	for any truncation point $j\geq0$. We now choose a particular $\theta^*$ of generalised additive model form. By Lemma 2 in \cite{schmidt-hieberNonparametricRegressionUsing2020} (which, as remarked after Theorem 4 in that paper, holds for all $\alpha\leq S$ where $S$ is the regularity of the wavelet basis), for any $j\geq0$ there exists $F^*_j\in H^{\alpha}_0([-d,d])$ such that $\|F^*\|_{H^{\alpha}([-d,d])}\leq K$ and for $\theta^*(x)=F^*_j(x_1+\ldots+x_d)$,
	\begin{equation}\label{Eq: hard to approximate GAM}
		 |\langle \theta^*,\psi_{jk}\rangle| = cK2^{-\frac{j}{2}(2\alpha+d)} 
	\end{equation}
	for $m2^{jd}$ values of $k$, where the constants $c$ and $m$ depend only on $d$ and the wavelet basis. Moreover, $F_j^*$ is sufficiently regular at the boundary of $[-d,d]$ (it is locally a polynomial) that it may be extended by zero outside of $[-d,d]$ to give an element of $H^{\alpha}(\R)$. By \eqref{Eq: Sobolev approx lower bound}, we see that for this $\theta^*$,
	$$ \varphi^*(\zeta_n) \gtrsim \rho_n^2 2^{2j\tau}(2^{-2j\alpha} - \zeta_n^2).$$
	The slowest choice of $\zeta_n$ such that this remains nonnegative is $\zeta_n \simeq 2^{-j\alpha}$; it remains to select the truncation point $j$, which must be chosen to satisfy
	\begin{equation}\label{Eq: truncation point selection criterion}
		\rho_n^2 2^{2j(\tau-\alpha)} \gtrsim nr_n^2.
	\end{equation}
	
	When $\alpha<\tau<\alpha + \frac{d}{2}$, the choice of $\rho_n,r_n$ from Lemma \ref{Lemma: rescaled GP prior small ball rates} yield the inequality $2^j\gtrsim n^{\frac{1}{2\alpha+2+d}}$ and thus the slowest rate $\zeta_n$ is
	$$ \zeta_n \simeq n^{-\frac{\alpha}{2\alpha+2+d}}. $$
	
	When $\tau\geq\alpha+\frac{d}{2}$, we instead obtain the inequality $2^j\gtrsim \rho_n^{-\frac{1}{\tau+1}}n^{\frac{1}{2\tau+2}}$ which gives
	$$ \zeta_n \simeq \rho_n^{\frac{\alpha}{\tau+1}}n^{-\frac{\alpha}{2\tau+2}}. $$
	
	Finally, we must argue that it is sufficient to consider the choice of $\rho_n$ prescribed by Lemma \ref{Lemma: rescaled GP prior small ball rates}, that is, other choices of $\rho_n$ (subject to the conditions in the statement of Theorem \ref{Thm: contraction rate lower bound}) yield lower bounds slower than stated in the proposition. For the remainder of the proof, we denote by $r_n^*,\rho_n^*$ the optimal small ball rate and rescaling rate from Lemma \ref{Lemma: rescaled GP prior small ball rates} and by $\zeta_n^*$ the contraction rate lower bounds given in the statement of Proposition \ref{Prop: specific crlbs}. We consider the prior $\tilde{\Pi}^{\tau,\rho}$ where we write $\rho_n = m_n\rho_n^*$ for a sequence $m_n\to0$ or $m_n\to\infty$. Write $\mathbb{H}^{\tau}$ for the RKHS of $\Pi^{\tau}$ (this is the prior of which $\tilde{\Pi}^{\tau,\rho}$ is a rescaled version), described in \eqref{Eq: rescaled prior RKHS}. We establish the best small ball rate $r_n$ achieved by $\tilde{\Pi}^{\tau,\rho}$ and then compare the resulting contraction rate lower bound $\zeta_n$ to $\zeta_n^*$, where $\zeta_n$ satisfies
	\begin{equation}\label{Eq: crlb criterion, general scaling}
		\varphi^*(\zeta_n) = \frac{1}{2}m_n^2(\rho_n^*)^2 \inf_{h\in \mathbb{H}^{\tau},\|h - \theta^*\|_{L^2}\leq\zeta_n}\|h\|_{\mathbb{H}^{\tau}}^2 + (m_n\rho_n^*\zeta_n)^{-\frac{2d}{2\tau-d}} \gtrsim nr_n^2
	\end{equation}
	for a sufficiently large constant.
	
	First we consider the case $\tau\leq\alpha$. By \eqref{Eq: small ball rate ineq, undersmooth prior}, it suffices to choose $r_n$ such that
	$$ m_n^2(\rho_n^*)^2 \lesssim nr_n^2, \quad (m_n\rho_n^*r_n)^{-\frac{2d}{2\tau+2-d}} \lesssim nr_n^2. $$
	Recall that in this case, $(\rho_n^*)^2 = n(r_n^*)^2$; thus solving each of these individually gives the bounds
	\begin{equation}\label{Eq: small ball rate ineq 2, undersmooth prior}
		r_n \gtrsim m_n r_n^*,\quad r_n\gtrsim m_n^{-\frac{d}{2(\tau+1)}} r_n^*.
	\end{equation}
	
	Since when $\tau\leq\alpha$ we assume that $\rho_n\lesssim \rho_n^*= n^{\frac{d}{4\tau+4+d}}$, we need only consider the case $m_n\to0$. Then the best possible choice of $r_n$ satisfying \eqref{Eq: small ball rate ineq 2, undersmooth prior} is $r_n\simeq m_n^{-\frac{d}{2(\tau+1)}}r_n^*$. We now wish to find $\zeta_n$ satisfying \eqref{Eq: crlb criterion, general scaling}. It suffices to ignore the first term and choose $\zeta_n$ such that
	\begin{align*}
		(m_n\rho_n^*\zeta_n)^{-\frac{2d}{2\tau-d}} &\gtrsim n\left[ m_n^{-\frac{d}{2(\tau+1)}}r_n^*\right]^2 \\
		\Leftrightarrow \zeta_n &\lesssim m_n^{-\frac{d+2}{2(\tau+1)}}\zeta_n^*;
	\end{align*}
	since $m_n\to0$, we can choose $\zeta_n$ to be slower than $\zeta_n^*$.
	
	
	Next we consider the case $\tau>\alpha$. By \eqref{Eq: small ball rate ineq, oversmooth prior}, we must now choose $r_n$ to satisfy
	$$  m_n^2(\rho_n^*)^2(r_n)^{-{\frac{2(\tau-\alpha)}{\alpha+1}}} \lesssim nr_n^2,\quad (m_n\rho_n^*r_n)^{-\frac{2d}{2\tau+ 2 -d}} \lesssim nr_n^2. $$
	Using the relationship between $\rho_n^*,r_n^*$ established in Lemma \ref{Lemma: rescaled GP prior small ball rates}, we may solve these individually to give
	\begin{equation}\label{Eq: small ball rate ineq 2, oversmooth prior}
		r_n \gtrsim m_n^{\frac{\alpha+1}{\tau+1}} r_n^*,\quad r_n\gtrsim m_n^{-\frac{d}{\tau+1}} r_n^*.
	\end{equation}
	
	Suppose that $m_n\to0$. Then we choose $r_n\simeq m_n^{-\frac{d}{\tau+1}} r_n^*$ which, analogously to above, yields
	$$ \zeta_n \lesssim m_n^{-\frac{d+1}{\tau+1}}\zeta_n^*, $$
	and since $m_n\to0$ we may choose $\zeta_n$ slower than $\zeta_n^*$.
	
	If instead $m_n\to\infty$, we choose $r_n \simeq m_n^{\frac{\alpha+1}{\tau+1}} r_n^*$ and then, arguing analogously to how we obtained \eqref{Eq: truncation point selection criterion}, we can take $\zeta_n\simeq 2^{-j\alpha}$ where $j$ must be chosen to satisfy
	$$ m_n^2(\rho_n^*)^2 2^{2j(\tau-\alpha)} \gtrsim m_n^{\frac{2(\alpha+1)}{\tau+1}}r_n^*; $$
	choosing the smallest such $j$ (regardless of whether $\tau<\alpha+d/2$ or not) leads to
	$$ \zeta_n = m_n^{\frac{\alpha}{\tau+1}}\zeta_n^*, $$
	which is slower than $\zeta_n$ since $m_n\to\infty$.
	
	This concludes the proof of Proposition \ref{Prop: specific crlbs}, and hence Theorem \ref{Thm: contraction rate lower bound}. \qed

	\begin{remark}[Upper and lower bound when $\tau\leq\alpha$, $\rho_n\gg \rho_n^*$] \label{Remark: suboptimal lower bound}
		Proposition \ref{Prop: specific crlbs} does not address the case where $\tau\leq\alpha$ and $\rho_n$ is faster than $\rho_n^*$. In this case, writing $\rho_n = m_n\rho_n^*$ for a sequence $m_n\to\infty$, the proof of Lemma \ref{Lemma: rescaled GP prior small ball rates} tells us that the small ball rate for $\tilde{\Pi}^{\tau,\rho}$ is
		$$ r_n \simeq m_n r_n^*. $$
		We see that $r_n$ satisfies the relationship $\rho_n^2 = nr_n^2$; using this fact, following the argument of Theorem 2.2.2 from \cite{nicklBayesianNonlinearStatistical2023} one deduces that $r_n$ is a contraction rate in prediction risk for $\tilde{\Pi}^{\tau,\rho}$. The theorem further implies that $r_n^{\frac{\beta-1}{\beta+1}}$ is an $L^2$-contraction rate. Both of these upper bounds are slower than the rates obtained by using the best rescaling $\rho_n^*$, which are $r_n^*$ and $(r_n^*)^{\frac{\beta-1}{\beta+1}}$ respectively. 
		
		The conventional wisdom in Bayesian nonparametrics is that such a small ball rate is sharp for a rescaled Gaussian prior (see, for example, \cite{vandervaartRatesContractionPosterior2008} and Section 11.5 in \cite{ghosalFundamentalsNonparametricBayesian2017}), and should lead to a matching lower bound. However, our proof technique using Proposition \ref{Prop: general contraction rate LB} only allows us to obtain the lower bound
		$$ \zeta_n \simeq m_n^{-\frac{2\tau}{d}}\zeta_n^*. $$
		We believe that this is an artefact of our proof, and that there is no accelerated rate for undersmooth priors with fast rescaling.
	\end{remark}

	
	\section{PDE Results for Inverse Problems}\label{Appendix: PDE facts}
	
	In this appendix, we give definitions of the function spaces used in the paper, and confirm Conditions \ref{Cond: G unif boundedness}, \ref{Cond: forward Lipschitz} and \ref{Cond: stability estimate} for suitable parameter choices when $\G$ is the forward map defined by \eqref{Eq: forward map}, where $f_{\theta}$ is given by \eqref{Eq: link function} in Darcy's problem or by \eqref{Eq: link function 2} in the Schr\"{o}dinger potential problem.
	
	\subsection{Function Spaces}
	
	In this section, $\mathcal{X}$ stands for either a smooth domain $\mathcal{O}\subset\Rd$ (that is, a non-empty, open, bounded set with smooth boundary $\partial\mathcal{O}$) or the unit cube $[-1,1]^d$. For $x\in\mathcal{X}$, let $|x|$ denote the Euclidean norm of $x$.
	
	Given $\beta\in\mathbb{N}$, we let $C^{\beta}(\mathcal{X})$ denote the space of $\beta$-times differentiable functions $\mathcal{X}\to\R$ with uniformly continuous derivatives, endowed with the norm
	$$ \|f\|_{C^{\beta}} = \sum_{|i|\leq\beta} \sup_{x\in\mathcal{X}}|D^if(x)|, $$
	where for any multi-index $i\in\mathbb{Z}_{\geq0}^d$, $D^i$ denotes the $i^{th}$ partial differential operator.
	Next, for any $\gamma\in(0,1)$ we define the H\"{o}lder semi-norm
	$$ |f|_{\gamma} = \sup_{x,y\in\mathcal{X},x\neq y}\frac{|f(x)-f(y)|}{|x-y|^{\gamma}}. $$
	For general $\beta>0$, let $\lfloor\beta\rfloor$ be the largest integer less than or equal to $\beta$;  define the H\"{o}lder norm
	$$ \|f\|_{C^{\beta}} = \|f\|_{C^{\lfloor\beta\rfloor}} + \sum_{|i|=\lfloor\beta\rfloor} |D^if|_{\beta-\lfloor\beta\rfloor} $$
	with the convention $|\cdot|_0 \equiv 0$, and the H\"{o}lder space
	$$ C^{\beta}(\mathcal{X}) = \left\{ f\in C(\mathcal{X}) : \|f\|_{C^{\beta}}<\infty \right\} $$
	normed by $\|\cdot\|_{C^{\beta}}$. Let $C^{\infty}(\mathcal{X}) = \cap_{\beta\geq0}C^{\beta}(\mathcal{X})$ denote the space of smooth functions on $\mathcal{X}$.
	
	We denote by $L^2(\mathcal{X})$ the Hilbert space of square-integrable functions $\mathcal{X}\to\R$, endowed with its usual inner product $\langle\cdot,\cdot\rangle_{L^2}$. For integer $\alpha\geq0$, we define the $\alpha$-smooth Sobolev space on $\mathcal{X}$ as
	$$ H^{\alpha}(\mathcal{X}) = \left\{ f\in L^2(\mathcal{X}) : \forall |i|\leq \alpha, \,\exists D^if \in L^2(\mathcal{X}) \right\}. $$
	This is a separable Hilbert space when endowed by the inner product
	$$ \langle f,g \rangle_{H^{\alpha}(\mathcal{X})} = \sum_{|i|\leq\alpha} \langle D^if,D^ig \rangle_{L^2}; $$
	write $\|\cdot\|_{H^{\alpha}(\mathcal{X})}$ for the associated Hilbert norm. For general $\alpha\geq0$, we define $H^{\alpha}(\mathcal{X})$ by interpolation (see, for example, \cite{lionsNonHomogeneousBoundaryValue1972}). Given $\alpha>\frac{d}{2}$, we have the Sobolev embedding $H^{\alpha}(\mathcal{X})\subset C^{\alpha-\frac{d}{2}}(\mathcal{X})$. We also recall the multiplicative inequality
	\begin{equation}\label{Eq: multiplicative inequality}
		\|fg\|_{H^{\alpha}} \lesssim \|f\|_{C^{\alpha}}\|g\|_{H^{\alpha}}, \quad\alpha\geq0
	\end{equation}
	which holds for all $f,g$ in the appropriate spaces (see Theorem 2.8.2 and p143 of \cite{triebelTheoryFunctionSpaces1983}).
	
	\subsection{Regularity Conditions on $\G$}
	
	We may now confirm the requisite conditions on $\G$ for the two specific inverse problems studied in this paper. The following draws heavily on Section 5 of \cite{nicklConvergenceRatesPenalized2020}, and we refer the interested reader to this reference for a more detailed exposition of the arguments presented below. We also note Section 2.1 of \cite{nicklBayesianNonlinearStatistical2023}, which introduces and checks these conditions using Sobolev spaces $H^{\beta}$ in place of the H\"{o}lder spaces $C^{\beta}$ as regularisation spaces. As discussed previously (see Remark \ref{Remark: can't use Sobolev}), Sobolev norms are not compatible with the compositional structures considered in this paper, but the PDE arguments are largely the same.
	
	The following result on link functions is standard, and we state it here to obtain explicit constants.
	\begin{lemma}\label{Lemma: link function properties}
		Consider the link function $\theta\mapsto f_{\theta}$ defined in \eqref{Eq: link function} or \eqref{Eq: link function 2}. Then given $M>0$, for $\theta_1,\theta_2\in C(\mathcal{O})$ with $\|\theta_i\|_{\infty}\leq M$, we have that
		$$ e^{-M}\|\theta_1 - \theta_2\|_{L^2} \leq \|f_{\theta_1} - f_{\theta_2}\|_{L^2} \leq e^{M}\|\theta_1-\theta_2\|_{L^2} $$
		and
		$$ e^{-M}\|\theta_1 - \theta_2\|_{\infty} \leq \|f_{\theta_1} - f_{\theta_2}\|_{\infty} \leq e^{M}\|\theta_1-\theta_2\|_{\infty} .$$
		Moreover, for any integer $\beta>0$, we have that if $\theta\in B_{C^{\beta}}(M)$ then
		\begin{itemize}
			\item for $f$ defined by \eqref{Eq: link function}, $\|f_{\theta}\|_{C^{\beta}} \leq M^{\beta}e^{M}+ K_{\min};$
			\item for $f$ defined by \eqref{Eq: link function 2}, $\|f_{\theta}\|_{C^{\beta}} \leq M^{\beta}e^{M}.$
		\end{itemize}
%
	\end{lemma}
	The lemma means that we may check Conditions 1-3 for $f_{\theta}$ in place of $\theta$, which we now do below.
	
	\subsubsection{Darcy's Problem}
	
	For $f\in C^1(\bar{\mathcal{O}})$ with $f\geq K_{\min}>0$, define the differential operator
	$$ L_f:H^2(\mathcal{O})\to L^2(\mathcal{O}), \quad L_f[u] = \nabla\cdot(f\nabla u). $$
	Standard elliptic PDE theory (e.g. Chapter 8 of \cite{gilbargEllipticPartialDifferential2001}) tells us that there exists a bounded linear inverse operator $V_f:L^2(\mathcal{O})\to H^2_0(\mathcal{O})$, such that for any $\psi\in L^2(\mathcal{O})$, $V_f[\psi]$ weakly solves the Dirichlet problem
	\begin{equation}
		\begin{aligned}
			L_f[u] &= \psi \quad \text{ on }\mathcal{O}, \\
			u &= 0 \quad \text{ on }\partial\mathcal{O}.
		\end{aligned}
	\end{equation}
	Recall that $\G(\theta) = G(f_{\theta}) = V_{f_{\theta}}[g]$, where $g$ is the known, smooth source term. Then Lemma 20 of \cite{nicklConvergenceRatesPenalized2020} (which really only requires $f\in C^1$) immediately yields that for any $\theta\in\Theta\subset C^1(\mathcal{O})$, 
	\begin{equation}\label{Eq: Darcy uniform boundedness}
		\|\G(\theta)\|_{\infty} \leq C\|g\|_{\infty}
	\end{equation}
	where $C>0$ depends only on $\mathcal{O}$ and $K_{\min}$. This establishes Condition \ref{Cond: G unif boundedness}.
	
	Next, we check the Lipschitz condition \eqref{Eq: forward Lipschitz condition}. Fix $\beta\geq1$ and assume that $\theta_1,\theta_2 \in C^{\beta}(\mathcal{O})$, with $\|\theta_i\|_{C^{\beta}}\leq M$ for $i=1,2$. We follow the proof of Theorem 9 in \cite{nicklConvergenceRatesPenalized2020}: observe that $u_{f_{\theta_1}} - u_{f_{\theta_2}} = 0$ on $\partial\mathcal{O}$ and on $\mathcal{O}$,
	$$ L_{f_{\theta_1}}[u_{f_{\theta_1}} - u_{f_{\theta_2}}] = g - g + \left(L_{f_{\theta_1}}-L_{f_{\theta_2}}\right)u_{f_{\theta_2}} = \nabla\cdot\left([f_{\theta_1}-f_{\theta_2}]\nabla u_{f_{\theta_2}}\right).$$
	This right-hand side is clearly in $L^2(\mathcal{O})$ (indeed, it is continuous) so by Lemma 21 in \cite{nicklConvergenceRatesPenalized2020}, we have for some constant $C = C(\mathcal{O},K_{\min})$ that
	$$ \|\G(\theta_1) - \G(\theta_2)\|_{L^2} \leq C\left(1 + \|f_{\theta_1}\|_{C^1}\right)\left\|\nabla\cdot\left([f_{\theta_1}-f_{\theta_2}]\nabla u_{f_{\theta_2}}\right)\right\|_{(H^2_0)^*}. $$
	As $\|f_{\theta_1}\|_{C^1}$ is bounded by Lemma \ref{Lemma: link function properties} and the fact that $\|\theta_1\|_{C^1}\leq M$, it suffices to bound the final norm suitably. Observe that by using the divergence theorem twice we obtain
	\begin{align*}
		\left\|\nabla\cdot\left([f_{\theta_1}-f_{\theta_2}]\nabla u_{f_{\theta_2}}\right)\right\|_{(H^2_0)^*}
		&= \sup_{\varphi\in H^2_0,\|\varphi\|_{H^2}\leq 1}\left|\int_{\mathcal{O}} \varphi \nabla\cdot\left([f_{\theta_1}-f_{\theta_2}]\nabla u_{f_{\theta_2}}\right)\right| \\
		&= \sup_{\varphi\in H^2_0,\|\varphi\|_{H^2}\leq 1}\left|\int_{\mathcal{O}} [f_{\theta_1}-f_{\theta_2}]\nabla\varphi\cdot \nabla u_{f_{\theta_2}}\right| \\
		&\leq \|f_{\theta_1}-f_{\theta_2}\|_{\infty}\sup_{\varphi\in H^2_0,\|\varphi\|_{H^2}\leq 1}\left|\int_{\mathcal{O}} \nabla\varphi\cdot \nabla u_{f_{\theta_2}}\right| \\
		&=  \|f_{\theta_1}-f_{\theta_2}\|_{\infty}\sup_{\varphi\in H^2_0,\|\varphi\|_{H^2}\leq 1}\left|\int_{\mathcal{O}} u_{f_{\theta_2}}\Delta\varphi \right| \\
		&\leq \|u_{f_{\theta_2}} \|_{\infty}\|f_{\theta_1}-f_{\theta_2}\|_{\infty},
	\end{align*}
	and by \eqref{Eq: Darcy uniform boundedness}, $\|u_{f_{\theta}}\|_{\infty}$ is bounded by a constant depending only on $g, K_{\min},\mathcal{O}$. This proves Condition \ref{Cond: forward Lipschitz}.
	
	It remains to show that the stability estimate \eqref{Eq: stability estimate} holds for a suitable choice of $\beta$. Let $\beta>1$, and let $\theta^*, \theta \in B_{C^{\beta}}(M)$. Note that Proposition 2.1.5 in \cite{nicklBayesianNonlinearStatistical2023} holds for $\theta\in C^{\beta}(\mathcal{O}),\beta>1$ rather than just $\theta\in H^{\beta}(\mathcal{O}),\beta>d/2 + 1$ since for $\theta\in C^{\beta}$, $\G(\theta)\in C^{\beta+1}\subset C^2$ and we can use the multiplicative inequality \eqref{Eq: multiplicative inequality} for \emph{any} positive smoothness rather than the version for Sobolev norms (which requires $\alpha>d/2$). This yields that
	\begin{equation}\label{Eq: Darcy L2-H2 stability estimate}
		\|\theta - \theta^*\|_{L^2} \leq C\|u_{f_{\theta}} - u_{f_{\theta^*}}\|_{H^2},
	\end{equation}
	where $C=C(\mathcal{O},g,K_{\min},M)>0$. By the interpolation inequality for Sobolev spaces, we have that
	$$\|u_{f_{\theta}} - u_{f_{\theta^*}}\|_{H^2} \lesssim \|u_{f_{\theta}} - u_{f_{\theta^*}}\|_{L^2}^{\frac{\beta-1}{\beta+1}}\|u_{f_{\theta}} - u_{f_{\theta^*}}\|_{H^{\beta+1}}^{\frac{2}{\beta+1}}, $$
	for a constant depending on $\beta,\mathcal{O}$ only and so \eqref{Eq: stability estimate} follows if we can bound the final Sobolev norm; it clearly suffices to bound $\|u_{f_{\theta}}\|_{H^{\beta+1}}$ and $\|u_{f_{\theta^*}}\|_{H^{\beta+1}}$. We prove this by following the method of Lemma 23 in \cite{nicklConvergenceRatesPenalized2020}. Let $f\in C^{\beta}(\mathcal{O})$. As the Laplacian $\Delta$ is a linear isomorphism $H^{\beta+1}\to H^{\beta-1}$, by rearranging the PDE \eqref{Eq: Darcy's problem} we have that for a constant depending only on $\mathcal{O}$,
	$$ \|u_f\|_{H^{\beta+1}} \lesssim \left\| f^{-1}(g-\nabla f\cdot\nabla u_{f}) \right\|_{H^{\beta-1}}. $$
	Using the multiplicative inequality \eqref{Eq: multiplicative inequality}, this is further bounded by
	$$ \|f^{-1}\|_{C^{\beta-1}}\|g-\nabla f\cdot\nabla u_{f}\|_{H^{\beta-1}}. $$
	By Lemma 29 in \cite{nicklConvergenceRatesPenalized2020} applied to $x\mapsto x^{-1}, x\in(K_{\min},\infty)$ we have for integer $\beta\geq0$ that
	$$ \|f^{-1}\|_{C^{\beta-1}} \leq C(\beta,K_{\min})(1 + \|f\|_{C^{\beta-1}}^{\beta-1}), $$
	and so again using the multiplicative inequality \eqref{Eq: multiplicative inequality} and the interpolation inequality,
	\begin{align*}
		\|u_f\|_{H^{\beta+1}} &\lesssim \left(1 + \|f\|_{C^{\beta-1}}^{\beta-1}\right)\left(\|g\|_{H^{\beta-1}} + \|f\|_{C^{\beta}}\|u_f\|_{H^{\beta}}\right) \\
		&\lesssim \left( 1 + \|f\|_{C^{\beta}}^{\beta}\right)\left( 1 + \|u_f\|_{H^{\beta+1}}^{\frac{\beta}{\beta+1}}\|u_f\|_{L^2}^{\frac{1}{\beta+1}}\right),
	\end{align*}
	for a constant depending only on $\mathcal{O}, K_{\min}, \beta$ and $g$. By \cite[Lemma 20]{nicklConvergenceRatesPenalized2020}, $\|u_f\|_{L^2}$ is bounded by a constant multiple of $\|g\|_{L^2}$. Thus rearranging the above inequality gives
	$$ \|u_f\|_{H^{\beta+1}} \lesssim 1 + \|f\|_{C^{\beta}}^{\beta(\beta+1)}. $$
	Since $\theta,\theta^*\in B_{C^{\beta}}(M)$ implies that $f_{\theta},f_{\theta^*}\in B_{C^{\beta}}(M')$ for some $M'>0$ by Lemma \ref{Lemma: link function properties}, this establishes Condition \ref{Cond: stability estimate} for any integer $\beta>1$ with $L'$ the constant from the previous inequality (depending only on $\mathcal{O},K_{\min},\beta,g$), $\xi = \beta(\beta+1)$ and
	$$ \zeta = \frac{\beta-1}{\beta + 1}. $$

	\subsubsection{Schr\"{o}dinger Problem}
	
	For $f\in C(\bar{\mathcal{O}}), f\geq0$, define the differential operator
	$$ L_f: H^2(\mathcal{O})\to L^2(\mathcal{O}), \quad L_f[u] = \frac{1}{2}\Delta u - fu. $$
	Then as in the previous case, standard elliptic PDE theory implies the existence of a bounded linear inverse operator $V_f$ such that for $\psi\in L^2(\mathcal{O})$, $V_f[\psi]$ solves the inhomogeneous equation
	\begin{equation}
		\begin{aligned}
			L_f[u] &= \psi \quad \text{ on }\mathcal{O}, \\
			u &= 0 \quad \text{ on }\partial\mathcal{O}.
		\end{aligned}
	\end{equation}
	As before, $\G(\theta) = G(f_{\theta}) = V_{f_{\theta}}[h]$.
	
	The Feynman-Kac formula instantly verifies Condition \ref{Cond: G unif boundedness} for $\G$ with $U=\|h\|_{\infty}$: see equation (2.6) and the surrounding discussion in \cite{nicklBayesianNonlinearStatistical2023}.
	
	To check the Lipschitz condition \ref{Eq: forward Lipschitz condition}, we proceed similarly to before. Note that for any $\theta_1,\theta_2\in C(\mathcal{O})$, we have that $u_{f_{\theta_1}} - u_{f_{\theta_2}} = h-h = 0$ on $\partial\mathcal{O}$, and on $\mathcal{O}$,
	$$ L_{f_{\theta_1}}[u_{f_{\theta_1}} - u_{f_{\theta_2}}] = (f_{\theta_1} - f_{\theta_2})u_{f_{\theta_2}}. $$
	Combining this with Lemma 25 in \cite{nicklConvergenceRatesPenalized2020} then gives
	$$ \|\G(\theta_1) - \G(\theta_2)\|_{L^2} \leq C \|(f_{\theta_1} - f_{\theta_2})u_{f_{\theta_2}}\|_{L^2} \leq C\|h\|_{\infty}\|f_{\theta_1} - f_{\theta_2}\|_{\infty},$$
	where we have used the uniform boundedness condition established previously to bound $\|u_{f_{\theta_2}}\|_{\infty}$. This confirms Condition \ref{Cond: forward Lipschitz} for any choice of $\beta\geq0$.
	
	Lastly, we must show that the stability estimate \eqref{Eq: stability estimate} holds for a suitable choice of $\beta$. We follow the scheme of Lemma 28 in \cite{nicklConvergenceRatesPenalized2020}. Let $f\in C(\bar{\mathcal{O}})$. From the Feynman-Kac formula, one obtains that
	\begin{equation}\label{Eq: F-K lower bound for u}
		\inf_{x\in\mathcal{O}}\, u_f(x) \geq h_{\min}e^{-c\|f\|_{\infty}}
	\end{equation}
	for some $c>0$ depending only on $\mathcal{O}$. By rearranging the PDE \eqref{Eq: Schrodinger problem}, we have that $f = (\Delta u_f)/2u_f$ on $\mathcal{O}$. Thus, using \eqref{Eq: F-K lower bound for u} and \eqref{Eq: multiplicative inequality}, we have that for $f_1,f_2\in C(\bar{\mathcal{O}})$
	\begin{align}
		\|f_1 - f_2\|_{L^2} &= \frac{1}{2}\left\|\frac{\Delta u_{f_1}}{u_{f_1}} - \frac{\Delta u_{f_2}}{u_{f_2}}\right\|_{L^2} \nonumber \\
		&\lesssim \left\| \frac{\Delta u_{f_1} - \Delta u_{f_2}}{u_{f_1}}\right\|_{L^2} + \left\|\Delta u_{f_2}\left(\frac{1}{u_{f_1}} - \frac{1}{u_{f_2}}\right)\right\|_{L^2} \nonumber \\
		&\lesssim h_{\min}^{-1}e^{c\|f_1\|_{\infty}}\|u_{f_1} - u_{f_2}\|_{H^2} + \|u_{f_2}\|_{C^2}\left\|\frac{1}{u_{f_1}} - \frac{1}{u_{f_2}}\right\|_{L^2}. \label{Eq: Schrodinger H2-L2 stability}
	\end{align}
	Again using \eqref{Eq: F-K lower bound for u} and the mean value theorem, we have that
	$$ \left\|\frac{1}{u_{f_1}} - \frac{1}{u_{f_2}}\right\|_{L^2} \leq h_{\min}^{-2}e^{c(\|f_1\|_{\infty} + \|f_2\|_{\infty})}\|u_{f_1} - u_{f_2}\|_{L^2}. $$
	Also, the first part of Lemma 27 in \cite{nicklConvergenceRatesPenalized2020} (which only requires $f\in C(\bar{\mathcal{O}})$) yields
	$$ \|u_{f_2}\|_{C^2} \leq C(1+\|f_2\|_{\infty})\|h\|_{C^2(\partial\mathcal{O})},$$
	where $C>0$ depends on $\mathcal{O}$ only. Plugging these two bounds into \eqref{Eq: Schrodinger H2-L2 stability}, one obtains for any $f_1,f_2\in B_{C(\bar{\mathcal{O}})}(M)$ that
	\begin{equation}\label{Eq: Schrodinger H2-L2 stability 2}
		\|f_1 - f_2\|_{L^2} \leq C\|u_{f_1} - u_{f_2}\|_{H^2}
	\end{equation}
	for a constant $C>0$ depending on $M,\mathcal{O},h_{\min}$.
	
	Now assume that $f_i\in C^{\beta}(\mathcal{O})$ for some $\beta>0$. By the Sobolev interpolation inequality we have that
	$$ \|f_1 - f_2\|_{L^2} \lesssim \|u_{f_1} - u_{f_2}\|_{H^2} \lesssim \|u_{f_1}-u_{f_2}\|_{L^2}^{\frac{\beta}{\beta+2}}\|u_{f_1}-u_{f_2}\|_{H^{\beta+2}}^{\frac{2}{\beta+2}}, $$
	and so appealing to Lemma \ref{Lemma: link function properties} as before, to establish \eqref{Eq: stability estimate} it suffices to show that $\|u_{f_i}\|_{H^{\beta+2}}, i=1,2$ are bounded. The argument follows the method of the second part of \cite[Lemma 27]{nicklConvergenceRatesPenalized2020}: since $\Delta$ is an isomorphism between Sobolev spaces we have that, by rearranging the PDE and using the interpolation and multiplicative inequalities,
	\begin{align*}
		\|u_f\|_{H^{\beta+2}} &\lesssim \|fu_f\|_{H^{\beta}} + \|g\|_{C^{\beta+1}(\partial\mathcal{O})} \\
		&\lesssim \|f\|_{C^{\beta}}\|u_f\|_{H^{\beta}} + 1 \\
		&\leq 1 + \|f\|_{C^{\beta}}\|u_f\|_{L^2}^{\frac{2}{\beta+2}}\|u_f\|_{H^{\beta+2}}^{\frac{\beta}{\beta+2}}\\
		\Rightarrow \|u_f\|_{H^{\beta+2}} &\lesssim 1 + \|u_f\|_{L^2}\|f\|_{C^{\beta}}^{\frac{\beta+2}{2}} \\
		&\lesssim 1 + \|f\|_{C^{\beta}}^{\frac{\beta+2}{2}}
	\end{align*}
	for a constant depending only on $g,\mathcal{O},\beta$, where in the final line we used the uniform boundedness property. This establishes the stability estimate Condition \ref{Cond: stability estimate} for any choice of $\beta>0$ with $L'$ the constant from the previous inequality, $\xi = \beta/2 + 1$, and
	$$ \zeta = \frac{\beta}{\beta+ 2}. $$

\end{document}